\documentclass[10pt]{amsart}
\usepackage{amsfonts,amssymb}
\RequirePackage{ifpdf}
\usepackage{amsmath}
\usepackage{color}
\usepackage{pstcol,pst-node}
\usepackage{graphics}
\usepackage{epic}
\usepackage{curves}
\usepackage[matrix,arrow,frame,curve]{xy}
\CompileMatrices

\def\appendix#1{
\addtocounter{section}{1} \setcounter{equation}{0}
\renewcommand{\thesection}{\Alph{section}}
\section*{Appendix \thesection\protect\indent\quad
#1}
}
\renewcommand{\theequation}{\thesection.\arabic{equation}}

\catcode`\@=11
\def\marginnote#1{}

\newcount\hour
\newcount\minute
\newtoks\amorpm
\hour=\time\divide\hour by60 \minute=\time{\multiply\hour by60
\global\advance\minute by-\hour}
\edef\standardtime{{\ifnum\hour<12 \global\amorpm={am}%
        \else\global\amorpm={pm}\advance\hour by-12 \fi
        \ifnum\hour=0 \hour=12 \fi
        \number\hour:\ifnum\minute<10 0\fi\number\minute\the\amorpm}}
\edef\militarytime{\number\hour:\ifnum\minute<100\fi\number\minute}

\def\draftlabel#1{{\@bsphack\if@filesw {\let\thepage\relax
      \xdef\@gtempa{\write\@auxout{\string
          \newlabel{#1}{{\@currentlabel}{\thepage}}}}}\@gtempa \if@nobreak
    \ifvmode\nobreak\fi\fi\fi\@esphack} \gdef\@eqnlabel{#1}}
    \def\@eqnlabel{}
\def\@vacuum{}
\def\draftmarginnote#1{\marginpar{\raggedright\scriptsize\tt#1}}

\def\draft{
  \oddsidemargin -.5truein
  \def\@oddfoot{\footnotesize \sl preliminary draft \hfil
    \rm\thepage\hfil\sl\today\quad\militarytime}
  \let\@evenfoot\@oddfoot \overfullrule 3pt
    \let\label=\draftlabel
    \let\marginnote=\draftmarginnote
  \def\@eqnnum{(\theequation)\rlap{\kern\marginparsep\tt\@eqnlabel}%
    \global\let\@eqnlabel\@vacuum}

  }
\newcommand{\beq}{\begin{equation}}
\newcommand{\eeq}{\end{equation}}

\def\be{\begin{equation}}
\def\ee{\end{equation}}
\def\bea{\begin{eqnarray}}
\def\eea{\end{eqnarray}}
\def\<{\langle}
\def\>{\rangle}

\def\nn{\nonumber}

\def\C{{\Bbb C}}
\def\P{{\Bbb P}}

\def\one#1{#1^{\raise5pt\hbox{$\scriptstyle\!\!\!\!1$}}\,{}}
\def\two#1{#1^{\raise5pt\hbox{$\scriptstyle\!\!\!\!2$}}\,{}}
\def\onetwo#1{#1^{\raise5pt\hbox{$\scriptstyle\!\!\!\!\!{12}$}}\,{}}

\newtheorem{theorem}{Theorem}[section]
\newtheorem{lm}[theorem]{Lemma}
\newtheorem{prop}[theorem]{Proposition}
\theoremstyle{definition}
\newtheorem{df}[theorem]{Definition}
\newtheorem{example}[theorem]{Example}
\newtheorem{remark}[theorem]{Remark}

\newtheorem{conj}[theorem]{Conjecture}

\theoremstyle{remark}

\allowdisplaybreaks

\begin{document}
\title[{Quantised Painlev\'e monodromies, Sklyanin and CY algebras.}]
{Quantised Painlev\'e monodromy manifolds, Sklyanin and Calabi-Yau algebras}
\author{Leonid Chekhov, Marta Mazzocco, Vladimir Rubtsov}

\maketitle

\maketitle

%\dedicatory{
\begin{center}
{\footnotesize\textit{{to Boris Dubrovin.}}}
\end{center}
%}

\begin{abstract}
In this paper we 
study quantum del Pezzo surfaces belonging to a certain class.  In particular we introduce the generalised Sklyanin-Painlev\'e algebra and characterise its PBW/PHS/Koszul properties. This algebra contains as limiting cases the generalised Sklyanin algebra, Etingof-Ginzburg and Etingof-Oblomkov-Rains quantum del Pezzo and the quantum monodromy manifolds of the Painlev\'e equations. 
\end{abstract}

%\tableofcontents

\section{Introduction}

In recent years, studying non-commutative rings through the methods of quantum algebraic geometry has sparked enormous interest due to its applications in mirror symmetry. The work by Gross-Hacking and Keel \cite{GrossHK} associates to Looijenga pairs on the A-side, i.e. pairs $(Y,D)$ where $Y$ is a smooth projective surface and $D$ is an anti-canonical cycle of rational curves, a mirror family  on the B-side constructed as the spectrum of an explicit algebra structure on a vector space. The elements of the basis of global sections uniformise such a spectrum and are called theta functions.

Interestingly, the A-side is equipped with a symplectic structure,  and it is quantised by geometric quantisation  {within the SYZ formalism \cite{SYZ}}, while the B side is naturally quantised by deformation quantisation.

In this paper we study a certain class of del Pezzo surfaces that can be put on either side of the mirror construction, or in other words, whose geometric and deformation quantisation coincide. In particular, we study the quantisation of a family of Poisson manifolds defined by the zero locus $\mathcal M_\phi$ of a degree $d$ polynomial $\phi\in \mathbb C[x_1,x_2,x_3]$ of the form
\begin{equation}\label{eq:pol-gen}
\phi(x_1,x_2,x_3)=x_1 x_2 x_3 + \phi_1(x_1)+ \phi_2(x_2)+ \phi_3(x_3)
\end{equation}
where $\phi_i(x_i)$ for $i=1,2,3$ is a polynomial of degree $\leq d$ in the variable $x_i$ only.

From an algebro-geometric point of view (under certain conditions on the degrees of each polynomial $\phi_i$, $i=1,2,3$) the projective completion $\overline{\mathcal M}_\phi$ in the weighted projective space  $\mathbb W \mathbb P^3$ of the affine surface $\mathcal M_\phi\subset\mathbb C^3$ is a (possibly degenerate)  del Pezzo surface.  In other words, the pair $(\overline{\mathcal M}_\phi,D_\infty)$, where $D_\infty$ is the divisor at infinity, is a Loojenga pair and $\mathcal M_\phi=\overline{\mathcal M}_\phi\setminus D_\infty$.
At the same time, each affine del Pezzo surface can be considered as  {an affine cone} ${\rm Spec}\left(\mathbb C[x_1,x_2,x_3]/\langle\phi\rangle\right)$ {over the projective curve $\mathbb P{\mathcal M}_\phi\subset \mathbb P^2$ corresponding to the divisor $D_\infty$,
with the trivial line bundle $L=\mathcal M_\phi\setminus\{0\}$ whose section ring is} $\oplus_{k\geq 0} H^0(\mathbb P{\mathcal M}_\phi, L^{\otimes k})$.  {Here, we call the curve  $\mathbb P{\mathcal M}_\phi$ {\it the projectivisation of  $\mathcal M_\phi$.}}\/ 

A quantisation of a del Pezzo surface of this type appeared in the work of Oblomkov \cite{Obl} as the spherical sub-algebra of the  $\check{C}C_1$ double affine Hecke algebra (DAHA). Then 
Etingof, Oblomkov and Rains proposed a notion of generalised DAHA for every simply laced affine Dynkin diagram and showed that their spherical sub-algebras quantise  the coordinate rings of affine surfaces obtained by removing a nodal $\mathbb P^1$ from a weighted  projective del Pezzo surface of degrees $3$, $2$ and $1$ respectively for $E_6^{(1)}$, $E_7^{(1)}$ and $E_8^{(1)}$ or by removing a triangle  from a  projective del Pezzo surface of degree $3$ in the case $D_4^{(1)}$.
In the same paper, the authors defined a holomorphic (but not algebraic) map from the mini-versal deformation of the corresponding Kleinian singularity $SL(2,\mathbb C)/\Gamma$ (where $\Gamma\in SL(2,\mathbb C)$ is the finite subgroup corresponding to the Dynkin diagram  $D_4$, $E_6$, $E_7$ and $E_8$ respectively via the McKay correspondence) to the family of surfaces $\mathcal M_\phi$ where $\phi$ is in our form:
\bea\label{eq:EOR}
D_4^{(1)}&x_1 x_2 x_3 + x_1^2 + x_2^2 + x_3^2 + \eta x_1 +\sigma x_2 + \rho x_3 + \omega , \nn\\
E_6^{(1)}& x_1 x_2 x_3 + x_1^3 + x_2^3 + x_3^2 +  \eta_2 x_1^2 + \eta_1 x_1+\sigma_2 x_2^2  +\sigma_1 x_2 + \rho x_3 + \omega,\nn\\
\qquad E_7^{(1)} & x_1 x_2 x_3 + x_1^4 + x_3^2 + x_3^2 +  \eta_3 x_1^3 +\dots + \eta_1 x_1+\sigma x_2 + \rho x_3 + \omega, \\
E_8^{(1)} & x_1 x_2 x_3 + x_1^5 + x_3^2 + x_3^2 +   \eta_4 x_1^4 +\dots + \eta_1 x_1+\sigma x_2 + \rho x_3 + \omega.  \nn
\eea
Following this work, 
P. Etingof and V. Ginzburg \cite{EtGinzb} have proposed a quantum  description of del Pezzo surfaces based on the flat deformation of cubic affine cone surfaces with an isolated elliptic singularity of type $\widetilde E_6,\widetilde E_7$ and $\widetilde E_8$ in (weighted) projective planes:
\be\label{eq:EG}
\begin{split}
\widetilde E_6 &\qquad\qquad \tau x_1 x_2 x_3+ \frac{x_1^3}{3}+\frac{x_2^3}{3}+\frac{x_3^3}{3}+   \eta_2 x_1^2 + \eta_1 x_1+ \\
 &\qquad\qquad \qquad\qquad +\sigma_2 x_2^2  +\sigma_1 x_2 + \rho_2 x_3^2+\rho_1 x_3 +\omega, \\
\widetilde E_7& \qquad\qquad\tau x_1 x_2 x_3+\frac{x_1^4}{4}+\frac{x_2^4}{4}+\frac{x_3^2}{2} +   \eta_3x_1^3 +\dots + \eta_1 x_1+ \\
&\qquad \qquad \qquad\qquad+\sigma_3 x_2^3+\dots  +\sigma_1 x_2 + \rho_2 x_3^2+\rho_1 x_3 +\omega, \\
\widetilde E_8 & \qquad\qquad\tau x_1 x_2 x_3+\frac{x_1^6}{6}+\frac{x_2^3}{3}+\frac{x_3^2}{2} +   \eta_5 x_1^5 +\dots+  \eta_2 x_1^2 + \eta_1 x_1+\\
&\qquad \qquad \qquad\qquad +\sigma_2 x_2^2  +\sigma_1 x_2 + \rho_2 x_3^2+\rho_1 x_3 +\omega.
\end{split}\ee
Their result gives a family of Calabi-Yau algebras parametrised by a complex number and a triple of polynomials of specifically chosen degrees. Interestingly, as far as we know, nobody has proved a similar result for the polynomials \eqref{eq:EOR}.   

Poisson manifolds defined by the zero locus $\mathcal M_\phi$ of a degree $3$ polynomial $\phi\in \mathbb C[x_1,x_2,x_3]$ of the form \eqref{eq:pol-gen} where $\phi_i(x_i)$ for $i=1,2,3$ is a polynomial of degree $2$ appear in the theory of the Painlev\'e differential equations as monodromy manifolds  \cite{SvdP}. Indeed,  
the  Painlev\'e sixth  monodromy manifold is precisely the affine surface that appeared in Oblomkov \cite{Obl} (see also \cite{EOR}) as the spectrum of the center of the Cherednik algebra of type $\check{C}C_1$ for $q=1$.  

This result was generalised in  \cite{M2}, where seven new algebras were produced as Whittaker degenerations of the Cherednik algebra of type $\check{C}C_1$ in such a way that their spherical--sub-algebras tend in the semi-classical limit to the monodromy manifolds of the respective Painlev\'e differential equations. 

In the present paper we give a quantisation of the Painlev\'e monodromy manifolds that fits into the scheme proposed by Etingof and Ginzburg {for the elliptic singularity $\widetilde E_6$} (see Theorem \ref{th-main-P}). 
Namely, {for an appropriate quantisation $\Phi$ of $\phi$},  we define an associative algebra $A_{\Phi}$, which is a flat deformation of the coordinate ring $\C[x_1,x_2,x_3]$ or, more precisely, the quantisation of the corresponding Poisson algebra $A_{\phi}=(\C[x_1,x_2,x_3], \{\cdot,\cdot\}_{\phi})$ where 
$$\{p,q\}_{\phi}=\frac{dp\wedge dq\wedge d\phi}{dx_1\wedge dx_2\wedge dx_3}$$ 
is the  Poisson-Nambu structure \eqref{eq:nambu} on $\C^3$ for $p,q\in   \C[x_1,x_2,x_3].$

The algebra $A_{\Phi}$ has three non-commuting generators $X_i, i=1,2,3$ subject to the relations 
$$X_iX_j -{q} X_jX_i = \phi_k (X_k), \quad (i,j,k) = (1,2,3)$$
with $\phi_k\in \C[X_k]$ and ${q} \in \C^{*}.$
One can consider the following diagram where the left and right column arrows are natural surjections and the horizontal arrows denote flat deformations or quantisations
of the corresponding Poisson algebras $A_{\phi}$ and $A_{\phi}/\langle\phi\rangle:$
\beq\label{diag1}
\xymatrix{ &A_{\phi}\ar[d]\ar@{~>}[r]^{\txt{fl. def.}}& A_{\Phi}^{{q}}\ar[d]\\
& A_{\phi}/\langle\phi\rangle\ar@{~>}[r]^{\txt{fl. def.}}&  A_{\Phi}^{{q}}/\langle \Omega-\Omega^0\rangle.\\ }
\eeq
Following the idea of \cite{EtGinzb}, we construct the bottom-right corner algebra  as a quotient of the (family of) associative algebras $A_{\Phi}^{{q}}$ by the bilateral
ideal $\langle \Omega-\Omega^0\rangle$ generated by evaluating a central element $\Omega \in A_{\Phi}^{{q}}$ for all $\phi$ corresponding to the Painlev\'e monodromy manifolds. 
As a result, we obtain a (family of) non-commutative $3$-Calabi-Yau algebras  that we denote by $\mathcal U\mathcal Z$ and their non-commutative $2$-dimensional quotients  as a quantum del Pezzo surfaces.

More precisely we give the following:

\begin{df}\label{df:UP-formal}
Given any scalars $\epsilon_1,\epsilon_2,\epsilon_3$, and $q$, $q^m\neq 1$ for any integer $m$, the {\it universal Painlev\'e algebra}\/  $\mathcal U\mathcal P$ is the non-commutative algebra with generators $X_1,X_2,X_3,\Omega_1,\Omega_2 ,\Omega_3$ defined by the relations:
\bea\label{eq:q-comm}
&&
q^{-1/2}X_1 X_2 -q^{1/2}X_2 X_1 - (q^{-1}-q) \epsilon_{3} X_3 +(q^{-1/2}-q^{1/2})\Omega_{3} =0, \nn\\
&&
q^{-1/2}X_2 X_3 -q^{1/2}X_3 X_2 - (q^{-1}-q) \epsilon_{1} X_1 + (q^{-1/2}-q^{1/2})\Omega_{1} =0,\\
&&
q^{-1/2}X_3 X_1 -q^{1/2}X_1 X_3 - (q^{-1}-q) \epsilon_{2} X_2 + (q^{-1/2}-q^{1/2})\Omega_{2} =0,\nn\\
&&
[\Omega_{i}  ,\cdot]=0, \qquad i=1,2,3.\nn
\eea
\end{df}

\begin{remark}
The name {\it universal}\/ has been chosen  because in the case 
$\epsilon_1=\epsilon_2=\epsilon_3=1$, this algebra corresponds to the Universal Askey-Wilson algebra \cite{TerU}.
\end{remark}

\begin{df}\label{df:UZ-formal} For any choice of three scalars $\Omega_i^0$, $i=1,2,3$, the {\it confluent Zhedanov algebra}\/  $\mathcal U\mathcal Z$ is the quotient 
$$
\mathcal U\mathcal P\slash \langle\Omega_1-\Omega_1^0,\Omega_2-\Omega_2^0,\Omega_3-\Omega_3^0\rangle.
$$ \end{df}

\begin{remark}
The name {\it confluent Zhedanov}\/ has been chosen   because  for different choices of the scalars $\epsilon_1,\epsilon_2,\epsilon_3$, the algebra $\mathcal U\mathcal Z$ is coincides with the confluent Zhedanov algebras studied in \cite{M2}.
\end{remark}

\begin{theorem}\label{th-main-P} 
The confluent Zhedanov algebra $\mathcal U\mathcal Z$ satisfies the following properties:
\begin{enumerate}
\item It is a Poincar\'e-Birkhoff-Witt (PBW) type deformation of the homogeneous quadratic $\C $-algebra with three generators $X_1,X_2,X_3$ and the
relations:
\bea\label{eq:q-comm-quad}
&&
q^{-1/2}X_1  X_2 -q^{1/2}X_2  X_1 = 0, \nn\\
&&
q^{-1/2}X_2   X_3 -q^{1/2}X_3   X_2 = 0,\\
&&
q^{-1/2}X_3   X_1 -q^{1/2}X_1   X_3 = 0.\nn
\eea
\item It is a family of $3$-Calabi-Yau algebras  with potential 
\bea
\label{qEGcub}
\Phi_{\mathcal U\mathcal Z} &:=& X_1X_2X_3 - q X_2X_1X_3 + \frac{q^2-1}{2\sqrt{q}}(\epsilon_1X_1^2 +\epsilon_2X_2^2 +\epsilon_3X_3^2) +\nn \\
& & +(1-q)(\Omega_1X_1 + \Omega_2X_2 +\Omega_3X_3).
\eea
\item Its  center $Z(\mathcal U\mathcal Z)$ is generated by 
\be\label{q-cubics}
{\Omega_4:= \sqrt{q}  X_3   X_2  X_1 - q  \epsilon_1  X_1^2 - \frac{ \epsilon_2}{q} X_2^2-q\epsilon_3 X_3^2+\sqrt{q}  \Omega_1 X_1+\frac{ \Omega_2}{\sqrt{q}} X_2+\sqrt{q} \Omega_3  X_3.}
\ee 
\end{enumerate}
\end{theorem}

The proof of this theorem is obtained by the combining Propositions  \ref{prop:PBW-UP}, \ref{prop:UZ-central-EG} and \ref{prop:CY-UZ-K}.
The construction of  the quotient  $\mathcal U\mathcal Z\slash\langle \Omega_4-\Omega_4^0\rangle$ within the Etingof-Ginzburg framework is carried out in Theorem  \ref{th:EG-UZ}.

Our quantisation is compatible with the Whittaker degeneration of generalised DAHA proposed in \cite{M2} - see Theorem \ref{prop-main-pq1} here below.
In particular, we show that the Kleinian case $D_4$ arises as a limit of the elliptic singularity case $\widetilde E_6$ - all other Kleinian cases follow as special limits as well as shown in \cite{ChMR}. Inspired by this, we study a broad class of 
degenerations of Poisson algebras in terms of rational degenerations of elliptic curves.   

Moreover, we connect with the work of Gross, Hacking and Keel \cite{GrossHK}, namely for each $\phi$ in the form \eqref{eq:pol-gen} we produce a Looijenga pair  $(Y,D)$ where  $Y$ is  the smooth {weighted} projective completion  of our affine surface $\mathcal M_\phi\subset\mathbb C^3$ and $D$ is some reduced effective normal crossing anticanonical divisor on $Y$  given by the divisor at infinity $D_\infty$.  This is equipped with a symplectic structure obtained by taking the Poincar\'e residue of the global $3$-form in {weighted projective space} $\mathbb W\mathbb P^3$ along the divisor $D_\infty$. This form is symplectic on $Y\setminus D_\infty=\mathcal M_\phi$ - this gives rise to the Nambu bracket on $\mathcal M_\phi$. At the same time, {the coordinate ring of}
 each affine del Pezzo $\mathcal M_\phi$ 
can be seen as {the graded ring} $\oplus_{k\geq 0} H^0(\mathbb P{\mathcal M}_\phi, L^{\otimes k})$, where $\mathbb P{\mathcal M}_\phi$ is the projectivisation of $\mathcal M_\phi$, and $L$ is a line bundle of an appropriate degree, defined by the anticanonical divisor so that the equation $\phi=0$ can be seen as  a relation  between some analogues of $theta$-functions related to toric mirror data on log-Calabi-Yau surfaces. 

Due to the fact that the Calabi Yau algebra associated to $\widetilde E_6$ specialises to the Sklyanin algebra with three generators \eqref{eq:sklyanin-algebra}, 
we provide a unified Jacobian algebra, that we call {\it generalised Sklyanin-Painlev\'e algebra,} which for different values of the parameters specialises to the generalised Sklyanin algebra \eqref{genSkl} of Iyudu and Shkarin, or to the $\widetilde E_6$-Calabi-Yau algebra of Etingof and Ginzburg or to our algebra  $\mathcal U\mathcal Z$. {This is a first step towards the more ambitious goal of describing all GHK theta functions in the non-commutative world.}

\begin{df} For any choice of the scalars $a,b,c,\alpha,\beta,\gamma, a_1,b_1,c_1,a_2,b_2,c_2\in\mathbb C$, such that $a,b,c$ are not roots of unity, the {\it generalised Sklyanin-Painlev\'e algebra}\/ is the non-commutative algebra with generators $X_1,X_2,X_3$ defined by the relations:
\bea\label{eq:q-comm-GSP}
&&
X_2X_3 -aX_3X_2 -\alpha X_1^2+a_1 X_1+a_2=0,\nn\\
&&
 X_3X_1 -bX_1X_3 -\beta X_2^2+b_1   X_2+b_2=0,\\
 &&
 X_1X_2 -cX_2X_1 -\gamma X_3^2 +c_1   X_3+c_2=0.\nn
 \eea
 \end{df}
 
 We fully characterise for which cases the generalised Sklyanin-Painlev\'e algebra is a Calabi Yau algebra with  Poincar\'e Birkhoff Witt (PBW) or Koszul properties or with a polynomial growth Hilbert series (PHS):

\begin{theorem}
For specific choices of the parameters as follows:
\begin{enumerate}
\item  $a=b=c\neq 0$ and $(a^3,\alpha\beta\gamma)\neq (-1,1)$,
\item $(a,b,c)\neq(0,0,0)$ and either $\alpha=\beta=a-b=0$ or $\gamma=\alpha=c-a=0$ or $\beta=\gamma=b-c=0$,
\item $\alpha=\beta=\gamma=0$ and  $(a,b,c)\neq(0,0,0)$,\end{enumerate}
the generalised Sklyanin-Painlev\'e algebra is potential, PHS and Koszul. 
\end{theorem}

Finally, in Theorem \ref{compar}, we deal with the question by P. Bousseau  whether his deformation quantisation of function algebras on certain affine varieties related to Looijenga pairs, proposed in the recent paper \cite{Bouss}, 
can be compared to Etingof and Ginzburg approach.

\vskip 4mm
This paper is organised as follows. In Section \ref{se:mfds-q} provide some background on the Painlev\'e monodromy manifolds and produce their quantisation in Theorem \ref{th-main-P}. In particular we introduce the family of  non-commutative algebras $\mathcal U\mathcal Z$ as the algebra generated by $\langle X_1,X_2,X_3\rangle$ and with relations \eqref{eq:q-comm}. In Section \ref{se:PBW} we discuss the notions of PBW, PHS and Koszul  property and show in what way the algebra $\mathcal U\mathcal Z$ satisfies them. In Section \ref{se:CY-Sk}, we discuss the notions of Calabi Yau algebra, the Etingof and Ginzburg construction and the Sklyanin algebra. We introduce the generalised Sklyanin-Painlev\'e algebra (see subsection \ref{se:pain-pot}) and characterise its PBW/PHS/Koszul properties. In Section \ref{se:GHK-theta}, we discuss the affine del Pezzo surfaces $\mathcal M_\phi$ for different choices of $\phi$ and their degenerations in terms of rational degenerations of elliptic curves. In Section \ref{se:NC-QFT} we provide the quantum version of such elliptic degenerations. Finally in Section \ref{se:conclusion} we provide several tables that resume all these results and discuss some open questions.

\vskip 2mm \noindent{\bf Acknowledgements.} The authors are grateful to
Yu. Berest, R. Berger,  F. Eshmatov, P. Etingof, D. Gurevich, N. Iyudu, T. Kelly, T. Koornwinder, M. Gross,  B. Pym, V. Sokolov, P. Terwilliger, A. Zhedanov for helpful discussions. 
Our special thanks to Geoffrey Powell who carefully read our first version and made several useful remarks {and to the referee whose remarks were truly helpful.}
This research was supported by the EPSRC Research Grant $EP/P021913/1$, by the Hausdorff Institute, by ANR DIADEMS and MPIM (Bonn) and SISSA (Trieste). V.R. was partly supported by the project IPaDEGAN (H2020-MSCA-RISE-2017), Grant Number 778010, and by the Russian Foundation for Basic Research under the Grants RFBR 18-01-00461 and 16-51-53034- 716 GFEN.

\vskip 2mm
{\it Boris, your strength, energy, optimism and courage during the last months of your life are a testament to the great man you are. Goodbye dear friend, teacher.}

\section{Painlev\'e monodromy manifolds and their algebraic quantisation}\label{se:mfds-q}

The Painlev\'e differential equations are nonlinear second order ordinary differential equations of the type:   
$$
y_{tt}={\mathcal R}(t,y,y_t),
$$
where ${\mathcal R}$ is rational in $y$,  $t$ and $y_t$, such that 
the general solution $y(t;c_1,c_2)$
satisfies the following two important properties (see \cite{Pain}):
\begin{enumerate}
\item  {\it Painlev\'e property:} The solutions have no movable 
critical points, i.e. the locations of multi-valued singularities of any of the solutions are independent of the particular solution chosen. 
\item {\it Irreducibility:} For generic values of the integration constants $c_1,c_2$, the solution $y(t;c_1,c_2)$ cannot be expressed via elementary or classical transcendental functions.
\end{enumerate}

The Painlev\'e differential equations possess many beautiful properties, for example they are ``integrable", i.e. they can be written as the compatibility condition 
\begin{equation}\label{eq:iso}
\frac{\partial A}{\partial t}-\frac{\partial B}{\partial \lambda}=[B,A],
\end{equation}
between an auxiliary $2\times 2$ linear system $\frac{\partial Y}{\partial \lambda}=A(\lambda;t) Y$ and an associated deformation system, under the condition that the monodromy data of the auxiliary system are constant under deformation. 

Moreover the Painlev\'e differential equations admit symmetries under affine Weyl groups which are related to the associated B\"acklund transformations. Taking these into account, to each Painlev\'e differential equation corresponds a {\it monodromy manifold,}\/ i.e. the set of monodromy data up to global conjugation and affine Weyl group symmetries. The so-called {\it Riemann--Hilbert correspondence}\/ associates to each solution of a Painlev\'e differential equation (up to B\"acklund transformations) a point in its monodromy manifold.

Each monodromy manifold is an affine cubic surface in $\mathbb C^3$ defined by the zero locus of the corresponding polynomial in
 $\mathbb C[x_1,x_2,x_3]$ given in Table 1, where  $\omega_1,\dots,\omega_4$ are some constants (algebraically dependent in all cases except PVI) related to the parameters appearing in the corresponding Painlev\'e equation. 
 
 \begin{table}[h]
\begin{center} 
\begin{tabular}{|c||c|c|} \hline 
P-eqs & Polynomials \\ \hline 
$PVI$ & $x_1 x_2 x_3 -x_1^2-x_2^2-x_3 ^2+\omega_1 x_1+\omega_2 x_2+\omega_3 x_3+\omega_4$ \\ \hline 
$PV$ & $x_1 x_2 x_3 -x_1^2-x_2^2+\omega_1 x_1+\omega_2 x_2+\omega_3 x_3+\omega_4$ \\ \hline
$PV_{deg}$ & $x_1 x_2 x_3 -x_1^2-x_2^2+\omega_1 x_1+\omega_2 x_2+\omega_4$\\ \hline 
$PIV$ & $x_1 x_2 x_3 -x_1^2+\omega_1 x_1+\omega_2 x_2+\omega_3 x_3+\omega_4$\\ \hline 
${PIII^{D_6}}$ & $x_1 x_2 x_3 -x_1^2-x_2^2+\omega_1 x_1+\omega_2 x_2+\omega_4$\\ \hline 
$PIII^{D_7}$ & $x_1 x_2 x_3 -x_1^2-x_2^2+\omega_1 x_1-x_2$\\ \hline 
$PIII^{D_8}$ & $x_1 x_2 x_3 - x_1^2-x_2^2-x_2$\\ \hline 
$PII^{JM}$ & $x_1 x_2 x_3 - x_1+ \omega_2 x_2- x_3+\omega_4$\\ \hline 
$PII^{FN}$& $ x_1 x_2 x_3 - x_1^2 +\omega_1  x_1-x_2-1$\\ \hline 
$PI$& $x_1 x_2 x_3-x_1-x_2+1$ \\ \hline
\end{tabular}
\vspace{0.2cm}
\end{center}
\caption{Painlev\'e monodromy manifolds}
\label{tab:sing}
\end{table}

Note that in Table 1, we distinguish ten different monodromy manifolds, the $PIII^{D_6}$, $PIII^{D_7}$ and $PIII^{D_8}$ correspond to the three different cases of the third 
Painlev\'e equation according to Sakai's classification \cite{sakai}, and  the two monodromy manifolds  $PII^{FN}$ and $PII^{JM}$ associated to the  second Painlev\'e equation correspond to the two different isomonodromy problems found by Flaschka--Newell  \cite{FN} and Jimbo--Miwa \cite{MJ1} respectively.

Each cubic surface $\mathcal M_\phi:= {{\rm Spec}}(\mathbb C[x_1,x_2,x_3]\slash\langle\phi=0\rangle)$, is endowed with the natural Poisson bracket defined by:
\begin{equation}\label{eq:nambu}
\{x_1,x_2\}=\frac{\partial\phi}{\partial x_3},\quad \{x_2,x_3\}=\frac{\partial\phi}{\partial x_1},\quad 
\{x_3,x_1\}=\frac{\partial\phi}{\partial x_2}.
\end{equation}

In the case of PVI, this Poisson bracket is induced by the Goldman bracket on the $SL_2(\mathbb C)$ character variety of a 4 holed Riemann sphere, or  is given by the Chekhov--Fock Poisson bracket on the complexified Thurston shear coordinates. In \cite{ChMR}, all cubic surfaces were parameterised in terms of Thurston shear coordinates $s_1,s_2,s_3$ and parameters $p_1,p_2,p_3$ such that the Poisson bracket (\ref{eq:nambu}) is induced by the following flat one:
\begin{equation}\label{eq:PB-shear}
\begin{array}{lr}
\begin{array}{l}
\{s_1,s_2\}=\{s_2,s_3\}=\{s_3,s_1\}=1,\\
 \{p_1,\cdot\}=\{p_2,\cdot\}=\{p_3,\cdot\}=0,\\
\end{array}
& \hbox{for } PVI,PV,PV_{deg}, PIV, PII,PI,\\
\\
\begin{array}{l}
\{s_1,s_2\}=\{p_2,s_1\}=\{s_3,s_2\}=\{p_2,s_3\}=1,\\
\{s_2,p_2\}=2, \{s_1,s_3\}= \{p_1,\cdot\}=\{p_3,\cdot\}=0,
\end{array}
&\hbox{for } PIII^{D_6},PIII^{D_7},PIII^{D_8}.\\
\end{array}
\end{equation}

We give this parameterisation in Table \ref{table:shear}, where we think of all the monodromy manifolds as having the form\footnote{Note that in the current paper we have inverted the signs of $x_1,x_2,x_3$ compared to \cite{ChMR}.}:
\be\label{eq:mon-mf}
{\phi_{P}^{(d)}=} x_1 x_2 x_3 - \epsilon_1^{(d)} x_1^2- \epsilon_2^{(d)} x_2^2- \epsilon_3^{(d)} x_3^2 +\omega_1^{(d)} x_1  +
 \omega_2^{(d)} x_2 +\omega_3^{(d)} x_3+
\omega_4^{(d)}=0,
\ee
where $d$ is an index running on the list of the Painlev\'e cubics $PVI, PV,PV_{deg},PIV$, $PIII^{D_6},PIII^{D_7}$, $PIII^{D_8},PII^{JM},PII^{FN},PI$  and the parameters $ \epsilon^{(d)}_{i},\, \omega^{(d)}_{i}$, $i=1,2,3$ are given by:
 \be\label{eq:epsilon}
\begin{split}
 \epsilon^{(d)}_{1} =\left\{\begin{array}{ll}
 1&\hbox{ for } d= PVI, PV,PIII^{D_6},PV_{deg}, PIII^{D_7}, PIII^{D_8},PIV, PII^{FN} ,\\
0&\hbox{ for } d= PII^{JM} , PI ,\\
 \end{array}\right.\\
 \\
 \epsilon^{(d)}_{2} =\left\{\begin{array}{ll}
 1&\hbox{ for } d= PVI, PV,PIII^{D_6},PV_{deg}, PIII^{D_7}, PIII^{D_8}\\
0&\hbox{ for } d= PIV, PII^{FN} , PII^{JM} , PI ,
 \end{array}\right.\\
  \epsilon^{(d)}_{3} =\left\{\begin{array}{ll}
 1&\hbox{ for } d= PVI,\\
0&\hbox{ for } d= PV,PIII^{D_6},PV_{deg}, PIII^{D_7}, PIII^{D_8},PIV, PII^{FN} , PII^{JM} , PI .
 \end{array}\right.
 \end{split}
 \ee
and
 \bea\label{eq:omega}
 &&
 \omega^{(d)}_{1} =
-g_1^{(d)}g_\infty^{(d)}-\epsilon^{(d)}_1 g_2^{(d)} g_3^{(d)} ,\quad \omega^{(d)}_{2} =
-g_2^{(d)}g_\infty^{(d)}-\epsilon^{(d)}_2 g_1^{(d)} g_3^{(d)},\nn\\
&&
 \omega^{(d)}_{3} =
-g_3^{(d)}g_\infty^{(d)}-\epsilon^{(d)}_3 g_1^{(d)} g_2^{(d)},\\
&&
 \omega^{(d)}_{4} =\epsilon^{(d)}_2\epsilon^{(d)}_3\left(g_1^{(d)}\right)^2+\epsilon^{(d)}_1\epsilon^{(d)}_3\left(g_2^{(d)}\right)^2
 +\epsilon^{(d)}_1\epsilon^{(d)}_2\left(g_3^{(d)}\right)^2+\left(g_\infty^{(d)}\right)^2+\nn\\
 &&\qquad\quad+g_1^{(d)} g_2^{(d)}
 g_3^{(d)} g_\infty^{(d)}-4\epsilon^{(d)}_1\epsilon^{(d)}_2\epsilon^{(d)}_3,
\nn\end{eqnarray}
where $g_1^{(d)}, g_2^{(d)}, g_3^{(d)}, g_\infty^{(d)}$ are  constants related to the parameters appearing in the Painlev\'e equations as described in Section 2 of \cite{ChMR} (note that in that paper capital letters are used for the $g_i^{(d)}$).

\begin{table}[h]\label{table:shear}
\begin{center} 
\begin{tabular}{|c|c|c|} \hline 
{\tiny{P-eqs }}& Flat coordinates \\ \hline 
{\tiny{$PVI$}} &
$
\begin{array}{l}
g_1=e^{\frac{p_1}{2}}+e^{-\frac{p_1}{2}},\qquad g_2=e^{\frac{p_2}{2}}+e^{-\frac{p_2}{2}}, \qquad g_3=e^{\frac{p_3}{2}}+e^{-\frac{p_3}{2}},\\
g_\infty=e^{s_1+s_2+s_3+\frac{p_1}{2}+\frac{p_2}{2}+\frac{p_3}{2}}+e^{-s_1-s_2-s_3-\frac{p_1}{2}-\frac{p_3}{2}-\frac{p_3}{2}},\\
x_1=e^{s_2+ s_3+\frac{p_2}{2}+\frac{p_3}{2}}+e^{-  s_2-s_3-\frac{p_2}{2}-\frac{p_3}{2}}+e^{s_2- s_3+\frac{p_2}{2}-\frac{p_3}{2}}+
g_2e^{-s_3-\frac{p_3}{2}}+g_3 e^{s_2+\frac{p_2}{2}},\\
x_2=e^{s_3+s_1+\frac{p_3}{2}+\frac{p_1}{2}}+e^{-s_3-s_1-\frac{p_3}{2}-\frac{p_1}{2}}+
e^{s_3-s_1+\frac{p_3}{2}-\frac{p_1}{2}}+g_3 e^{-s_1-\frac{p_1}{2}}+g_1e^{s_3+\frac{p_3}{2}},\\
x_3=e^{s_1+ s_2+\frac{p_1}{2}+\frac{p_2}{2}}+e^{-s_1- s_2-\frac{p_1}{2}-\frac{p_2}{2}}+e^{s_1-s_2+\frac{p_1}{2}-\frac{p_2}{2}}+
g_1e^{-s_2-\frac{p_2}{2}}+
g_2 e^{s_1+\frac{p_1}{2}}.\\
\end{array}
$
\\ \hline 
{\tiny{$PV$}} & $
\begin{array}{l}
g_1=e^{\frac{p_1}{2}}+e^{-\frac{p_1}{2}},\quad g_2=e^{\frac{p_2}{2}}+e^{-\frac{p_2}{2}}, \quad g_3=e^{-s_1- s_2-s_3-\frac{p_1}{2}-\frac{p_2}{2}},
\quad
g_\infty=1,\\
x_1=e^{-s_1-\frac{p_1}{2}}+g_3 e^{s_2+\frac{p_2}{2}},\\
x_2=e^{-s_2-\frac{p_2}{2}}+e^{-s_2-2 s_1-\frac{p_2}{2}-{p_1}}+g_3 e^{-s_1-\frac{p_1}{2}}+g_1e^{-s_1- s_2-\frac{p_1}{2}-\frac{p_2}{2}},\\
x_3=e^{s_1+ s_2+\frac{p_1}{2}+\frac{p_2}{2}}+e^{-s_1- s_2-\frac{p_1}{2}-\frac{p_2}{2}}+e^{s_1-s_2+\frac{p_1}{2}-\frac{p_2}{2}}+
g_1e^{-s_2-\frac{p_2}{2}}+
g_2 e^{s_1+\frac{p_1}{2}}.\\
\end{array}
$
 \\ \hline
{\tiny{$PV_{deg}$}} & $
\begin{array}{l}
g_1=e^{\frac{p_1}{2}}+e^{-\frac{p_1}{2}},\qquad g_2=e^{\frac{p_2}{2}}+e^{-\frac{p_2}{2}}, \qquad g_3=0,\qquad 
g_\infty=1\\
x_1=e^{-s_1-\frac{p_1}{2}},\qquad
x_2=e^{-s_2-\frac{p_2}{2}}+e^{-2 s_1-s_2-p_1-\frac{p_2}{2}}+g_1 e^{-s_1-s_2-\frac{p_2}{2}-\frac{p_1}{2}},\\
x_3=e^{s_1+ s_2+\frac{p_1}{2}+\frac{p_2}{2}}+e^{-s_1- s_2-\frac{p_1}{2}-\frac{p_2}{2}}+e^{s_1-s_2+\frac{p_1}{2}-\frac{p_2}{2}}+
g_1e^{-s_2-\frac{p_2}{2}}+
g_2 e^{s_1+\frac{p_1}{2}}.\\
\end{array}
$\\ \hline 
{\tiny{$PIV$}} & $
\begin{array}{l}
g_1=e^{\frac{p_1}{2}}+e^{-\frac{p_1}{2}},\qquad g_2=e^{+\frac{p_2}{2}}, \qquad g_3=0,\qquad g_\infty=e^{-s_1-s_2-s_3-\frac{p_1}{2}},\\
x_1=e^{-2s_1-s_2-2 s_3 -p_1}+e^{-2 s_1 -s_2 -s_3-p_1},\\
x_2=e^{-2 s_1-s_2-p_1}+e^{-s_2}+e^{-2 s_1-s_2-s_3-p_1}+g_1e^{-s_1-s_2-\frac{p_1}{2}},\\
x_3=e^{-s_3} +g_2 e^{s_1+\frac{p_1}{2}}.\\
\end{array}
$
\\ \hline 
{\tiny{${PIII^{D_6}}$}} & $
\begin{array}{l}
g_1=g_3=1,\qquad g_2=e^{s_1+\frac{s_2}{2}+\frac{p_2}{2}}, \qquad 
g_\infty=e^{\frac{s_2}{2}+s_3+\frac{p_2}{2}},\\
x_1=e^{-\frac{s_2}{2}+\frac{p_2}{2}}+ e^{s_1-\frac{s_2}{2}+\frac{p_2}{2}}+e^{-\frac{s_2}{2}+s_3+\frac{p_2}{2}}+e^{s_1-\frac{s_2}{2}+s_3+\frac{p_2}{2}}+e^{s_1+\frac{s_2}{2}+s_3+\frac{p_2}{2}},\\
x_2=e^{s_1}+e^{s_1-s_2}- e^{-s_2}+e^{s_3}+2 e^{s_1+s_3}+e^{-s_2+s_3}+e^{s_1-s_2+s_3}+e^{s_1+s_2+s_3},\\
x_3=e^{-\frac{s_2}{2}-\frac{p_2}{2}}+e^{\frac{s_2}{2}-\frac{p_2}{2}}+e^{\frac{s_2}{2}+\frac{p_2}{2}}.\\
\end{array}
$\\ \hline 
{\tiny{$PIII^{D_7}$}} & $
\begin{array}{l}
g_1=1,\qquad g_2=e^{s_1+\frac{s_2}{2}+\frac{p_2}{2}}, 
\qquad x_1=e^{-\frac{s_2}{2}+\frac{p_2}{2}}+e^{s_1-\frac{s_2}{2}+\frac{p_2}{2}}+e^{s_1+\frac{s_2}{2}+\frac{p_2}{2}},\\
g_3=0,\qquad
g_\infty=e^{s_2+\frac{p_2}{2}}, \qquad 
x_2=1+2 e^{s_1}+e^{s_1-s_2}+e^{-s_2}+e^{s_1+s_2},\\
x_3=e^{-\frac{s_2}{2}-\frac{p_2}{2}}+e^{\frac{s_2}{2}-\frac{p_2}{2}}+e^{\frac{s_2}{2}+\frac{p_2}{2}}.\\
\end{array}
$\\ \hline 
{\tiny{$PIII^{D_8}$}} & $
\begin{array}{l}
g_1=1,\qquad g_2=g_\infty=e^{\frac{s_2}{2}+\frac{p_2}{2}}, \qquad g_3=0,
\qquad x_1= e^{-\frac{s_2}{2}+\frac{p_2}{2}}+e^{\frac{s_2}{2}+\frac{p_2}{2}},\\
x_2=2 +e^{-s_2}+e^{s_2},\qquad
x_3=e^{-\frac{s_2}{2}-\frac{p_2}{2}}+e^{\frac{s_2}{2}-\frac{p_2}{2}}+e^{\frac{s_2}{2}+\frac{p_2}{2}}.\\
\end{array}
$\\ \hline 
{\tiny{$PII^{JM}$ }}& $
\begin{array}{l}
g_1= g_3=g_\infty=1,\qquad g_2=e^{+\frac{p_2}{2}} ,\\
x_1=e^{-s_1}+e^{-s_1-s_3},\qquad
x_2=e^{s_3}+e^{s_1+s_3},\qquad
x_3=e^{-s_2- s_3}+e^{-s_3}.\\
\end{array}
$\\ \hline 
{\tiny{$PII^{FN}$}}& $
\begin{array}{l}
g_1=e^{-s_1-s_2-s_3},\qquad g_2=g_\infty=1, \qquad g_3=0,\qquad x_1=e^{s_2+ s_3},\\
x_2=e^{2 s_3+s_1+s_2}+e^{2 s_3+s_2}+e^{-s_1-s_2}+e^{-s_2},\qquad
x_3=e^{-s_3}+ e^{-s_2-s_3}.\\
\end{array}
$\\ \hline 
{\tiny{$PI$}}&$
\begin{array}{l}
g_1=g_2=g_\infty=1, \qquad g_3=0,\\
x_1=e^{-s_1},\qquad
x_2=e^{-s_1-s_2}+e^{-s_2},\qquad
x_3=e^{s_1+ s_2}+ e^{s_1}.\\
\end{array}
$ \\ \hline
\end{tabular}
\end{center}\caption{Flat coordinates on the Painlev\'e monodromy manifolds}
\label{tab:sing}
\end{table}

We recall that the celebrated confluence scheme of the Painlev\'e differential equations is the following diagram,
$$
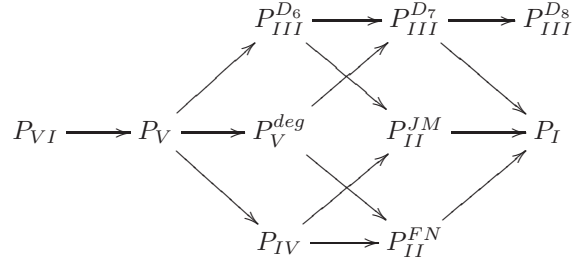
\xymatrix
 {&&P_{III}^{D_6}\ar[dr]\ar[r]&P_{III}^{D_7}\ar[dr]\ar[r]&P_{III}^{D_8}\\
 P_{VI}\ar[r]&P_{V}\ar[r]\ar[ur]\ar[dr]&P_{V}^{deg}\ar[dr]\ar[ur]& P_{II}^{JM}\ar[r]&P_I\\
 &&P_{IV}\ar[ur]\ar[r] & P_{II}^{FN}\ar[ur]\\
}
$$ 
where the arrows represent {\it confluences,} i.e. degeneration procedures where the independent variable, the dependent variable and the parameters are rescaled by suitable powers of $\varepsilon$ and then the limit $\varepsilon\to 0$ is taken. This was studied on the level of monodromy manifolds in \cite{ChMR}.

Here, we provide the quantisation of all the Painlev\'e cubics and produce the corresponding {\it quantum confluence} in such a way that quantisation and confluence commute.

To produce the quantum Painlev\'e cubics, we introduce the Hermitian operators $S_1,S_2,S_3,P_1,P_2,P_3$ subject to the commutation rule 
inherited from the Poisson bracket of $ s_1,\dots,p_3$:
$$
[P_j,\cdot]=0,\qquad [S_j,S_{j+1}]=i\pi \hbar \{ s_j,s_{j+1}\}\quad j=1,2,3,\ j+3\equiv j.
$$
Observe that  the commutators $[S_i,S_{j}]$ are always numbers and therefore we have
$$
\exp\left({a S_i}\right) \exp\left({b S_j}\right) =
\exp\left(a {S_i}+b {S_i}+\frac{ab}{2}[S_i,S_{j}]\right) ,
$$
for any two constants $a,b$. Therefore we have the Weyl ordering:
\be\label{eq:WO}
e^{S_{1}+S_{2}}=q^{\frac{1}{2}}e^{S_{1}}e^{S_{2}}=q^{-\frac{1}{2}}e^{S_{2}}e^{S_{1}},\quad q\equiv e^{-i\pi\hbar}.
\ee
After quantisation, the parameters $g_1^{(d)},\dots,g_\infty^{(d)}$ that are not equal to $0$ or $1$ become Hermitian operators $G_1^{(d)},\dots,G_\infty^{(d)}$ and are automatically Casimirs. We define the operators $\Omega_i^{(d)}$ in terms of $G_1^{(d)},\dots,G_\infty^{(d)}$ by the same formulae \eqref{eq:omega} that link the $\omega_i^{(d)}$ to the $g_i^{(d)}$'s - these are also Casimirs.  The parameters $\epsilon^{(d)}_i$ are scalars, and they remain scalar under quantisation.

{We introduce the Hermitian operators $X_1,X_2,X_3$ as follows: consider the classical expressions for $x_1,x_2,x_3$ is terms of $s_1, s_2, s_3$ and $p_1,p_2,p_3$ as in Table \ref{table:shear} in which each additive term is written as exponential of a sum of integer linear combinations of  $s_1, s_2, s_3$ and $p_1,p_2,p_3$ . Replace those exponents by their quantum version. For example the quantum version of 
$x_1$ in the $PVI$ case is 
\begin{equation*}
\begin{split}
X_1=e^{S_2+ S_3+\frac{P_2}{2}+\frac{P_3}{2}}+e^{- S_2-S_3-\frac{P_2}{2}-\frac{P_3}{2}}+e^{S_2- S_3+\frac{P_2}{2}-\frac{P_3}{2}}+\\
+(e^{\frac{P_2}{2}}+e^{-\frac{P_2}{2}})e^{-S_3-\frac{P_3}{2}}
+(e^{\frac{P_3}{2}}+e^{-\frac{P_3}{2}})e^{S_2+\frac{P_2}{2}}.
\end{split}
\end{equation*}
In this way, the resulting operators are Hermitian, because $S_1,S_2,S_3,P_1,P_2,P_3$  are.}
Then, the following result  establishes a relation between the quantisation of the Painlev\'e monodromy manifolds and the confluent Zhedanov algebra given in Definition \ref{df:UP-formal}:

\begin{prop}\label{prop-main-pq}
The Hermitian operators $X_1,X_2,X_3,\Omega_1^{(d)},\Omega_2^{(d)},\Omega_3^{(d)}$ generate the algebra $\C\langle X_1,X_2,X_3,\Omega_1^{(d)},\Omega_2^{(d)},\Omega_3^{(d)}\rangle\slash\langle J_1,J_2,J_3,J_4\rangle$ with
\bea\label{eq:q-comm-(d)}
&&
J_1 = q^{-1/2}X_1 X_2 -q^{1/2}X_2 X_1 - (q^{-1}-q) \epsilon^{(d)}_{3} X_3 +(q^{-1/2}-q^{1/2})\Omega_{3}^{(d)}, \nn\\
&&
J_2 = q^{-1/2}X_2 X_3 -q^{1/2}X_3 X_2 - (q^{-1}-q) \epsilon^{(d)}_{1} X_1 + (q^{-1/2}-q^{1/2})\Omega_{1}^{(d)},\\
&&
J_3 = q^{-1/2}X_3 X_1 -q^{1/2}X_1 X_3 - (q^{-1}-q) \epsilon^{(d)}_{2} X_2 + (q^{-1/2}-q^{1/2})\Omega_{2}^{(d)},\nn\\
&&
J_4=[\Omega_{i}^{(d)},\cdot], \qquad i=1,2,3.\nn
\eea
where $\epsilon^{(d)}_i$ are the same as in the classical case.
\end{prop}

\proof The proof of this result is obtained by direct computation by using the definitions of the quantum operators $X_1,X_2$ and $X_3$ in terms of $S_1,S_2,S_3$. By applying the quantum commutation relations for  $S_1,S_2,S_3$  (\ref{eq:WO}), relations (\ref{eq:q-comm}) follow.
\endproof

In \cite{ChMR}, we showed that the confluence procedure for the Painlev\'e differential equations corresponds to certain limits of the shear coordinates, for example for PVI to PV is obtained by the substitution $p_3\to p_3-2 \log\varepsilon$ in the limit $\varepsilon\to 0$. We define the {\it quantum confluence} by the same rescaling  the quantum Hermitian operators by  $\varepsilon$ and taking  the same limit as $\varepsilon\to 0$. For example, by imposing exactly the same limiting procedure on $P_3$, we obtain a limiting procedure on the quantum operators  $X_1,X_2,X_3,\Omega_1^{(VI)},\Omega_2^{(VI)},\Omega_3^{(VI)}$  satisfying relations \eqref{eq:q-comm-(d)} for $d=VI$, that produces some new quantum operators  $X_1,X_2,X_3,\Omega_1^{(V)},\Omega_2^{(V)},\Omega_3^{(V)}$. By construction, these operators satisfy relations \eqref{eq:q-comm-(d)} for $d=V$. The same construction can be repeated for every $d$.
 Therefore, we have the following:
 
\begin{theorem}\label{prop-main-pq1} The confluence of the Painlev\'e equations commutes with their quantisation.\end{theorem}

\section{Poincar\'e-Birkhoff-Witt (PBW)-deformation properties of the quantum algebra (\ref{eq:q-comm})}\label{se:PBW}

In this section we study the algebraic properties of the quantum algebra $\mathcal U\mathcal Z$. The basic observation is that when all constants $\epsilon_i$ and  all values of the Casimirs $\Omega_i $, $i=1,2,3$, are zero, then \eqref{eq:q-comm} are standard quantum commutation relations defining a graded algebra that is a PBW deformation of the polynomial algebra in three variables. Here we adapt the work of \cite{PP} and \cite{BG} to check that  $\mathcal U\mathcal Z$ is a PBW type deformation for all cases of $\epsilon_i$ and all values of the Casimirs $\Omega_i$, $i=1,2,3$.

\subsection{To PBW or not to PBW}\label{suse:reviewPBW}
Here we discuss the definition of the PBW, PHS and Koszul properties.
 
Let $V$ be a finite-dimensional ${\mathbb K}$-vector space  of dimension $n$ with basis $\{x_i\}_{i=1}^n$. Consider the tensor algebra  $T^\bullet(V)$  of $V$ over ${\mathbb K}$ - this is the
free associative algebra $T^\bullet(V)={\mathbb K}\langle x_1,\ldots,x_n\rangle$. For any pair of integers $1\leq i< k\leq n$ we choose an element $J_{i,k} \in T^\bullet(V)$ such that $\deg J_{i,k}  \leq 2.$ 
Let $J$ be the union of the bilateral ideals  
$$
x_i\otimes x_k - x_k\otimes x_i - J_{i,k}
$$  
in $T^\bullet(V)$. Then the
 quotient  algebra $A = T^\bullet(V)/\langle J\rangle$ is equipped with the ascending filtration $\{F_k\}, k\geq -1; F_{-1}=0$ (i.e. $F_{k-1}\subset F_k$ ) such that $F_k$ consists of all elements of degree $\leq k$ in $x_1,\dots,x_n.$
 
\begin{df}\label{PBW}
The (filtered) unital associative algebra $A$ is said to satisfy the PBW property if
there is an isomorphism of graded algebras
$$\oplus_{k\geq 0} F_k/F_{k-1} \simeq S(V),$$
where $S(V)$ is the symmetric algebra of $V$. \cite{PP}
\end{df}

Given a filtered algebra $A $ with filtration by finite-dimensional vector spaces, we write
$$P_t(A) := \sum_{k\in \mathbb Z} \dim(A_k)t^k \in {\mathbb Z}[[t]]$$
for the {\it Hilbert-Poincar\'e series of the associated graded algebra}
$${\rm gr}(A) = \oplus_{k\geq 0}A_k :=\oplus_{k\geq 0} F_k/F_{k-1}.$$

For the purposes of this paper, we distinguish the case of $n=3$ and give the following definition:

\begin{df}\label{PHS}
The algebra $A$ is said to satisfy the PHS property if its Poincar\'e-Hilbert series of $A$
coincides with $\frac{1}{(1-t)^3}.$ We shall call  \it{PHS-algebras} $3$-algebras with this property. \cite{IS}
\end{df}

In  the case  of a Lie algebra $\mathfrak g$ of dimension $n$ with a basis $\{x_1, . . . , x_n\}$,  there is a natural reformulation of the PBW-property for the universal enveloping algebra 
$U(\mathfrak g)$ in terms of the map $\sigma : S(\mathfrak g) \to {\rm gr}(U(\mathfrak g))$ where 
 $S(\mathfrak g)$ is the symmetric algebra of the Lie algebra $\mathfrak g$ and ${\rm gr}(U(\mathfrak g))$ is the associated graded algebra of the filtered algebra
 $T^{\bullet}(\mathfrak g).$ This map is defined due to the universality of $U(\mathfrak g)$  from the relation $\sigma\circ \tau = \phi$ where $\tau: T^{\bullet}(\mathfrak g) \to S(\mathfrak g)$ is the canonical projection and $\phi: T^{\bullet}(\mathfrak g)\to {\rm gr}(U(\mathfrak g))$
 is the surjective morphism of graded algebras induced by the canonical projection of $T^{\bullet}(\mathfrak g)\to U(\mathfrak g)$.
 
{In this case, the following three statements are equivalent (see \cite{Gri}):
 \begin{itemize}
 \item the homomorphism $\sigma : S(\mathfrak g)\to {\rm gr} (U(\mathfrak g))$ is a graded algebra isomorphism;
\item if $\mathfrak g$ admits a totally ordered basis $\{x_{\lambda}\}_{\lambda\in \Lambda}$ then the subset 
$$\{1\}\cup\{x_1\ldots x_{\lambda_n}\mid (\lambda_1,\ldots,\lambda_n)\in {\Lambda}^n,\quad \lambda_1\leq\ldots \leq\lambda_n,\quad n\geq 1\}$$
 gives a basis of $U(\mathfrak g)$;
\item the canonical  map $\mathfrak g \to U(\mathfrak g )$ is an injection.
 \end{itemize}
We shall use these reformulations of PBW  to choose among them a  form which is convenient to our aims.} 

We conclude this subsection by recalling the definition of Koszul algebra. 

Let $A$ be a graded algebra over a field  $\mathbb K$ of characteristic $0$:
$$
A = \oplus_{k=0}^{\infty}A_k,
$$
its  augmentation ideal $A^{+}$ is by definition
$$
A^{+} := \oplus_{k=1}^{\infty}A_k
$$ 
and the canonical projection 
$$
\pi :A\mapsto A_0=A\slash A^{+},
$$
is called augmentation map. By the augmentation map, 
$A_0$ can be considered as an $A$-module:
$$
A\times A_0 \to A_0, \quad (a,x) = \pi(a)x.
$$

\begin{df} (Koszul Algebras). 
A Koszul algebra $A$ is an $\mathbb N$-graded algebra  $A = \oplus_{k=0}^{\infty}A_k$ over a field  $\mathbb K$ that satisfies following conditions:
\begin{itemize}
\item  $A_0 = \mathbb K.$
\item $ A_0\simeq A\slash A^{+},$ considered as a graded $A$-module, admits a graded projective resolution
$$\dots \to P^{(2)} \to P^{(1)} \to P^{(0)} \to A_0 \to 0.$$
such that $P^{(i)}$ is generated as a $\mathbb Z$-graded $A$-module by its degree $i$ component, i.e., for the decomposition of $A$- modules:
$$P^{(i)} = \oplus_{j\in \mathbb Z} P^{(i)}_j$$
one has that $P^{(i)}=AP^{(i)}_i$.
\end{itemize}
\end{df}

Standard examples of Koszul algebra are the symmetric algebra $S(V)$ and the exterior algebra $\Lambda(V)$ of  an $n$ dimensional ${\mathbb K}$-vector space  $V$.

Given a Koszul algebra $A = \oplus_{k=0}^{\infty}A_k$, consider the tensor algebra $ T(A_1)$ and the map
$$
\mu : T(A_1) \to A, \quad \mu(x_1\otimes \dots \otimes x_k) := x_1\ldots x_k.
$$
A classical theorem (see for example \cite{Prid}) states that every Koszul algebra is quadratic, namely, 
$$
A \simeq T(A_1)\slash \langle I \rangle
$$
where  $\langle I \rangle$  us the ideal generated by the quadratic relation:
$$  
(\ker \mu)\cap (A_1\otimes_{\mathbb K} A_1).
$$
The inverse statement is not always true. Priddy (\cite {Prid}) proved that if a homogeneous quadratic algebra has a PBW basis, then it is
Koszul.

\subsection{PBW-type algebra structure}  In this sub-section, we follow the work by  Braverman-Gaitsgory \cite{BG} to adapt the ideas of sub-section \ref{suse:reviewPBW} to the case of non homogeneous algebras such as our 
quantum algebra $\mathcal {UZ}$.

The free  {non-commutative} polynomial associative algebra $\mathbb C\langle X_1, X_2, X_3\rangle$ can be considered as the tensor algebra $T^\bullet(V)$, where  $V={\rm Vect}\langle X_1,X_2,X_3\rangle$, that is {\it filtered} by the natural filtration:
$$
F^k (T^\bullet(V)) = \{ \oplus_{j\leq k}T^j(V)\}.
$$  
We are now going to explain how this filtration descends to the quotient. 

Fix a subspace ${\hat I}\subset F^2(T^\bullet(V)) = \C\oplus V\oplus (V\otimes V)$ and let  $I\subset V\otimes V$ be the image of ${\hat I}:  I = \pi({\hat I})$ under the natural projection $\pi : F^2(T^\bullet(V)) \to V\otimes V$. 
There is a epimorphism of graded algebras (denoted by the same letter)  $\pi : T^\bullet(V)/\langle I\rangle \to {\rm gr}(T^\bullet(V)/\langle \hat I\rangle) $.

\begin{df}\cite{BG}
The non-homogeneous quadratic algebra $\hat{\mathcal A}=T^\bullet(V)/\langle \hat I\rangle$ {\it  is a PBW-type deformation} of  ${\mathcal A} := T^\bullet(V)/\langle I\rangle$ if the projection $\pi$ is an isomorphism of
graded algebras.
\end{df}

{\begin{remark} If the subspace ${\hat I}\subset F^2(T^\bullet(V))$ (as it is in our case)  the algebra ${\mathcal A} = T (V )/\langle I\rangle$ is usually called a 2-homogeneous algebra. The deformation 
${\hat {\mathcal A}} = T (V )/\langle{\hat  I}\rangle$ is defined by two linear maps $\l_1 : I \to V$ and $\l_2 : I \to \mathbb C$ in the sense that 
$${\hat {\mathcal A}} = T (V )/\langle{\hat  I}\rangle \simeq T(V)/ (u-l_1(u) -l_2(u) \vert u\in I). $$ 
The linear maps $l_{1,2}$ satisfy some Jacobian type conditions as in Theorem \ref{th:BGPP} here below.
\end{remark}}
Roughly speaking, this means that the graded algebra ${\rm gr}(\hat{\mathcal A})$ associated to the
filtered non-homogeneous quadratic algebra $\hat{\mathcal A}=T^\bullet(V)/\langle \hat I\rangle$ is the   {2-}homogeneous quadratic ${\mathcal A} = T^\bullet(V)/{\langle I\rangle}$.

To show that  our quantum algebra $\mathcal{UZ}$ is a PBW-type  deformation, the  first step is to show that it  admits a natural filtration. This is obtained by considering it as a quotient of the free polynomial associative algebra with three generators $\mathbb C\langle X_1, X_2, X_3\rangle $ by non-homogeneous  relations { with linear and affine terms.}  Then, we need to prove that $\pi$ is indeed an isomorphism, namely we need to prove the first statement of Theorem \ref{th-main-P}:

\begin{prop}\label{prop:PBW-UP}
The quantum algebra  ${\mathcal U\mathcal Z}$   is a PBW type -deformation of the homogeneous quadratic $\C $-algebra with three generators $X_1,X_2,X_3$ and the
relations \eqref{eq:q-comm-quad}. 
\end{prop}

\proof The demonstration consists of  two steps. First, we drop linear and constant terms and consider  the ``purely" quadratic algebra $\mathcal A$ and show that it is a standard PBW-deformation of the polynomial free algebra 
$\mathbb C\langle X_1, X_2, X_3\rangle$ with three generators. By quotienting out the relations \eqref{eq:q-comm-quad} we obtain a 
graded algebra and one can easily see (choosing, for example, the base of ordered monomials $X_1^pX_2^sX_3^r$) that the dimension of the homogeneous components
of this algebra for different $s\neq 0$ is constant (flat-deformation).

As second step, we  consider the  non homogeneous algebra  $\hat{\mathcal A}$ (for generic $q$). This is based on the application of the following theorem due to Braverman-Gaitsgory \cite{BG} and Polishchuk-Positselsky \cite{PP} to the 
homogeneous ideal  ${\mathcal I}(I)$ generated by the relations $I\subset V\otimes V$ (\ref{eq:q-comm-quad}):

\begin{theorem}\label{th:BGPP}{\cite{BG}}
Let $\hat{\mathcal A}$ be a non homogeneous quadratic algebra, $\hat{\mathcal A} = T^\bullet(V)/{\langle \hat I\rangle}$, and ${\mathcal A} = T^\bullet(V)/\langle I\rangle)$ its corresponding  {2-}homogeneous quadratic algebra. Suppose $\mathcal A$ is a Koszul algebra. Then $\hat{\mathcal A}$ is a PBW-type deformation of $\mathcal A$ if and only if there exist linear functions $l_1: {I}\to V, \quad l_2 : {I}\to \C$ for which
$$
\hat I = \{u-l_1(u) -l_2(u)\mid u\in I\}.
$$
and the following conditions are satisfied
\begin{itemize}
\item  ${\rm{Im}}( l_1\otimes{\rm  Id} - {\rm  Id}\otimes l_1)\subseteq I$;
\item  $l_1  (l_1 \otimes {\rm  Id} -{\rm  Id}\otimes l_1) = -(l_2 \otimes {\rm  Id} - {\rm  Id} \otimes l_2)$,
\item  $l_2(l_1 \otimes {\rm  Id} - {\rm  Id} \otimes l_1) = 0,$
\end{itemize}
where  the maps ${l_1 \otimes {\rm  Id} -{\rm  Id} \otimes l_1}$ and ${l_2 \otimes {\rm  Id} -{\rm  Id} \otimes l_2}$ are defined on the subspace $(I\otimes V)\cap (V\otimes I) \subset T^\bullet(V).$
\end{theorem}

{\begin{remark}For the case of the finite-dimensional Lie algebra $\mathfrak g$, one has $\hat{\mathcal A} = U(\mathfrak g)$ and $\mathcal A = S(\mathfrak g)$, the symmetric algebra of $\mathfrak g.$
Consider  $I\subset \mathfrak g\otimes \mathfrak g$ defined as $I=\{x_1 \otimes x_2 - x_2\otimes x_1,\,x_1, x_2 \in \mathfrak g\}$. Then $l_1 (x_1 \otimes x_2 - x_2\otimes x_1) := [x_1,x_2],\quad l_2 :=0.$
The three conditions in Theorem \ref{th:BGPP} are equivalent to the Jacobi identity.\end{remark}}

Polishchshuk and Positselsky studied the conditions for PBW property for quadratic algebras in a more general 
setting  (\cite{PP}). We shall reformulate the conditions in theorem \ref{th:BGPP} in a form that is easy to verify in our case (see  Theorem 2.1 ch.5 in \cite{PP}) , i.e. in terms of  the bracket operator $\lbrack \cdot  ,\cdot  \rbrack:  I\subset V\otimes V  \to V$ satisfying two conditions

\beq\label{first}
  \lbrack \cdot  ,\cdot  \rbrack_{12} -  \lbrack \cdot  ,\cdot  \rbrack_{23} : (I\otimes V) \cap (V\otimes  I)\to I
\eeq

\beq\label{second}
  \lbrack \cdot  ,\cdot  \rbrack  (\lbrack \cdot  ,\cdot  \rbrack_{12} -  \lbrack \cdot  ,\cdot  \rbrack_{23}): { (I\otimes V)\cap (V\otimes I)}\mapsto 0.
\eeq

We remark that the subspace 
$$
(I\otimes V) \cap ( V\otimes  I) \subset V\otimes V\otimes V
$$
defines an analog of the space of symmetric elements of degree $3$. The bracket operation $\lbrack \cdot  ,\cdot  \rbrack:  I\subset V\otimes V  \to V$ is defined only on the subspace $I$  that is  why, due to the first condition ensures that the bracket maps $I\otimes V\cap V\otimes I$ again into $I$
and we can apply it once more.

In this setting, the map ${l_1 \otimes {\rm  Id} -{\rm  Id} \otimes l_1}$ is given by $ \lbrack \cdot  ,\cdot  \rbrack_{12} -  \lbrack \cdot  ,\cdot  \rbrack_{23}$ while the map ${l_2 \otimes {\rm  Id} -{\rm  Id} \otimes l_2}$ is
 $ \lbrack \cdot  ,\cdot  \rbrack  (\lbrack \cdot  ,\cdot  \rbrack_{12} -  \lbrack \cdot  ,\cdot  \rbrack_{23})$.

As mentioned before,  if the quadratic algebra ${\mathcal A} = T^\bullet(V)/{\langle I\rangle}$  is Koszul then the associated graded
algebra ${\rm gr}(\hat{\mathcal A})$ where $\hat{\mathcal A} = T^\bullet(V)/\langle I - [\cdot,\cdot]I \rangle$ is isomorphic to  ${\rm gr}({\mathcal A}).$ Here, $I-  [\cdot,\cdot]{I}$  means the space of elements $u - [\cdot,\cdot]u, \, u\in I$ and the ideal $\langle I - [\cdot,\cdot]I \rangle$
co-incides with the non-homogeneous ideal $\langle\hat{I}\rangle.$ 

\vskip 1mm \noindent {Let us check the conditions \eqref{first}, \eqref{second}} in the case when  $V={\mathbb C}X_1\oplus {\mathbb C}X_2 \oplus {\mathbb C}X_3$ and  $T^\bullet(V) = {\mathbb C}\langle X_1,X_2,X_3\rangle $
and $\langle\hat{ I}\rangle$ is the ideal generated by relations \eqref{eq:q-comm}:
\beq\label{qlc}
\hat {\mathcal A}= T^\bullet(V)/\langle\hat{I}\rangle.
\eeq

The first  condition \eqref{first} is valid straightforwardly. The second condition \eqref{second} follows from the following equality
\begin{equation}\label{eq:jac1}
\begin{split}
(X_1 X_2 -qX_2 X_1)X_3 + (X_2 X_3 -qX_3 X_2)X_1 + (X_3 X_1 -qX_1 X_3)X_2 =\\
X_3(X_1 X_2 -qX_2 X_1) +  X_1(X_2 X_3 -qX_3 X_2) + X_2(X_3 X_1 -qX_1 X_3),\\
\end{split}\end{equation}
that is proved by replacing the quadratic terms in the brackets by $L_1,L_2,L_3$,
where
\bea\label{eq:jac2}
&&
L_1 := (q^{-1/2}-q^{3/2}) \epsilon^{(d)}_{3} X_3-(1-q){{\Omega_{3}^{(d)}}},\nn\\
&&
L_2 := (q^{-1/2}-q^{3/2}) \epsilon^{(d)}_{1} X_1 -(1-q){{\Omega_{1}^{(d)}}},\\
&&
 L_3 := (q^{-1/2}-q^{3/2}) \epsilon^{(d)}_{2} X_2-(1-q){{\Omega_{2}^{(d)}}},\nn
\eea
leading to the identity
$$
L_1X_3 +L_2X_1 + L_3X_2 =  X_3L_1 + X_1L_2 + X_2 L_3,
$$
that is trivially satisfied due to the fact that $[L_i,X_i]=0$.

To conclude, the ``pure quadratic" part  $\mathcal A$ is Koszul hence,  the non-homogeneous algebra $\hat{\mathcal A}$ is a flat deformation of
the polynomial algebra ${\mathbb C}[X_1,X_2,X_3].$

\begin{remark}
Braverman and Gaitsgory gave a fairly simple proof that the Koszul property of $\mathcal A$ and the conditions (\ref{first}) and (\ref{second}) i.e. PBW-property
imply the existence of a graded deformation ${\mathcal A}_{\hbar}$ of $\hat{\mathcal A}$ such that at $\hbar =1$ it is {\it canonically} isomorphic to
$\hat{\mathcal A}.$ 
This is what we shall understand under ``good" (or ``flat") deformation properties.
\end{remark}

{\subsection{Zhedanov algebra and its degenerations}
As explained in Section \ref{se:mfds-q}, the quantum algebras of definition \ref{df:UZ-formal} are quantisations of the monodromy manifolds of the Painlev\'e differential equations. The  Painlev\'e sixth  monodromy manifold appeared in the paper by Oblomkov \cite{Obl} (see also \cite{EOR}) as the spectrum of the center of the Cherednik algebra of type $\check{C_1}C_1$ for $q=1$.  
This result was generalised in \cite{ChM} where this affine cubic surface was explicitly quantised leading to the Zhedanov algebra, which is isomorphic \cite{K1} to the spherical sub--algebra of the Cherednik algebra of type $\check{C_1}C_1$. In \cite{M2}, seven new algebras were produced as confluences of the Cherednik algebra of type $\check{C_1}C_1$ in such a way that their spherical--sub-algebras tend in the semi-classical limit to the monodromy manifolds of all other Painlev\'e differential equations. The quantum algebras defined by relations (\ref{eq:q-comm-quad}) are isomorphic to the spherical--sub-algebras introduced in \cite{M2}, in the same way in which the  Zhedanov algebra is isomorphic to the spherical sub--algebra of the Cherednik algebra of type $\check{C_1}C_1$. Our Theorem \ref{th-main-P}  shows that the Zhedanov algebra and its degenerations are flat deformations of the   polynomial algebra ${\mathbb C}[X_1,X_2,X_3].$}

\section{Relation with Calabi-Yau and Sklyanin algebras}\label{se:CY-Sk}

The aim of this section is to clarify the relations of our quantum  algebra $\mathcal{UZ}$ with the quantum analogues of del Pezzo surfaces introduced by Etingof and Ginzburg \cite{EtGinzb}. The latter  are elements of a very general class of non-commutative algebras related to the twisted Calabi-Yau algebras introduced by V. Ginzburg \cite{GinzbCY}. We start by recalling these notions here.

\subsection{Calabi-Yau algebras and potentials}\label{suse:CY-algebra}

Let $A$ be a finite dimensional, associative and graded $\C$-algebra . We say that $A$ is  $d$-Calabi-Yau of dimension $d$ if ${\rm Ext}^d_A (A, A\otimes A) \simeq A$ as a bimodule and otherwise ($n\neq d)$) 
${\rm Ext}^n_A(A, A\otimes A)=0$.  In this paper, we will focus  on the case of $3$-Calabi-Yau algebras. Ginzburg  has argued that most 3-Calabi-Yau algebras arise as a certain quotient of the free
associative algebra. More precisely,  let $V$ be a $\C$-vector space with  base $X_1,X_2,X_3$; its tensor algebra $T^\bullet(V)$ is the free associative graded algebra $A:=\C\langle X_1,X_2,X_3\rangle .$ One can consider the elements
of $\C\langle X_1,X_2,X_3\rangle $ as non-commutative words obtained from the variables $X_1,X_2,X_3.$  The quotient $T^\bullet(V)/[T^\bullet(V),T^\bullet(V)]$ is the space of cyclic words or ``traces". This is the $0$-degree
Hochschild homology of  the free algebra $\C\langle X_1,X_2,X_3\rangle$.  We shall use in what follows the usual notation for the quotient of an associative algebra by the space of commutators, $A_{\natural}:=A/[A,A].$ 

One can define cyclic derivatives $\partial_j \equiv \partial_{X_j}$ for any $\Phi \in A_{\natural}$ by
\beq
\partial_j \Phi :=\sum_{k|i_k=j}X_{i_{k}+1}X_{i_{k}+2}...X_{i_N}X_{i_1}X_{i_2}...X_{i_{k }-1}\in A,
\eeq
where $j=1,2,3$ and all indices $i_1,\ldots , i_N\in (1,2,3) .$
%This notion easily extends to a map from $T^\bullet(V)_{\natural}$ to $T^\bullet(V)$ and to the $T^\bullet(V)$- bimodule $T^\bullet(V)\otimes T^\bullet(V).$

The two-sided  ideal  $J_{\Phi} =<\partial_1\Phi , \partial_2\Phi, \partial_3 \Phi >$ in  $A$ is a non-commutative analogue of the Jacobian ideal and we can pass to the quotient
\beq\label {CYpot}
A_{\Phi} := A/J_{\Phi}.
\eeq
We say that the an element $\Phi\in( F^3 T^\bullet(V))_{\natural}$ is a Calabi-Yau potential if $A_{\Phi} $ is a $3$-CY-algebra.

\subsection{Etingof-Ginzburg quantisation} \label{suse:EGq}

Given a polynomial  $\phi\in \C[x_1,x_2,x_3]$ {such that the compactifying divisor of its zero set is smooth,}  Etingof and Ginzburg constructed an associative algebra $A_{\Phi}$ which is a {\it flat deformation} of the coordinate ring $\C[x_1,x_2,x_3]$ or, more precisely, the quantisation of the corresponding Poisson algebra $A_{\phi}=(\C[x_1,x_2,x_3], \{\cdot,\cdot\}_{\phi})$ where 
$$\{P,Q\}_{\phi}=\frac{dP\wedge dQ\wedge d\phi}{dx_1\wedge dx_2\wedge dx_3}$$ 
is the  Poisson-Nambu structure \eqref{eq:nambu} on $\C^3$ for $P,Q\in   \C[x_1,x_2,x_3]$  \cite{EtGinzb}.  

Let us remind that the flat deformations of a Poisson algebra $(A,\pi)$ are governed by the second group of Poisson cohomology $HP^2(A)$ and a flatness of the Poisson algebra means also a flatness of a deformation of $A$ as a commutative algebra. The flat deformations considered by 
Etingof and Ginzburg are  {\it semiuniversal deformations}  with smooth parameter scheme such that the Kodaira-Spencer map is a vector space isomorphism.

As a consequence of the computations in \cite{POR} (see Proposition 3.2), the family of affine Poisson brackets \eqref{eq:nambu} is a family of {\it unimodular} Poisson
brackets, so by the result of Dolgushev (\cite{Dolg}) the quantisation $A_{\Phi}$  is a Calabi-Yau algebra generated  by three non-commutative generators $X_i, i=1,2,3$  subject to the relations 
$$
\frac { \partial \Phi}{\partial X_1}=\frac { \partial \Phi}{\partial X_2}=\frac { \partial \Phi}{\partial X_3}=0,
$$
where $\Phi$ is a potential whose non commutative Jacobian ideal is a suitable quantum analogue of the classical Jacobian ideal in the local algebra of $\phi$.

In \cite{EtGinzb}, the authors quantise the natural Poisson structure on the hyper-surface in $\mathbb C^3$ with an isolated  elliptic singularity of type $\tilde E_r$, $r=6,7,8$. Such hyper-surfaces are the zero locus of the weighted homogeneous part of the  polynomials \eqref{eq:EG} in $\mathbb P^2$, $\mathbb{WP}_{1,1,2}$ and $\mathbb{WP}_{1,2,3}$ respectively:
\be\label{eq:EG0}
\begin{split}
\widetilde E_6 &\qquad\qquad \phi_\infty^{(6)} = \tau x_1 x_2 x_3+ \frac{x_1^3}{3}+\frac{x_2^3}{3}+\frac{x_3^3}{3}, \\
\widetilde E_7& \qquad\qquad \phi_\infty^{(7)} = \tau x_1 x_2 x_3+\frac{x_1^4}{4}+\frac{x_2^4}{4}+\frac{x_3^2}{2}, \\
\widetilde E_8 & \qquad\qquad \phi_\infty^{(8)}= \tau x_1 x_2 x_3+\frac{x_1^6}{6}+\frac{x_2^3}{3}+\frac{x_3^2}{2}.
\end{split}\ee
In each case, their quantisation produces a $3$-Calabi-Yau algebra $A_{\Phi_\infty^{(r)}}$ defined by a suitable quantum potential
$\Phi_\infty^{(r)}$, $r=6,7,8$.  Motivated by  the study of miniversal deformations of elliptic singularities,  Etingof and Ginzburg study deformations of the potential $\Phi_\infty^{(r)}$ by adding a term of the form
\be\label{eq:gen-psi-r}
\Psi_r= P(X_1) +Q(X_2) +R(X_3),
\ee
where the polynomials $P, Q$ and $R$ depend of a total of $\mu$ arbitrary parameters, $\mu$ being the Milnor number of the elliptic singularity, and have smaller degree than $\Phi_\infty^{(r)}$ has in the variable $X_1$, $X_2$ and $X_3$ respectively. They prove that for such choice of $\Psi_r$ the sum potential
$\Phi_r:=\Phi_\infty^{(r)} + \Psi_r$ also defines $3$-Calabi-Yau algebra $A_{\Phi_r}$ with central element $\Omega_r$ that in the classical limit tends to the full polynomial $\phi_r$ in \eqref{eq:EG}.

This central element, let us drop the index $r$ to keep the discussion general,  $\Omega$ is used in \cite{EtGinzb} as a non-commutative analogue of the polynomial $\phi$ and the quotient $A_{\Phi}\slash (\Omega)$ is a non-commutative analogue of the Poisson algebra $A_{\phi}\slash(\phi).$ As a consequence, 
the authors consider the following  commutative diagram where the left and right column arrows are natural surjections and the wave-like arrows denote flat deformations (or quantisations)
of the corresponding Poisson algebras $A_{\phi}$ and $A_{\phi}/(\phi):$
\beq\label{diag2}
\xymatrix{ &A_{\phi}\ar[d]\ar@{~>}[r]^{\txt{fl. def.}}& A_{\Phi}^{{q}}\ar[d]\\
& A_{\phi}/\langle\phi\rangle\ar@{~>}[r]^{\txt{fl. def.}}&  A_{\Phi}^{{q}}/\langle \Omega-\Omega^0\rangle.\\ }
\eeq
The idea of \cite{EtGinzb} is to construct the bottom-right corner algebra  as a quotient of the (family of) associative algebras $A_{\Phi}^{{q}}$ by the bilateral
ideal $\langle \Omega-\Omega^0\rangle$ generated by a central element $\Omega \in Z(A_{\Phi}^{{q}})$.

At the quantum level, the difficulty is that the potential $\Phi$ and the central element $\Omega$ are different, even though, in their classical limit, they produce the same  polynomial $\phi$. As a consequence, to complete the construction of $ A_{\Phi}^{{q}}/(\Omega)$ one needs to find the explicit expression for $\Omega$. In \cite{EtGinzb} this was done explicitly for the elliptic singularities of type $\widetilde E_6,\widetilde E_7$ and $\widetilde E_8$. 

Let us describe the $\tilde E_6$ case in some detail - as we shall see,  this is the specific case that in certain {singular} limit produces the monodromy manifolds of the Painlev\'e equations {- note that the compactifying divisor of the corresponding affine del Pezzo surfaces is always given by a triangle of lines, and therefore is not smooth. For this reason, the proofs of sub-section \ref{suse:EG-UP} are not a straight consequence of the results in  \cite{EtGinzb}}. 

It is convenient to recast the polynomial $\phi_6$ in the form
$$
\phi^{\tau,t}_{{\bf a,b,c},d} =\tau x_1x_2x_3 + \frac{t}{3}(x_1^3 +x_2^3 +x_3^3)  
+ \frac{1}{2}(a_1x_1^2 +b_1x_2^2 +c_1x_3^2) + a_2x_1 + b_2x_2 +c_2x_3 + d.
$$ 
Let us denote by $A_{\phi^{\tau,t}_{{\bf a,b,c},d}}/(\phi^{\tau,t}_{{\bf a,b,c},d})$ the coordinate ring defined by  the affine Poisson surface $\phi^{\tau,t}_{{\bf a,b,c},d} =0$ in $\C^3$.

We note that $\phi_P^{(d)}$ defined in  (\ref{eq:mon-mf}) is a  specialisation of $\phi^{\tau,t}_{{\bf a,b,c},d}$ corresponding to the choice of parameters:
\beq\label{cubic_choice}
\tau =1,\,  t=0,\, {\bf a}=(-2\epsilon_1^{(d)},\omega_1^{(d)}),\quad {\bf b}=(-2\epsilon_2^{(d)},\omega_2^{(d)}), \quad {\bf c}=(-2\epsilon_3^{(d)},\omega_3^{(d)}),\quad d=\omega_4^{(d)},
\eeq
so that the $\mathbb P^1$ bundle over the projectivisation $\mathbb P M_{\phi_P}$ of the the surface $M_P$ coincides with our general isomonodromic cubic surface  (\ref{eq:mon-mf}).

%\phi^{d}_{t,\bf a,b,c}({q})=X_1X_2X_3 - {q} X_2X_1X_3 - \frac{t}{3}(X_1^3 +X_2^3 +X_3^3)  
%+ \frac{1}{2}(a_1X_1^2 +b_1X_2^2 +c_1X_3^2) + a_2X_1 + b_2X_2 +c_2X_3 + d

\vskip2mm

Etingof and Ginzburg consider the family of homogeneous potentials
$\Phi_{EG}\in \C\langle X_1,X_2,X_3\rangle_{\natural}$ 
\be\label{homSklcub}
\Phi_{EG}=X_1X_2X_3 - {q} X_2X_1X_3 - \frac{t}{3}(X_1^3 +X_2^3 +X_3^3) ,
\ee
and show that the filtered algebra  $A^{{q}}_{\Phi_{EG}}$ with generators $X_1,X_2,X_3$ subject to the relations
\begin{align*}
&X_1  X_2-{q}   X_2  X_1 =  t  X_3^2 ,\\
&X_2  X_3-{q}  X_3  X_2  =  t  X_1^2 ,\\
&X_3  X_1-{q}  X_1  X_3   = t  X_2^2,
\end{align*}
is a $3$-CY algebra.

They then add a deformation potential  $\Psi_{EG}$ where
\be\label{homSklcubp}
\Psi_{EG}=
 \frac{1}{2}(a_1X_1^2 +b_1X_2^2 +c_1X_3^2) + a_2X_1 + b_2X_2 +c_2X_3 + d.
\ee
Note that this deformation potential is precisely in the form \eqref{eq:gen-psi-r} with
$$
P = \frac{1}{2} a_1X_1^2+ a_2X_1+ \frac{1}{3}d,\quad Q= \frac{1}{2} b_1X_2^2+ b_2X_2+ \frac{1}{3}d,
\quad R = \frac{1}{2} c_1X_3^2+ c_2X_3+ \frac{1}{3}d. $$
The sum $\Phi_{EG} + \Psi_{EG}$ depends on $q$ and further $8$
parameters ${\bf a}=(a_1,a_2),{\bf b}=(b_1,b_2),{\bf c}=(c_1,c_2), q$ and $t$ (the Milnor number of the 
corresponding elliptic Gorenstein singularity  is $8$).

The Jacobian  of the potential  $ {\Phi_{EG} + \Psi_{EG}}$ gives the following  relations 
\begin{align}
\label{eq:t-comm}
&X_1  X_2-{q}   X_2  X_1- t   X_3^2+c_1 X_3+c_2=0,\nn\\
&X_2 X_3-{q}  X_3  X_2- t  X_1^2+a_1 X_1+a_2=0,\\
&X_3  X_1-{q}  X_1  X_3- t  X_2^2+b_1 X_2+b_2=0,\nn
\end{align}
that define the family of algebras  $A^{{q}}_{ {\Phi_{EG} + \Psi_{EG}}}$.

The Theorem 3.4.4 from \cite{EtGinzb} claims that for generic values of  the  parameters $q,a_1,a_2,b_1,b_2,c_1,c_2,d,t$,  the family of algebras  $A^{{q}}_{ {\Phi_{EG} + \Psi_{EG}}}$ is Calabi-Yau. 

In fact, they prove a more general statement; in each case we may choose
\be\label{eq:EG_678j}
\Phi_{EG}=\left\{\begin{array}{lc} X_1X_2X_3 - {q} X_2X_1X_3 - \frac{t}{3}(X_1^3 +X_2^3 +X_3^3), &\hbox{for } E_6\\
 X_1X_2X_3 - {q} X_2X_1X_3 -t \left( \frac{1}{4} X_1^4 + \frac{1}{4} X_2^4 + \frac{1}{2} X_3^2\right), &\hbox{for }  E_7\\
 X_1X_2X_3 - {q} X_2X_1X_3 -t \left( \frac{1}{6} X_1^6 + \frac{1}{3} X_2^3 + \frac{1}{2} X_3^2\right), &\hbox{for }  E_8\\
\end{array}\right.
\ee
and taking $\Psi_{EG}= P(X_1) +Q(X_2) +R(X_3)$ depending on generic  $\mu +1$ parameters with $P$, $Q$, $R$ non-homogeneous polynomials of degree:
$$
{\rm deg}(P)=\left\{\begin{array}{lc} 
2, & \hbox{for } E_6,\\
 3, &\hbox{for }  E_7,\\
 5, &\hbox{for }  E_8,\\
\end{array}\right.\quad {\rm deg}(Q)=\left\{\begin{array}{lc} 
2, &\hbox{for }  E_6,\\
 3, &\hbox{for }  E_7,\\
 2, &\hbox{for }  E_8,\\
\end{array}\right. {\rm deg}(R)=\left\{\begin{array}{lc} 
2, & \hbox{for } E_6,\\
1, & \hbox{for } E_7,\\
1, &\hbox{for }  E_8,\\
\end{array}\right.
$$
the sum $\Phi_{EG}+\Psi_{EG}$ is a Calabi-Yau potential and its Jacobian defines a 
filtered family of associative $3$-Calabi-Yau algebras  with $\mu+1$ parameters, where $\mu$ is the Milnor number of the respective Gorenstein singularity.

In each case $\widetilde E_r$, $r=6,7,8$, this family of filtered algebras $A^{{q}}_{{\Phi_{EG} +\Psi_{EG} }}$ forms the {\it Rees algebras} of the corresponding algebras $A^{q}_{\Phi_{EG}}$ with homogeneous potentials $\Phi_{EG}$ given in \eqref{eq:EG_678j}.
The algebras $A^{{q}}_{{\Phi_{EG} +\Psi_{EG} }}\slash (\Omega)$ where $\Omega \in Z(A^{{q}}_{ {\Phi_{EG} + \Psi_{EG}}})$ is a non-scalar
central element, give a semi-universal family of associative algebras (depending on $q$ and $\mu$ parameters)  which are  $3$-Calabi-Yau as well.

The theorem 3.4.5 in \cite{EtGinzb} proves that the center  $Z(A^{{q}}_{\Phi_{EG}+ \Psi_{EG}})$ is the  polynomials algebra $\C[\Omega]$ and the quotient-algebra  $A^{{q}}_{ {\Phi_{EG} + \Psi_{EG}}}/\Omega$ of   $A^{{q}}_{ {\Phi_{EG} + \Psi_{EG}}}$ by the two-sided ideal 
$\langle\Omega\rangle$
gives  a flat deformation of $A_{\phi^{\tau,t}_{{\bf a,b,c},d}}/(\phi^{\tau,t}_{{\bf a,b,c},d}).$

The main difficulty in this description, as it remarked by Etingof and Ginzburg, is to compute
the explicit form of the Casimir $\Omega$.  In the case of $\widetilde E_6$, the central element is given by  \cite{EtGinzb,R}:

\begin{equation}\label{eq:omegaEG}
\begin{split}
\Omega_{EG} &= (-a_1^2 {q}^2 - a_2 {q} t - 2 a_2 {q}^2 t - a_2 {q}^3 t - b_1 c_1 {q} t^2) X_1 +\\
& t (-b_2 - 2 b_2 {q} - 2 b_2 {q}^2 - b_2 {q}^3 - a_1 c_1 {q} t + b_1^2 t^2 - b_2 t^3- b_2 {q} t^3) X_2 + \\
& t (-c_2 {q} - 2 c_2 {q}^2 - 2 c_2 {q}^3 - c_2 {q}^4 - a_1 b_1 {q} t - c_1^2 {q} t^2 + c_2 t^3 +c_2 {q} t^3) X_3 + \\
& (1 + {q}) t^2 c_1 {q} t X_2  X_1+t (-b_1 - b_1 {q} - b_1 {q}^2 - 2 b_1 t^3 - b_1 {q} t^3) X_2^2+\\
&(-a_1 {q}^2 + a_1 {q} t^3) X_2  X_3 + (1 + {q}) t^2 b_1 t X_3X_1  + (a_1 {q}^3 + a_1 {q} t^3) X_3  X_2 + \\
& t (-c_1 {q}^2 - c_1 {q}^3 - c_1 {q}^4 + c_1 t^3 + 2 c_1 {q} t^3) X_3^2 + (1 + {q}) t^2 (1 + t) (1 - t + t^2) X_2^3+\\
& (1 + {q}) t ({q}^3 - t^3) X_2  X_3X_1 - 
(1 + {q}) t (1 + t) (1 - t + t^2) {q} X_3  X_2 X_1 + ({q}^3 - t^3) (1 + {q}) t X_3^3.\\
\end{split}
\end{equation}

\subsection{Algebra $\mathcal{UZ}$ as singular limit of an Etingof-Ginzburg Calabi-Yau algebra}\label{suse:EG-UP}

In this section we prove some further nice properties of the algebra  $\mathcal{UZ}$ by showing that it is isomorphic to a singular limit of an Etingof-Ginzburg Calabi-Yau algebra. Indeed, the specialisation of relations \eqref{eq:t-comm} with
\be\label{eq:specialisation}
\begin{split}
a_1=\frac{(q^2-1)\epsilon_1^{(d)}}{\sqrt{q}},\quad b_1=\frac{(q^2-1)\epsilon_2^{(d)}}{\sqrt{q}}, \quad c_1=\frac{(q^2-1)\epsilon_3^{(d)}}{\sqrt{q}}\\
a_2=\Omega_1(1-q),\quad b_2=\Omega_2(1-q), \quad c_2=\Omega_3(1-q), \qquad t=0,
\end{split}\ee
gives the commutation relations \eqref{eq:q-comm}. The following result proves the third statement in Theorem \ref{th-main-P}:

\begin{prop}\label{prop:UZ-central-EG}
The  cubic Casimir $\Omega_4$ defined in \eqref{q-cubics} is a special limit of the Etingof-Ginzburg central element  $\Omega_{EG}$.
\end{prop}

\proof To deduce the central element  $\Omega_4$   as a limit of $\Omega_{EG}$, we first need to introduce a quadratic term $X_1^2$ in $\Omega_{EG}$ by applying the commutation relations \eqref{eq:t-comm}. Then, by taking the limit as $t\to 0$ of  $\frac{1}{t}(\Omega_{EG}-a_1 a_2(q^2+t^3))$   we obtain:
\be\label{eq:omega-eg0}
\Omega_{EG}^{t=0}:= (q^2-1) q X_3 X_2 X_1
-(q+1) (a_2 q X_1 +b_2  X_2 +c_2 q  X_3 ) - a_1 q^2 X_1^2- 
 b_1 X_2^2 - c_1 q^2 X_3^2.
\ee
The specialisation of $\Omega_{EG}^{t=0}$ with
\eqref{eq:specialisation}
is a central element in the algebra $\mathcal{UZ}$ that coincides with $(q^2-1)\sqrt{ q}\Omega_4$. 
\endproof

From this perspective, one can specialise the potential $\Phi_{EG} + \Psi_{EG}$ with the choice of parameters \eqref{eq:specialisation}. In this way, one obtains precisely the potential \eqref{qEGcub}. This potential can be decomposed as $\Phi_{\mathcal U\mathcal Z}={\Phi_{SP}}+ \Psi_{\mathcal{UZ}}$ where 
$$
\Phi_{SP} = X_1 X_2 X_3 - q X_2 X_1 X _3  \in \mathbb C\langle X_1,X_2,X_3\rangle _{\natural}
$$
is a homogeneous degree 3  potential that yields the skew polynomial algebra of three variables $X_1,X_2,X_3$ \eqref{eq:q-comm-quad} and
$$
\Psi_{\mathcal{UZ}}=\frac{(q^2-1)}{\sqrt{q}} \left(\epsilon^{(d)}_1X_1^2  + \epsilon^{(d)}_2X_2^2  +\epsilon^{(d)}_3X_3^2\right) 
+\left (q -1\right)(\Omega_3^{(d)} X_3  + \Omega_1^{(d)} X_1 + \Omega_2^{(d)} X_2),
$$
is the specialisation of ${\Psi_{EG}}$ with the choice of parameters \eqref{eq:specialisation}. Therefore we have
the following result from which the second statement of Theorem \ref{th-main-P} follows automatically:

\begin{prop}\label{prop:CY-UZ-K}
The associative algebra 
$$
A^q_{\Phi_{\mathcal U\mathcal Z} } := \C\langle X_1,X_2,X_3\rangle / \langle\partial_1 \Phi_{\mathcal U\mathcal Z},\partial_2 \Phi_{\mathcal U\mathcal Z},\partial_3 \Phi_{\mathcal U\mathcal Z}\rangle.
$$
coincides with $ \mathcal U\mathcal Z$ and is a non-homogeneous $3$-Calabi-Yau Koszul  algebra. 
\end{prop}

\proof 
To prove that $A^q_{\Phi_{\mathcal U\mathcal Z} } $ coincides with $ \mathcal U\mathcal Z$ we simply observe that 
the cyclic derivatives of the potential $\Phi_{\mathcal U\mathcal Z}$ give precisely 
the first three expressions in \eqref{eq:q-comm}. 
To prove that $A^q_{\Phi_{\mathcal U\mathcal Z} }$ is  a $3$-Calabi-Yau Koszul  algebra  we cannot apply Theorem 3.4.5 of \cite{EtGinzb}  directly  to the cubic potential \eqref{qEGcub} because of the fact that limit $t{\to}0$ is singular. Instead, 
we use the fact that, as proved in Proposition \ref{prop:PBW-UP}, this algebra is  a PBW deformation of the $3$-Calabi-Yau Koszul  algebra  $A_{{\Phi}_{SP}}$ with potential $\Phi_{SP}$ and apply  Theorem 3.1  in \cite{BT} that states that a non-homogeneous graded $3$-algebra is a Calabi-Yau Koszul  algebra if the homogeneous part is a $3$-graded Calabi-Yau Koszul  algebra.
\endproof

Then we can prove the following:

\begin{theorem}\label{th:EG-UZ}
{The quotient algebra $A^q_{\Phi_{\mathcal U\mathcal Z}}/(\Omega_{\mathcal U\mathcal Z})$ is a non-commutative  deformation of the Poisson quotient 
$A_{\phi^{1,0}_{\bf \epsilon,\omega}}/(\phi^{1,0}_{\bf \epsilon,\omega})$ of the algebra $A_{\phi^{1,0}_{\bf \epsilon,\omega}}$, where
$\phi^{1,0}_{\bf \epsilon,\omega}$ denotes the potential $\phi^{\tau,t}_{{\bf a,b,c},d} $ with the parameter choice \ref{cubic_choice},
and the following commutative diagram 
holds:}
\beq\label{diag2}
\xymatrix{ & A_{\phi^{1,0}_{\bf \epsilon,\omega}}\quad \ar[d]\ar@{~>}[r]^{\txt{fl. def.}}&\quad A^q_{{\Phi}_{\mathcal U\mathcal Z}}={\mathcal U\mathcal Z}\ar[d]\\
& A_{\phi^{1,0}_{\bf \epsilon,\omega}}/\langle \phi^{1,0}_{\bf \epsilon,\omega}\rangle\quad \ar@{~>}[r]^{\txt{fl. def.}}& \quad  A^q_{\Phi_{\mathcal U\mathcal Z}}/\langle \Omega_4-\Omega_4^0\rangle.\\
}
\eeq
\end{theorem}

\proof 
{The upper horizontal arrow  is a direct consequence of Proposition \ref{prop:CY-UZ-K}.  The left vertical downward arrow is   the Poisson algebra morphism between the Jacobian Nambu--Poisson algebra
$A_{\phi^{1,0}_{\bf \epsilon,\omega}}$ on $\mathbb C[x_1,x_2,x_3]$ and the Poisson algebra of functions on the affine cubic surface 
$\phi^{1,0}_{\bf \epsilon,\omega} =0$ in $\mathbb C^3.$  The lower horizontal and right vertical downward arrows are constructed by combining the main ideas of Theorem 3.4.5 of \cite{EtGinzb}, Theorem  \ref{th:BGPP}, Propositions \ref{prop:PBW-UP} and \ref{prop:UZ-central-EG} and Theorem \ref{th-main-P}. We remark that (after the computation of the central element by our singular limit (\ref{prop:UZ-central-EG}) the flatness argument verification is a routine procedure practically verbatim to the demonstration of propositions 7.3.1, Corollary 7.4.1  and Lemma 7.4.2 in \cite{EtGinzb}.
}\endproof

\subsection{Generalised Etingof-Ginzburg cubics}\label{suse:EG-gen}We now  replace the homogeneous part  $\Phi_{EG}$  given in  \eqref{homSklcub} by
$\Phi_{\alpha,\beta,\gamma}\in \C\langle X_1,X_2,X_3\rangle_{\natural}$ 
\be\label{nonhomSklcub}
\Phi_{\alpha,\beta,\gamma} =X_1X_2X_3 - {q} X_2X_1X_3 - \frac{1}{3}(\alpha X_1^3 +\beta X_2^3 +\gamma X_3^3)  
\ee
and consider the family of filtered algebras  $A^{{q}}_{\Phi_{\alpha,\beta,\gamma}}$
with generators $X_1,X_2,X_3$ subject to the relations
\begin{align*}
&X_1  X_2-{q}   X_2  X_1 =  \gamma X_3^2 ,\\
&X_2  X_3-{q}  X_3  X_2  =  \alpha X_1^2 ,\\
&X_3  X_1-{q}  X_1  X_3   =  \beta  X_2^2.
\end{align*}
Due to the results in \cite{IS, IS1} we know that algebras with homogeneous potentials from \eqref{nonhomSklcub} are non-commutative Koszul $3$-Calabi-Yau for certain choices of 
the parameters $\alpha,\beta,\gamma$ but  they are not always PBW or PHS.

\subsubsection{Digression}
Here we  list some alternative definitions of the PBW property used in the literature:\footnote{We are in debt to Natalia Iyudu for her patient explanation and clarification of different definitions of PBW property.}
\begin{df}\label{PBW-1}
The associative filtered algebra $A$ is a PBW-algebra if
\begin{enumerate}
\item The algebra $A$ is a Koszul and has Poincar\'e-Hilbert series $P_A (t) =\frac{1}{(1-t)^n}$.
\item The elements $x_1^{i_1}x_2^{i_2}\ldots,x_n^{i_n}$, where $i_1,\dots,i_n\in\mathbb Z$, form a linear basis.
\item There is an ordering on generators $x_1,\ldots ,x_n$ w.r.t. which the defining relations form a Gr\"obner basis.
\item The associated graded algebra is canonically isomorphic to the algebra generated by the homogeneous parts of quadratic relations.
\end{enumerate}
\end{df}

For example, the algebra of commutative polynomials satisfies (2) and in the case $n=3$ is a PHS algebra. Note that the fourth definition implies that any homogeneous algebra automatically has the PBW property.

\begin{example}\label{ex:3.4}
Let $A$ be the quantum algebra given by three generators $X_1,X_2,X_3$ and three relations
$$X_3^2 + aX_1X_2 + bX_2X_1,\quad X_2^2+ aX_3X_1 + bX_1X_3,\quad X_1^2 + aX_2X_3 + bX_3X_2$$
and the parameters 
$$(a,b)\neq (0,0),\quad(a^3,b^3)\neq (1,1),\quad (a+b)^3\neq -1.$$
This algebra (number P1, table VI in \cite{IS1}) is PBW with respect to the definitions (1) (2) and (4) in \ref{PBW-1} but not PBW in sense of the definition (3).
Conversely, the algebra $B$ given by three generators $Y_1,Y_2,Y_3$ and three relations
$$Y_1Y_2 + bY_2Y_1,\quad Y_3Y_1 + bY_1Y_3,\quad Y_2Y_3 + bY_3Y_2,\quad b\neq 0$$ (which is number PII, table VI in \cite{IS1})
is a PBW-algebra for all definitions in \ref{PBW-1}
\end{example}

\subsubsection{Calabi-Yau-Koszulity and PBW- properties of algebras whose potential is non homogeneous}
Here we consider  algebras whose potential has homogeneous
cubic part of ${\Phi_{\alpha,\beta,\gamma}}$ as well as non-homogeneous terms. 
Namely, we extend the family of algebras $A^{{q}}_{ {\Phi_{EG} + \Psi_{EG}}}$  by introducing  the   potential
$\Phi_{\alpha,\beta,\gamma} +{\Psi_{EG}}$
and considering 
the family $A^{{q}}_{\Phi_{\alpha,\beta,\gamma} +{\Psi_{EG}}}$  whose relations take the form
\be\label{eq:p,q,r-comm}
\begin{split}
X_1  X_2-{q}   X_2  X_1-    \gamma X_3^2+c_1   X_3+c_2=0,\\
X_2 X_3-{q}  X_3  X_2-  \alpha  X_1^2+a_1  X_1+a_2=0,\\
X_3  X_1-{q}  X_1  X_3-  \beta X_2^2+b_1 X_2+b_2=0.
\end{split}\ee
Inspired by  B. Shoikhet   \cite{BSh},  we call this generalised algebra family  by {\it Etingof-Ginzburg type algebras}.

The generalised Etingof-Ginzburg algebra \eqref{eq:p,q,r-comm} is a Koszul, $3$-Calabi-Yau for the cases when all constants $\alpha, \beta,\gamma$ are equal and non-zero, or only one of them is zero, or if two of the constants are equal and non-zero but $q=1$ (\cite{IS1}, Table VIII). 

Below, for $\gamma = 0$ we have computed the central element. We stress that the corresponding Etingof-Ginzburg algebras are {\it not Calabi-Yau} for generic values of $q,\, \alpha$ and $\beta$.

\begin{lm}
For $\gamma=0$, the element
\be\label{eq:cen-EG-gen}
\begin{split}
\Omega_{GEG}&= q (1 + q) (-1 + q^3) X_3   X_2   X_1+
q^3 (1 + q)\alpha X_1^3  + (1 + q)\beta X_2^3  -  a_1 q^2 (1 + q + q^2) X_1^2 \\
&   + c_1 q (1 + q)\alpha\beta X_2   X_1 - 
 b_1 (1 + q + q^2) X_2^2  - c_1 q^2 (1 + q + q^2)  X_3^2 -\\
&-q X_1 (a_2 (1 + 2 q + 2 q^2 + q^3) + b_1 c_1\alpha) + 
 X_2 (-b_2 (1 + 2 q + 2 q^2 + q^3) - a_1 c_1 q\beta) - \\
& q X_3 (c_2 (1 + 2 q + 2 q^2 + q^3) + c_1^2\alpha\beta)\\
 \end{split}
 \ee
is a central element in the algebra  \eqref{eq:p,q,r-comm}.
\end{lm}

As already mentioned, the algebra $A^{{q}}_{ {\Phi_{EG} + \Psi_{EG}}}$ is  a non-commutative Calabi-Yau algebra.
Moreover, in \cite{EtGinzb} it is shown  that the Hilbert-Poincar\'e polynomial of the algebra $A^{{q}}_{ {\Phi_{EG} + \Psi_{EG}}}$ is $\frac{1}{(1-t)^3},$ i.e. this is a PHS-algebra. Conversely, as   follows from  Example \ref{ex:3.4}, the homogeneous degree $3$ part $A^{{q}}_{\Phi^{d}_{t,\bf 0,0,0}}$ of this algebra is not a PBW-algebra in the sense of (4) in Definition \ref{PBW-1}.

%These results fit very well to the general scheme of PBW-deformations of Koszul Calabi-Yau algebras. Indeed, the results of F. Van Ostayen et all  (\cite{HVZ})  say that if  
%the algebra $T^\bullet(V)/<\mathcal I>$ is Koszul Calabi-Yau and ${\mathcal A}= T^\bullet(V)/\langle I \rangle $ is its  PBW-deformation then ${\mathcal A}$ is also a Calabi-Yau.

In the next subsection we discuss a known example of  the generalised Etingof--Ginzburg corresponding to $\alpha=\beta=\gamma=0$ for which the Etingof-Ginzburg type algebra  is ``good" Koszul Calabi-Yau.

%The examples in the subsequent  sections motivate us to call such 3-Calabi-Yau algebras {\it log-Calabi-Yau} algebras because of the ``toric" or "cluster variety" nature of the %underlying 3D affine Poisson varieties.

\subsection{\bf Odesskii algebra of Sklyanin type}\label{OdSklDeg}
In \cite{Od84}, 
Odesskii  defined a quadratic algebra ${\mathfrak O}_q$ with three generators $X_1,X_2,X_3$ satisfying the following relations:
\begin{equation}\label{eq:OdSklrel}
X_1X_2 - qX_2X_1 =X_3;\quad X_2X_3 - qX_3X_2 = X_1;\quad X_3X_1 - qX_1X_3 = X_2,
\end{equation}
and proved that, for generic $q$, the center ${\mathcal Z}({\mathfrak O}_q)$ is generated by the following element $\Omega^q := (q^2 -1)X_1X_2X_3 +X_1^2 + q^2X_2^2 + X_3^2$. When $q\to 1$ the algebra tends to the universal enveloping $U(\mathfrak{sl}_2)$. Odesskii  called the algebra ${\mathfrak O}_q$ a Sklyanin type algebra.

\begin{theorem}
The Odesskii algebra $\mathfrak O_q$ is a PBW deformation of the $3$-Calabi-Yau Koszul algebra of skew polynomials defined by  the potential 
\beq\label{ref:CYodess}
\Phi_O := \Phi_{SP} -\frac{1}{2}( X_1^2 + X_2^2 +X_3^2) \in \C\langle X_1,X_2,X_3\rangle _{\natural}.
\eeq
\end{theorem}

\proof This algebra is a PBW-algebra in the sense of all definitions in \ref{PBW-1} because of the good PBW-properties in all senses
of its homogeneous degree 3 part (see second case in the Example \ref{ex:3.4}).
To check these properties  we again apply the Theorem 3.1 in the case $N=2$  of R. Berger et R. Taillefer  \cite{BT}.
\endproof

\begin{remark}This algebra is related to the following version of a quantised universal enveloping algebra for $\mathfrak sl_2$ (\cite{MolRag}):
make a rotation in the $(X_1,X_2)$  plane:
$$X_1\to -X_2;\quad  X_2\to X_1; \quad X_3\to X_3$$ and then the rescaling
\be\label{eq:res-o}
X_1\to (q-q^{-1})X_1;\quad  X_2\to (q-q^{-1})X_2; \quad X_3\to (q-q^{-1})X_3,\ee
maps the Odesskii algebra to the algebra with relations
\beq\label{molrag}
qX_1X_2 - X_2X_1 =(q-q^{-1})X_3;
\eeq
\begin{align*}
&qX_2X_3 - X_3X_2 =(q-q^{-1}) X_1;\\
&qX_3X_1 - X_1X_3 = (q-q^{-1})X_2\\
\end{align*}
and with the Casimir
$${\tilde \Omega}_O := -qX_1X_2X_3 +q^2X_1^2 + X_2^2 + X_3^2.$$
\end{remark}

\begin{remark}
This quantum Casimir cubic goes to the famous {\it Markov cubic} in the limit $q\to 1.$ 
\end{remark}

\subsection{\bf Sklyanin algebra with three generators}\label{skl3}
One of the most famous examples of a $3$-Calabi-Yau algebra is the graded associative algebra $Q_3({\mathcal E, a,b,c})$ {which is related to  a  (possibly degenerate or singular) elliptic curve ${\mathcal E}$}
\beq\label{eq:sklyanin-algebra}
Q_3({\mathcal E, a,b,c})= \C\langle X_1,X_2,X_3\rangle /J_\Phi
\eeq
with 
$$
J_\Phi =\langle aX_2X_3 +bX_3X_2 +cX_1^2,\quad aX_3X_1 +bX_1X_3 +cX_2^2,\quad aX_1X_2 +bX_2X_1 +cX_3^2\rangle,
$$
where $(a,b,c)\in \C^3$ are some parameters.  This algebra is a special sub-case of  the one generated by $\Phi_{EG}$ with $q=\frac{b}{a}$ and $t=\frac{c}{a}$.

{Artin and Schelter \cite{AS} proved that, if the parameters  $(a,b,c)\in \C^3$ define the homogeneous coordinates of a point in $\mathcal E$, } this algebra satisfies  the Poincare-Birkhoff-Witt condition for all definitions in \ref{PBW-1} except {(3)} and hence  it can be considered as a deformation of the polynomial ring $\mathbb C[x_1,x_2,x_3]$. For this reason, this algebra is often called the Artin--Schelter--Tate--Sklyanin algebra with three generators, but in this paper, for brevity, we call   ``Sklyanin algebra"  any  graded associative algebra with quadratic relations which satisfies the Poincare-Birkhoff-Witt {or PHS-}conditions and that can be considered as a deformation of 
the polynomial ring $\mathbb C[x_1,x_2,x_3].$  
{Iyudu and Shkarin \cite{IS} have proved that this algebra is a CY algebra.}

\subsection{\bf Generalised Sklyanin algebras with three generators}

Iyudu and Shkarin (\cite{IS}) introduced the {\it generalised Sklyanin algebra} with three generators as the following quotient of the free associative algebra 
\beq\label{genSkl}
{\tilde Q_3}(a,b,c,\alpha,\beta,\gamma)= \C\langle X_1,X_2,X_3\rangle /J_{GS}
\eeq
where
$$J_{GS} =\langle X_2X_3 -aX_3X_2 -\alpha X_1^2,\, X_3X_1 -bX_1X_3 -\beta X_2^2,\, X_1X_2 -cX_2X_1 -\gamma X_3^2\rangle,$$
where  $(a,b,c,\alpha,\beta,\gamma)\in \C^6$ is a generic set of complex constants.

These generalised Sklyanin algebras are not always potential and in fact, for generic $(a,b,c,\alpha,\beta,\gamma)$ they have neither good PBW-properties nor Koszul properties.
However, for special values of the parameters they do and a complete classification is given in the following  result \cite{IS}: 

\begin{theorem}\label{gensklhom}
The generalised Sklyanin algebra is PHS if and only if at least one of the following conditions is satisfied:
\begin{enumerate}
\item For $a=b=c\neq 0$ and $(a^3,\alpha\beta\gamma)\neq (-1,1)$ - this case includes the quadratic Sklyanin algebra $Q_3(\mathcal E,a,c,\frac{\alpha}{3})$.
\item For $(a,b,c)\neq(0,0,0)$ and either $\alpha=\beta=a-b=0$ or $\gamma=\alpha=c-a=0$ or $\beta=\gamma=b-c=0$.
\item For a specific choice of all parameters in terms of a root of unity, it is a ``finite" algebra which is out of our interest.
\item For $a=b=c=0$ and $\alpha\beta\gamma\neq 0$. This algebra is potential without the cubic term $X_1 X_2 X_3$ and  is out of our interest.
\item For $\alpha=\beta=\gamma=0$ and  $(a,b,c)\neq(0,0,0)$, this is the case of the skew polynomial algebra.
\end{enumerate}
In all these cases, the  generalised Sklyanin algebra is potential and Koszul. The potential can be written as follows:
$$
\Phi_{GS}=  \frac{1}{3}(\alpha X_1^3 +\beta X_2^3 + \gamma X_3^3 ) + \tilde a X_1X_2X_3  +\tilde b X_2X_1X_3,
$$
where $\tilde a$ and $\tilde b$ depend on $a,b,c,q$.
\end{theorem}

\subsection{Generalised Sklyanin-Painlev\'e potential.}\label{se:pain-pot}
Motivated by the idea of merging together generalised Sklyanin algebra and our algebra \eqref{eq:q-comm}, we 
consider the following generalisation of the potential of Etingof and Ginzburg to include the first  two cases of Theorem \ref{gensklhom}:
\begin{equation}\label{eq:EGprep-gen}
\Phi=\Phi_{GS}+ \Psi_{EG}
\end{equation}
For the choice of parameters as in the cases of Theorem \ref{gensklhom}, the algebra
$$
{\mathcal A}^q :=C\langle X_1,X_2,X_3\rangle /J,
$$
where 
$$
J=\langle \partial_{X_1} \Phi, \partial_{X_2} \Phi, \partial_{X_3} \Phi\rangle
$$
is a {\it generalised Sklyanin-Painlev\'e algebra}\/ \eqref{eq:q-comm-GSP} and
gives a PHS- or PBW-type 3-Calabi-Yau deformation of $C\langle X_1,X_2,X_3\rangle /J_{\Phi_{GS}}$.

For special choices of the parameters $\alpha,\beta,\gamma,a,b,c,a_1,b_1,c_1,a_2,b_2,c_2$, the space of objects ${\mathcal  A}^q_\natural:= {\mathcal A}^q\slash[{\mathcal A}^q,{\mathcal A}^q]$  appears in the Physics literature to which section \ref{se:NC-QFT} is dedicated. 

\begin{remark}
It is interesting to observe the correspondence between the conditions on the constants $\alpha,\beta$ and $\gamma$ in the generalised Etingof-Ginzburg algebras and the {\it semi--classical} conditions on 
the existence Poisson-Nambu 3$d$ polynomial algebras in the the paper of L. Vinet and A. Zhedanov (\cite{VZ}) where they classify various Poisson analogues of the Askey-Wilson algebras AW(3).
This correspondence could be behind the fact that the potential and the central element in the generalised Etingof-Ginzburg algebras are, in general, different, despite having the same semi-classical limit.
\end{remark}

\section{Poisson structures and degenerations of elliptic curves}\label{se:GHK-theta}

Motivated by the observation by M. Gross, P. Hacking and S.Keel (see Example  6.13 of \cite{GrossHK})   that the family associated to 
 \eqref{eq:mon-mf} is a log-symplectic Calabi-Yau variety, or in other words, that the projective completion $Y$ of \eqref{eq:mon-mf}  with the cubic divisor {$D_\infty$} given by a triangle of lines, is an example of a  Looijenga pair,  in this section we study the degenerations of  a certain class of Looijenga pairs $(Y,D)$.

In this context we need to fix some notation and assumptions to make our discussion clear. We consider the polynomials $\phi\in \mathbb C[x_1,x_2,x_3]$ of the form \eqref{eq:EOR}, or \eqref{eq:EG} or belonging to Table 1. We list all such polynomials in the first column of Table 3.

\begin{table}[h]
\begin{center} 
\begin{tabular}{|c|c|c|c|} \hline 
 Polynomials $\phi$ & $\begin{array}{c}\delta\\ {\rm weights}\\ \end{array}$ & $\phi_\infty$ \\ \hline 
 $\frac{x_1^6}{6}+\frac{x_2^3}{3}+\frac{x_3^2}{2}+ \tau x_1 x_2 x_3+\eta_5 x_1^5+\dots +\omega ,$ & $\begin{array}{c} 1\\(1,2,3)\\ \end{array}$  & $\frac{x_1^6}{6}+\frac{x_2^3}{3}+\frac{x_3^2}{2}+ \tau x_1 x_2 x_3$  \\ \hline
 $\frac{x_1^4}{4}+\frac{x_2^4}{4}+\frac{x_3^2}{2}+ \tau x_1 x_2 x_3+\eta_3 x_1^3+\dots +\omega , $ &$\begin{array}{c} 2\\ (1,1,2)\\ \end{array}$  & $\frac{x_1^4}{4}+\frac{x_2^4}{4}+\frac{x_3^2}{2}+ \tau x_1 x_2 x_3$ \\\hline
 $\frac{x_1^3}{3}+\frac{x_2^3}{3}+\frac{x_3^3}{3}+ \tau x_1 x_2 x_3+\eta_2 x_1^2+\dots +\omega  ,$ & $\begin{array}{c} 3\\ (1,1,1)\\ \end{array}$  & $\frac{x_1^3}{3}+\frac{x_2^3}{3}+\frac{x_3^3}{3}+ \tau x_1 x_2 x_3$ \\\hline
$x_1 x_2 x_3 + x_1^5 + x_2^2 + x_3^2  +\eta_4 x_1^4+\dots +\omega$, &$\begin{array}{c}1\\ (2,5,3)\\ \end{array}$  & $x_1 x_2 x_3 + x_1^5 + x_2^2$ \\\hline
$x_1 x_2 x_3 + x_1^4 + x_2^2 + x_3^2 +\eta_3 x_1^3+\dots +\omega $, & $\begin{array}{c}2\\  (1,2,1)\\ \end{array}$ & $x_1 x_2 x_3 + x_1^4 + x_2^2 $ \\\hline
$x_1 x_2 x_3 + x_1^3 + x_2^3 + x_3^2  +\eta_2x_1^2+\dots +\omega$, & $\begin{array}{c}3\\ (1,1,1)\\ \end{array}$  & $x_1 x_2 x_3 + x_1^3 + x_2^3 $ \\\hline
$ x_1 x_2 x_3 + \sum_{k=1}^3(\omega_k x_k- \epsilon_k  x_k^2) +
\omega_4 $ &  $\begin{array}{c}3\\(1,1,1)\\ \end{array}$ & $x_1 x_2 x_3$\\ \hline 
\end{tabular}
\vspace{0.2cm}
\end{center}
\caption{del Pezzo surfaces as Loojenga pairs - in the last row we dropped the index ${}^{(d)}$.}
\label{tab:sing}
\end{table}

 In each case,  the projective completion $\overline{\mathcal M}_\phi$ of $\mathcal M_\phi$  in the weighted projective spaces $\mathbb W \mathbb P^3$  are del Pezzo surface of degree $\delta$ \cite{EOR,{EtGinzb}}. We denote by $(x_0,\dots,x_3)$ the weighted homogeneous coordinates in $\mathbb W \mathbb P^3$. We list the degree $\delta$ and the weights of the variables $(x_1,x_2,x_3)$ in the second column - we always assume the weight of the homogeneous coordinate $x_0$ to be $1$.

For each polynomial  $\phi\in \mathbb C[x_1,x_2,x_3]$ in Table 3, we take the weighted homogeneous part $\phi_\infty$ and list it in the third column. The equation $\phi_\infty=0$ defines a projective curve in $\mathbb W\mathbb P^2$. 
The pair $(\overline{\mathcal M}_\phi,D_\infty)$ is a Looijenga pair and $\overline{\mathcal M}_\phi\setminus D_\infty$ is the affine surface $\mathcal M_\phi\in\mathbb C^3$.

The projectivisation $\mathbb P\mathcal M_\phi$ of $\mathcal M_\phi$ is a projective manifold of dimension $1$ embedded in $\mathbb P^2$ by the linear system given by sections of a line bundle of degree $\delta$ - for degenerated cubic divisors of del Pezzo degree 2 and 1 such sections are expressed via  Gross-Hacking-Keel $\theta$-functions.

The coordinate ring of  $\overline{\mathcal M}_\phi\setminus D_\infty$ is 
$\mathbb C[x_1,x_2,x_3]/\langle\phi\rangle$, which corresponds to the cone over the projectivisation $\mathbb P\mathcal M_\phi\in\mathbb P^2$, namely
$$
\mathbb C[x_1,x_2,x_3]/\langle\phi\rangle = \oplus_k H^0(\mathbb P\mathcal M_\phi,  L^{\otimes k}),
$$
where $ L$ is the trivial bundle of degree $\delta$. By taking the generalisation of the Poincar\'e residue for weighted projective spaces (see for example, \cite{Den}) of the global $3$-form in $\mathbb W \mathbb P^3$ along the divisor $D_\infty$, one obtains a symplectic form on the quotient  $\mathbb C[x_1,x_2,x_3]/\langle\phi\rangle$ which descends from the Nambu bracket restricted to the symplectic leaves $\phi=0$. 

In this Section we carry out the above construction for each $\phi$ in Table 3. We also consider special cases and  {\it degenerations,}\/ namely singular limits obtained by rescaling the weighted homogeneous coordinates and taking limits of such rescaling to infinity.
We show that such degenerations correspond to  rational degenerations of elliptic curves.

\subsection{\bf Degenerations of the Sklyanin algebra with three generators}\label{skl3}

In this subsection we consider $\phi_\infty= \frac{x_1^3}{3}+\frac{x_2^3}{3}+\frac{x_3^3}{3}+ \tau x_1 x_2 x_3$, a special case of the third row of Table 3. This case is related to the semi--classical limit of the Sklyanin algebra  \eqref{eq:sklyanin-algebra};
namely, take $a,b$ such that $a+b$ is proportional to $1-q$, the  semi---classical limit  $q_3(\mathcal E,\tau)$ (where $\tau=\frac{c}{3}$) of the Sklyanin algebra $Q_3({\mathcal E, a,b,c})$  carries a Poisson structure (which is also called {\it Poisson Sklyanin algebra).}\/ In \cite{FO} and \cite{Pol} it was  shown that this Poisson algebra belongs to a {family of} Poisson structures on the moduli space of parabolic vector bundles of degree 3 and rank 2 on the projective space $\mathbb P^2$. The explicit expression for the elliptic Poisson brackets of $q_3(\mathcal E,\tau)$ is the natural one carried by the family of the Hesse cubics 
\begin{equation}\label{eq:hesse}
{\phi}_{\tau} = \frac{1}{3}(x_1^3 + x_2^3 + x_3^3) +\tau x_1x_2x_3  =0
\end{equation}
that define the embedding of $\mathcal E$ in $\mathbb P^2$.
Namely the quadratic brackets on the affine space $\mathbb C^3$ which define a quadratic Poisson algebra structure on 
$$
A_{{\phi}_{\tau}} = {\mathbb C}[x_1,x_2,x_3]\slash{\phi_\tau}= \oplus_{k\geq 0}H^0 ({\phi}_{\tau},L^{\otimes k}) 
$$ 
and $L$ is the degree 3 line bundle over the cubic curve ${\phi}_{\tau}$ are: 
\begin{equation}\label{eq:poi-tau}
\{x_1, x_2\} = x_3^2 + \tau x_1x_2;\quad \{x_2, x_3\} = x_1^2 +\tau  x_2x_3;\quad \{x_3, x_1\} = x_2^2 + \tau x_3x_1.
\end{equation}

It a straightforward computation to check that the algebra $q_3(\mathcal E,\tau)$ is invariant under the Heisenberg group $H_3$ and
unimodular (see \cite{POR}).

\subsubsection{Rational degenerations of Sklyanin Poisson algebra and triangular divisor of Painlev\'e projective surfaces}

 A. Odesskii in \cite{OdRat} proposed a description of all rational degenerations for  a generalisation  of elliptic algebras known as Sklyanin--Odesskii--Feigin algebras,  and their semi--classical counterparts - namely rational Poisson quadratic algebras.  We shall restrict ourselves to one example of it in the case of the  Poisson elliptic algebra {$q_3(\mathcal E,\tau).$}  It is shown in \cite{OdRat} that the center of the rational degeneration $R^1_{3}(-\frac{2}{3})$ of the Sklyanin algebra $Q_3(\mathcal E {,a,b,c})$ is generated by one polynomial
of degree 3  in {$\mathbb P^2$}
$$
{\tilde \phi}= \frac{1}{3}y_2^3+ y_1y_2y_3.
$$ 
Indeed, if  we take the Casimir element  $\phi$ of  $q_3(\mathcal E,\tau)$ given by the Hesse cubic \eqref{eq:hesse} and 
take  the rational limit $\tau\to \infty$ which gives us the triangle configuration in Figure 1
$$
\{x_1=0\}\cup \{x_2=0\} \cup \{x_3=0\}
$$ 
in the coordinates  $y_i, i=1,2,3$ defined as:
$$
y_1 = \sqrt{\tau}x_1,\quad  y_2 = x_2,\quad y_3 = \sqrt{\tau}x_3,
$$
we obtain $\tilde\phi$.

The same triangle configuration is the divisor at infinity of the projective completion ${\bar\phi_{P}} \hookrightarrow {\mathbb P}^3$  of the general Painleve cubic \eqref{eq:mon-mf}.

\begin{remark}
It is clear that, in the limit $\tau\to\infty$, the Poisson brackets \eqref{eq:poi-tau}
give  the  cluster Poisson structure (\cite{GShV})
$$
\{x_1, x_2\} =  x_1x_2;\quad \{x_2, x_3\} = x_2x_3;\quad \{x_3, x_1\} =  x_3x_1,
$$ 
but in the degenerated coordinates $y_1, y_2, y_3$ these brackets read
$$\{y_1,  y_2\} =  y_1 y_2;\quad \{ y_2,  y_3\} =  y_2  y_3;\quad \{ y_3,  y_1\} =  y_2^2 +  y_3y_1.
$$
Because these are   brackets on $\mathbb C^3$, they define a quadratic Poisson algebra structure on 
$$A_{\tilde{\phi}} = {\mathbb C}[y_1,y_2,y_3]/{\tilde{\phi}} = \oplus_{k\geq 0}H^0 ({\tilde\phi},L^{\otimes k})
$$ 
where $L$ is the degree $3$ line bundle over the cubic divisor $ \frac{1}{3}y_2^3+ y_1y_2y_3 = 0$  which is the union of the line $y_2 =0$ and the conic $ \frac{1}{3} y_2^2 +y_1y_3=0$.
The rational degeneration deforms the  cluster Poisson structure.
We will consider the quantum version of this in subsection \ref{suse:LBW}.
\end{remark}

\subsubsection{\bf Elliptic curves in weighted projective spaces, related Sklyanin Poisson structures and their rational degenerations.}\label{suse:E7-E8}

We deal first with the polynomial $\phi$ in the third row of Table 3. As discussed in subsection \ref{suse:EGq}, it is convenient to write this polynomial in the form
$$
\phi^{\tau,t}_{{\bf a,b,c},d} =\tau x_1x_2x_3 + \frac{t}{3}(x_1^3 +x_2^3 +x_3^3)  
+ \frac{1}{2}(a_1x_1^2 +b_1x_2^2 +c_1x_3^2) + a_2x_1 + b_2x_2 +c_2x_3 + d.
$$

%%%

The projectivisation $\P M_{\phi^{\tau,t}_{{\bf a,b,c},d}}$ of the hypersurface $M_{\phi^{\tau,t}_{{\bf a,b,c},d}}$
is a curve in $\mathbb P^2$ and $M_{\phi^{\tau,t}_{{\bf a,b,c},d}}$ can be seen as a line-bundle over $\P M_{\phi^{\tau,t}_{{\bf a,b,c},d}}$.
When $t=1, {\bf a}={\bf b}={\bf c}= d = 0$ the surface $M_{\phi^{\tau,1}_{{\bf 0,0,0},0}}$ is an affine cone over a {normally} embedded elliptic curve in $\P^2$ of degree $3$ given by the homogeneous cubic
$$
\bigl\{\phi_\infty = \tau x_1x_2x_3 + \frac{1}{3}(x_1^3 +x_2^3 +x_3^3)=0\bigr\}\subset \mathbb P^2.
$$

This cone surface $M_{\phi^{\tau,1}_{{\bf 0,0,0},0}}$ is an example of a simple elliptic  Gorenstein singularity ($\widetilde E_6$ case corresponding to the elliptic singularities list).
Note that the same formula for $\phi$  also defines a  hypersurface in $\mathbb C^3$ with a  triple point singularity in $0$. 

Let us now deal with the first two lines of Table 3. Denote by  $\phi_{1,1,2}$ and $\phi_{2,1,3}$ the $\phi_\infty$ in the second and first row respectively.

The surface
\begin{equation}\label{eqref: deg2}
\bigl\{ \phi_{1,1,2} = \tau_1 x_1x_2x_3 + \frac{1}{4}x_1^4 +\frac{1}{4}x_2^4 +\frac{1}{2}x_3^2 =0\bigr\} \subset \mathbb{C}^3
\end{equation}
has a  double point  in  $\mathbb C^3$ that is an  elliptic Gorenstein singularity of  type $\widetilde E_7$. It defines the affine cone over a  homogeneous degree $4$ elliptic curve in weighted projective space $\mathbb{WP}_{1,1,2}$ defined by the same equation $ \phi_{1,1,2} =0$. Similarly, the surface of  type $\widetilde E_8$ 
\begin{equation}\label{eqref: deg1}
\bigl\{ \phi_{2,1,3} = \tau_2 x_1x_2x_3 + \frac{1}{3}x_1^3 +\frac{1}{6}x_2^6 +\frac{1}{2}x_3^2=0\bigr\} \subset \mathbb{C}^3
\end{equation}
which is the affine cone over a homogeneous degree $6$ elliptic curve in weighted projective space $ \mathbb{WP}_{2,1,3}$ defined by the same equation $ \phi_{2,1,3} =0$.

From an algebraic point of view the coordinate rings $A_{\phi}$ discussed in subsection  \ref{suse:EGq}, {both} $A_{\phi_{1,1,2}}$ and $A_{\phi_{2,1,3}}$ are graded rings  such that

\begin{enumerate}
\item $A_{\phi} =\mathbb C[x_1, x_2,x_3]/\phi= \oplus_{k\geq 0} H^{0}(\phi, L^{\otimes k}) $ where $L$ is the degree 3 line bundle over the cubic curve $\phi$ and the sections of $L$ form the linear system\footnote{An explicit construction of linear systems defined by sections of $L$ for degree 2 and 1 in terms of appropriate theta functions similar to  this case  can be found, for example, in the paper \cite{RoanIsing}.} defining the embedding $\phi \hookrightarrow \mathbb P^2$;

\item $A_{\phi_{1,1,2}} = \oplus_{k\geq 0} H^{0}(\phi_{1,1,2}, L^{\otimes k}) =\mathbb C[x_1, x_2,x_3]/\phi_{1,1,2},$ where $L$ is the degree 2 line bundle over the nodal curve 
$\phi_{1,1,2}$ and the sections of $L$  define the embedding $\phi_{1,1,2} \hookrightarrow \mathbb{WP}_{1,1,2}$. 

\item $A_{\phi_{2,1,3}} = \oplus_{k\geq 0} H^{0}(\phi_{2,1,3}, L^{\otimes k}) =\mathbb C[x_1, x_2,x_3]/{\phi_{2,1,3}} $ where $L$ is the degree 1 line bundle over the nodal curve 
$\phi_{2,1,3}$ and the sections of $L$  define the embedding $\phi_{2,1,3} \hookrightarrow \mathbb{WP}_{2,1,3}.$ 
\end{enumerate}

{In the seminal paper by Gross, Hacking and Keel \cite{GrossHK}, the notion of vertex of degree $n$, $\mathbb V_n$,  was introduced. 
%For $n\geq 3$ the $n-$vertex $\mathbb V_n$ is a union of affine coordinate planes in  $\mathbb A^n:$ $$\mathbb V_n := \mathbb A^2_{y_1,y_2}\cup \mathbb A^2_{y_2,y_3}\cup ...\cup\mathbb A^2_{y_n,y_1} \subset \mathbb A^n_{y_1,..,y_n}.$$
These surfaces $\mathbb V_n$ are the most important singularities for moduli of smooth surfaces.  
We shall be specially interested in the degenerate cases of vertices $\mathbb V_{1,2}$ of degree $1$ and $2$ which can be defined as follows:
$$\mathbb V_1 = {\rm Spec}(\mathbb C[y_1,y_2,y_3]/(y_1y_2y_3 - y_1^2 -y_3^2))$$ and
$$\mathbb V_2 = \mathbb A^2_{y_1,y_2}\cup \mathbb A^2_{y_2,y_1} = {\rm Spec}(\mathbb C[y_1,y_2,y_3]/(y_1y_2y_3 -y_3^2)).$$
The first surface $\mathbb V_1$ is the affine cone over a nodal cubic curve in the weighted projective space $\mathbb{WP}_{2,1,3}$  and the second, $\mathbb V_2$, is the affine cone over the union of two rational curves embedded in the weighted projective space $\mathbb{WP}_{1,1,2}.$ }

{Both affine cones $V_1$ and $V_2$, being affine surfaces in $\mathbb C^3$ carry the natural  Poisson Jacobi unimodular structures defined by the potentials  $\phi_1 = y_1y_2y_3 - y_1^2 -y_3^2$ and 
$\phi_2 = y_1y_2y_3  -y_3^2$ respectively. P.Bousseau  has identified these Jacobian Poisson brackets with the natural Poisson structures on classical Gross-Siebert theta-functions (Th. 21 and Corollary 22 of \cite{Bouss}), defined in terms of so called ``broken lines" data.}

{\begin{prop}\label{prop:degen}
The vertex surfaces $\mathbb V_1$ and $\mathbb V_2$ coincide with degenerations of the affine cones $\mathcal M_{\phi_{2,1,3}}$ and $\mathcal M_{\phi_{1,1,2}}$ respectively. Their Jacobian Poisson algebras are isomorphic to
degenerations of elliptic Poisson Sklyanin algebras defined by the weighted cubic Casimirs $\phi_{2,1,3}$ and ${\phi_{1,1,2}}$ respectively.
\end{prop}}
\proof We  apply the same procedure of degeneration as above, namely we rescale
$$
x_1 \to y_1 3^{1/3},\quad  x_2 =  -\frac{ y_2}{\tau_2 2^{1/2}3^{1/3}},\quad x_3 = y_3 2^{1/2},
$$
and  take the limit $\tau_2 \to \infty$ to obtain
$$ 
  \phi_{{2,1,3} _0} = y_1^3 +  y_3^2  - y_1  y_2 y_3.$$
The corresponding Jacobian Poisson brackets read as
\beq\label{Grosstheta_1}
\{ y_1, y_2\} =  2 y_3 - y_1 y_2, \quad  \{ y_2,  y_3\} =    3y_1^2 -  y_3   y_2,\quad  \{  y_3, y_1\} =-y_1  y_3
\ee
and define a Posson algebra structure on the ring
 \beq\label{eqref: 123tilde}
 A_{  \phi_{{2,1,3}_0}}:=\mathbb C[y_1 ,y_2,  y_3]/ \phi_{{2,1,3}_0} = \oplus_{k\geq 0} H^{0}(  \phi_{{2,1,3} _0}, L^{\otimes k}) ,
 \eeq
where $L$ is degree 1 line bundle over the singular curve $\tilde \phi_{{2,1,3} _0}=0$, i.e.  the rational nodal cubic of arithmetic genus $1$ embedded in $\mathbb{WP}_{2,1,3}$ . 
Then
$$\mathcal M_{ \phi_{{2,1,3}_0}} :={\rm Spec}(A_{  \phi_{{2,1,3}_0}})$$ 
is the affine cone in $\mathbb C^3$ over the singular curve $  \phi_{{2,1,3} _0}=0$, {namely $\mathbb V_1.$}

Similarly by 
$$
  x_1 = -\frac{1}{2^{1/4}\sqrt{\tau_1}} y_1,\quad   x_2 =  \frac{1}{2^{1/4}\sqrt{\tau_1}}y_2,\quad  x_3 = \sqrt{2}y_3,
$$
in the limit $\tau_1 \to \infty$ one has
\be\label{eq:phi-112-0}
\phi_{{1,1,2}_0} = y_3^2  - y_1y_2y_3.
\ee
The corresponding Jacobian Poisson brackets read as
\beq\label{nodcub_3}
\{y_1, y_2\} =  2y_3  - y_1y_2, \quad  \{y_2,y_3\} = - y_3y_2,\quad  \{y_3,y_1\} =  - y_1y_3, 
\ee
and define a Posson algebra structure on the ring
 \beq\label{eqref: 112tilde}
 A_{ \phi_{{1,1,2}_0}} : =\mathbb C[y_1, y_2,y_3]/ \phi_{{1,1,2}_0}= \oplus_{k\geq 0} H^{0}( \phi_{1,1,2}, L^{\otimes k}),
 \eeq
where $L$ is degree $2$ line bundle over $\phi_{1,1,2} =0$, the union  of two rational curves $y_3 =0$ and $y_3 {-}y_1y_2= 0$ embedded in 
$\mathbb{WP}_{1,1,2}$ 
and 
$$\mathcal M_{\phi_{{1,1,2}_0}} ={\rm Spec}(A_{ \phi_{{1,1,2}_0}} )$$ 
is the affine cone  in $\mathbb C^3$ {which is nothing but the 2-vertex $\mathbb V_2.$}
\endproof

In section \ref{suse:qGS} we provide a quantisation of these two degenerate cases and calculate the central elements - the quantisation of the full non-degenerate case can be found in \cite{EtGinzb}.

{\begin{remark}
In the next section, we will  discuss the appearance of the same degenerated Jacobian Poisson structures as a result of semi--classical limit of the quantum generalized Sklyanin-Painlev\'e 
``marginal" deformation of a quiver theory prepotential.
\end{remark}
}

Note that in the weighted projective space, there are many different homogeneous polynomials $\phi$ of degree $4$ that define the same quotient by the Jacobian ideal. For example
\begin{equation}\label{eq:cas112}
\bigl\{\tilde \phi_{1,1,2} = \tau_1 \tilde x_1\tilde x_2\tilde x_3 + \frac{1}{3}(\tilde x_3^2 + \tilde x_3\tilde x_1^2+ \tilde x_1\tilde x_2^3) =0\bigr\} \subset \mathbb{WP}_{1,1,2},
\end{equation}
defines the same algebra as $ \phi_{1,1,2}$ and
\begin{equation}\label{eq:cas231}\bigl\{\tilde \phi_{2,1,3}  = \tau_2 \tilde x_1\tilde x_2\tilde x_3+ \frac{1}{3} (\tilde x_1^3 + \tilde x_2^3\tilde x_3 + \tilde x_3^2)  =0\bigr\} \subset \mathbb{WP}_{2,1,3},
\end{equation}
defines the same algebra as $ \phi_{2,1,3}$.

In \cite{OdRubpol}, A. Odesskii and the  third author described two {\it non-rational} Poisson morphisms between the Poisson algebra of Jacobian type associated with the
 Hesse cubic \eqref{eq:poi-tau} in the variables  and the two homogeneous polynomials $\phi_{1,1,2}$:
\begin{equation}\label{eq:ref nalg1}
\tilde x_1= x_1^{\frac{3}{2}} x_3^{-\frac{3}{4}},\quad \tilde x_2 =x_2  x_3^{\frac{1}{4}}x_1^{-\frac{1}{2}}, \quad y_3 = x_3^{\frac{1}{2}}
\end{equation} 
 and   $  \phi_{2,1,3}  $:
 \begin{equation}\label{eq:ref nalg2}
\tilde x_1= x_1,\quad \tilde x_2 = x_2x_3^{-\frac{1}{2}}, \quad \tilde x_3 = x_3^{\frac{3}{2}}.
\end{equation} 

As  discussed in \cite{OdRubpol}, the non-rational Poisson morphisms  \eqref{eq:ref nalg1}, \eqref{eq:ref nalg2} have their origin in the Calabi-Yau mirror symmetry dualities (\cite{GrPlessRoan})
and the question of  their ``quantum" interpretation  was posed. Because by rescaling $\tilde x_1,\tilde x_2,\tilde x_3$ in exactly the same way as $x_1,x_2,x_3$ one can produce the same rational limits $ \phi_{{2,1,3}_0}$, $\phi_{{1,1,2}_0}$,
the quantisation produced in  subsection \ref{suse:qGS}   gives a partial  answer to this question by providing a quantisation for some rational limits of $\tilde\phi_{1,1,2}$ and $\tilde\phi_{2,1,3}$.

\subsubsection{\bf Degnerate Sklyanin algebras with three generators}

The Sklyanin algebra $Q_3({\mathcal E}, a,b,c)$ has the following {\it degeneration locus}  
$$
\mathcal D = \lbrace{(1,0,0);(0,0,1);(0,0,1)\rbrace}\sqcup\lbrace{(a,b,c)\mid a^3=b^3=c^3\rbrace}.
$$
Following \cite{Smith}, we call  \emph{degenerate Sklyanin algebra} the algebra $Q_3({\mathcal E, a,b,c})$  with $(a,b,c)\in D$.

It was proven by P. Smith that such a degenerate Sklyanin algebra is isomorphic to
 $\C\langle u,v,w\rangle /J$ where the ideal $J$ is $J=\langle u^2=v^2=w^2=0\rangle$ if $a=b$, and 
 $J=\langle uv=vw=wu=0\rangle $ if $a\neq b$. In the semiclassical limit the latter case corresponds to $\phi=u v w$, which is  the decorated  character variety of $\pi_1(\P^1\setminus\lbrace{z_1,z_2,z_3\rbrace})$ \cite{ChMR}.

\begin{remark}The latter model has a  quiver representation  with potential $Q= uu^* +vv^* + ww^* - uvw - wvu$  
\cite{BridgeSmith}. \end{remark}

\section{Non-commutative cubics and QFT deformations}\label{se:NC-QFT} 

There is an  interesting similarity  between the formulae for the quantum potential $\Phi$ defined in \eqref{eq:EGprep-gen}  and the non-commutative potentials describing the {\it marginal and relevant} deformations of the $N=4$ super Yang-Mills (SYM)  theory in four dimensions  with gauge group $U(n)$ (see \cite{DBerLeigh} for a physical background)
This theory is written in terms of  the $N=1$ 
SYM theory with three adjoint chiral super-fields $X_1,X_2, X_3$ coupled by the potential:
$$\Phi_{smooth} = g{\rm Tr}([X_1,X_2]X_3)$$ 
with coupling constant $g$, where, following the physics literature ${\rm Tr}$ denotes the map $A\to A_\natural$. From now on we drop ${\rm Tr}$, i.e. we denote potentials and their images in  $A_\natural$ with the same symbol.  
  
The moduli space of supersymmetric gauge theories is an important and rather  well-studied object (a mathematical account  of this theory can be found in the recent paper of C. Walton \cite{ChW1}). The marginal deformations, which preserve some conformal symmetry, of the $N=4$ Superconformal Field Theory have many interesting applications. In particular, within the framework of   the AdS/CFT correspondence, they have a nice Supergravity dual descriptions. 

If one chooses to preserve $N=1$ Super Conformal Field Theory then the moduli space of the marginal deformations  is given by the potential:
\beq\label{ref:superpotmarg}
\Phi_{\rm marg} =  X_1X_2X_3 - qX_2X_1X_3+ \frac{1}{3}\lambda (X_1^3+X_2^3+X_3^3).
\eeq

Another important class of  deformations is provided by {\it relevant deformations} which describes the theory away from the Ultra-Violet  conformal fixed point:
$$\Phi_{\rm rel} = m_1 X_1^2+ m_2  (X_2^2 +  X_3^2)  + \sum_k d_k X_k.$$

The structure of the vacua of $D$-brane gauge theories relates to Non-Commutative Geometry via the potentials  ${\Phi}_{phys}$ by so called $F-term$ constraints:
\beq\label{vacalg}
\frac{\partial {\Phi}_{phys}}{\partial X_k} = 0, \quad k=1,2,3
\eeq
where ${\Phi}_{phys}= \Phi_{\rm marg} + \Phi_{\rm rel}$. This gives rise to the following non homogeneous relations:
\beq\label{deformvac}
\left\{
\begin{array}{cc}
X_1X_2 - qX_2X_1 = & -\Lambda X_3^2 - m_2X_3 - d_3 \\
X_2X_3 - qX_3X_2 = & -\Lambda X_1^2 - m_1X_1 - d_1\\
X_3X_1 - qX_1X_3 = & -\Lambda X_2^2 -m_2X_2 - d_2
\end{array}
\right.
\eeq
This algebra is a particular case of the algebra $A^{{q}}_{ {\Phi_{EG} + \Psi_{EG}}}$ studied in subsection \ref{suse:EG-gen} for $\alpha=\beta=\gamma=-\Lambda$,  $a_1=m_1$, $b_1=c_1=m_2$, $a_2=d_1$, $b_2=d_2$, $c_2=d_3$, or in other words, of the general algebra 
$\mathcal A^q$ introduced in subsection \ref{se:pain-pot}.

\begin{remark} We precise how this deformation algebra relates to previously studied:
\begin{itemize}
\item If {$\Lambda= 0$} and $m_1 =m_2=- \frac{1}{2},\quad e_i =0, \quad i=1,2,3$ then we have the potential of \eqref{ref:CYodess}
and this algebra coincides with the Odesskii degeneration of Sklyanin algebra in subsection \ref{OdSklDeg};
\item We see that for ${\Lambda}=0$, $m_1=m_2=\sqrt{q}(q^{-1}-q) $   and $d_i=(1-q)\Omega_i^{(VI)}$   the Poisson algebra \eqref{deformvac} has the form \eqref{eq:q-comm}.
\item If $\Lambda=\frac{c}{a}$ and $q=b/a$, then this Poisson algebra coincides with a deformation of the quadratic
Sklyanin Poisson algebra with three generators $q_{3}(\mathcal E)$. The latter can be obtained from this deformation  by setting the massess to zero: $m_1=m_2 =0.$
\end{itemize} 
\end{remark}

\subsection{Semi-classical limits.}
We now take the semi-classical limit of  \eqref{deformvac}  and compare it with the cubic surfaces $\mathcal M_\phi:= {{\rm Spec}}(\mathbb C[x_1,x_2,x_3]\slash\langle\phi=0\rangle)$, for $\phi$ in table 1. These cubics  are endowed with the natural Poisson bracket \eqref{eq:nambu}. By the correspondence principle
$$ \lim_{q\to 1}\frac{[X_1,X_2]}{1-q} = \{x_1,x_2\},
$$
and, applying the algebra relations 
$$
[X_1,X_2]=(q-1) X_2 X_1 -\Lambda X_3^2 - m_2X_3 - d_3
$$
so that 
$$\{x_1,x_2\} = x_1x_2  -   \lim_{q\to 1}\frac{{\Lambda}}{1-q}X_3^2 -\lim_{q\to 1}\frac{{m_2 X_3}}{1-q} + \lim_{q\to 1}\frac{d_3}{1-q},$$
and similarly
$$\{x_2,x_3\} = x_2x_3  -   \lim_{q\to 1}\frac{{\Lambda}}{1-q}X_1^2 -\lim_{q\to 1}\frac{{m_1 X_1}}{1-q} +\lim_{q\to 1}\frac{d_1}{1-q},$$
$$\{x_3,x_1\} = x_1x_3  -   \lim_{q\to 1}\frac{{\Lambda}}{1-q}X_2^2 -\lim_{q\to 1}\frac{{m_2 X_2}}{1-q} + \lim_{q\to 1}\frac{d_2}{1-q}.$$

By a slight abuse of notation, we denote the classical masses again by $m_1,m_2$, the classical limit of $\Lambda$ by $\lambda$ and  put $\delta_i=\lim_{q\to 1}\frac{d_i}{1-q}$, so that the Casimir function for this Poisson algebra is
$$\phi_{\rm cl,tot} (x_1,x_2,x_3) = x_1x_2x_3 -m_1 x_1^2 -m_2(x_2^2 +x_3^2) - \frac{\lambda}{3}(x_1^3 + x_2^3 + x_3^3) +\delta_1x_1+\delta_2x_2+\delta_3x_3. 
$$
We see that the corresponding Poisson algebras include all interesting families of quadratic-linear-constant Poisson brackets in $\C[x_1,x_2,x_3]$ and, in particular, for $\lambda =0$,
the family coincides with the Poisson structure on the Painlev\'e monodromy data cubics.  At the same time, by 
neglecting the terms of degree $<3$ in $\phi_{\rm cl,tot}$ we obtain 
 \beq\label{quadrdef Sklyanin}
\phi_{\rm cl,marg} (x_1,x_2,x_3) = x_1x_2x_3 -m_1 x_1^2 -m_2(x_2^2 +x_3^2) - \frac{{\lambda}}{3}(x_1^3 + x_2^3 + x_3^3),
\eeq
that is a perturbation of the classical Sklyanin algebra $q_{3,1}(\mathcal E)$ (see section \ref{skl3}).

\subsection{Degeneration of  quadratically perturbed  $q_3$-Sklyanin brackets and {vertex $\mathbb V_{1,2}$ varieties.}}\label{SklGrZieb}

Consider the special case of  \eqref{quadrdef Sklyanin} with $m_2 =0$ and $\lambda = \frac{3}{m_1^3}:$
\beq\label{quadrdef Sklyanin1}
\phi_{\rm cl,1} (x_1,x_2,x_3) = x_1x_2x_3 -m_1 x_1^2 - \frac{1}{m_1^3}(x_1^3 + x_2^3 + x_3^3). 
\eeq
This is an example of a central element for the classical Sklyanin algebra perturbed  by the quadratic term $m_1 x_1^2$.
{The following proposition is a rephrasing of Proposition \ref{prop:degen} with the ``infinite mass limit" interpretation.
 \begin{prop}
There are two  ``mass re-scalings"  of \eqref{quadrdef Sklyanin1} such that in the infinite mass limit $m_1\to \infty$, $\phi_{\rm cl,1}$ degenerates to the potential $\phi_{{2,1,3}_0}$ of the vertex $\mathbb V_1$ Poisson algebra in the first rescaling, and to the potential  $\phi_{{1,1,2}_0}$  of the Poisson algebra of the vertex $\mathbb V_2$ in the second rescaling.
\end{prop} 
}
\proof We introduce the coordinates $y_i, i=1,2,3$ connected  to $x_1,x_2,x_3$ by the
following relations:
$$x_1 = \frac{y_1}{\sqrt{m_1}},\quad x_2 = \frac{y_2}{\sqrt{m_1}},\quad  x_3 = m_1 y_3,$$
so that the Casimir now reads as
\beq\label{quadrdef Sklyanin2}
\phi_{\rm cl, 2} (y_1,y_2,y_3) = y_1y_2y_3  - y_1^2 - y_3^3 + \frac{(y_1^3 + y_2^3)}{\sqrt{m_1^9}}. 
\eeq

In the infinite mass limit $m_1\to \infty$, \eqref{quadrdef Sklyanin2} goes evidently to 
\beq\label{ref:quadrdef Sklyanin3}
\phi_{\rm cl, 3} (y_1,y_2,y_3) = y_1y_2y_3  - y_1^2 - y_3^3
\eeq
Note that up to permutations of $y_1,y_2,y_3$, $\phi_{\rm cl, 3}$ is the same as $\tilde\phi_{213_0}$, and therefore, as 
 discussed at the end of subsection \ref{suse:E7-E8}, 
the cubic surface ${\mathbb V_1 = }\mathcal M_{\phi_{\rm cl, 3}}\subset \C^3$ given by $\phi_{\rm cl, 3} (y_1,y_2,y_3) = y_1y_2y_3  - y_1^2 - y_3^3$  can be considered as an affine cone 
over a singular {genus one rational} curve  $\mathcal E_{sing}\subset \mathbb W\P(3,1,2)$. Its coordinate ring 
$$ 
\C[\mathcal M_{\phi_{\rm cl, 3}}]=\C[y_1, y_2,y_3]/( y_1y_2y_3  - y_1^2 - y_3^3)
$$ 
is isomorphic to the ring of sections $\oplus_{k\geq 0}H^0(\mathcal E_{sing},{\mathcal O}(k))$ of a degree $1$ line bundle ${\mathcal O}(1))$ on the nodal rational curve $\mathcal E_{{sing}}$
of arithmetic genus 1 (see \cite{GrossHK} ch.5). This cone is parametrized by toric theta-functions
$\vartheta_i$, $ i=1,2,3$ satisfying the relation
$$
\vartheta_1\vartheta_2\vartheta_3  = \vartheta_1^2 + \vartheta_3^3
$$
(see Theorem 2.34 of  \cite{GrossHK}).

Now we come back to the  Poisson algebra corresponding to \eqref{quadrdef Sklyanin1}:
$$\{x_1, x_2\} = -\frac{3x_3^2}{m_1^3} + x_1x_2;\quad \{x_2, x_3\} = -\frac{3x_1^2}{m_1^3} -2m_1 x_1 + x_2x_3;\quad \{x_3, x_1\} = \frac{3x_2^2}{m_1^3} + x_3x_1$$
which will be written in the degenerated coordinates $y_i, i=1,2,3$ as

\beq\label{ref:nodcub1}\{y_1,y_2\} = -3y_3^2 + y_1y_2;\quad \{y_2, y_3\} = - \frac{3y_1^2}{\sqrt{m_1^9}} -2y_1 + y_2y_3;\quad \{y_3, y_1\} =   \frac{3y_2^2}{\sqrt{m_1^9}} + y_3y_1.\eeq
From this,  in the infinite mass limit  we obtain once again (compare with \eqref{Grosstheta_1}) a perturbed cluster Poisson structure:
\beq\label{ref:nodcub2}\{y_1,y_2\} = -3y_3^2 + y_1y_2;\quad \{y_2, y_3\} =  -2y_1 + y_2y_3;\quad \{y_3, y_1\} =   y_3y_1\eeq
 which  defines the Poisson algebra  structure on the coordinate ring of the affine cone over the  curve $\mathcal E_{sing}$.

If, instead, we introduce the coordinates $\tilde y_i, i=1,2,3$ connected  to $x_1,x_2,x_3$ by the
following relations:
$$x_1 = \frac{\tilde y_1}{\sqrt{m_1}},\quad x_2 = \tilde y_2,\quad  x_3 = \sqrt{m_1} \tilde y_3,$$
the Casimir now reads as
\beq\label{quadrdef Sklyanin3}
\phi_{\rm cl, 4} (\tilde y_1,\tilde y_2,\tilde y_3) = \tilde y_1\tilde y_2\tilde y_3  - \tilde y_1^2 - \frac{1}{m_1^3}(\frac{\tilde y_1^3}{\sqrt{m_1^3}}+ \tilde y_2^3 +m_1^{3/2} \tilde y_3^3).
\eeq

In the infinite mass limit $m_1\to \infty$, \eqref{quadrdef Sklyanin3} goes evidently to 
\beq\label{ref:quadrdef Sklyanin3}
\phi_{\rm cl, 5} (\tilde y_1,\tilde y_2,\tilde y_3) = \tilde y_1\tilde y_2\tilde y_3  - \tilde y_1^2 
\eeq
Note that (up to the change of variable $\tilde y_1=y_1$, $\tilde y_2=y_3$, $\tilde y_3=y_2$) $\phi_{\rm cl, 5}$ is the same as $\tilde\phi_{112_0}$, and therefore, as before,
 the cubic surface ${\mathbb V_2 =}\mathcal M_{\phi_{\rm cl, 4}}\subset \C^3$ is an affine cone 
over the singular curve  $\tilde{\mathcal E}_{sing}\subset \mathbb W\P(2,1,1)$. Its coordinate ring 
$$ 
\C[\mathcal M_{\phi_{\rm cl, 4}}]=\C[\tilde y_1, \tilde y_2,\tilde y_3]/( \tilde y_1\tilde y_2\tilde y_3  - \tilde y_1^2)
$$ 
is isomorphic to the ring of sections $\oplus_{k\geq 0}H^0(\tilde{\mathcal E}_{sing},{\mathcal O}(k))$ of a degree $2$ line bundle ${\mathcal O}(1))$ on the degenerated curve 
$\tilde{\mathcal E}_{{sing}}$ which is a union of conic and a line. This cone is parametrised by Gross-Siebert toric theta-functions
$\vartheta_i$, $ i=1,2,3$ satisfying the relation
$$
(\vartheta_1\vartheta_2 - \vartheta_3 )\vartheta_3 =0
$$
(see Proposition 40 of  \cite{Bouss}).

By writing the  Poisson algebra corresponding to \eqref{quadrdef Sklyanin1} in the new coordinates $\tilde y_1,\tilde y_2,\tilde y_3$ and taking 
 the infinite mass limit  we obtain once again  (compare  with \ref{Grosstheta_1}) a perturbed cluster Poisson structure:
\beq\label{ref:nodcub3}\{\tilde y_1,\tilde y_2\} = \tilde y_1\tilde y_2;\quad \{\tilde y_2, \tilde y_3\} =  -2\tilde y_1 + \tilde y_2\tilde y_3;\quad \{\tilde y_3, \tilde y_1\} =   \tilde y_3\tilde y_1\eeq
 which  defines the Poisson algebra  structure on the coordinate ring of the affine cone over the  curve $\tilde{\mathcal E}_{sing}$.
\endproof

\subsection{Quantisation of Gross-Siebert theta functions}\label{suse:qGS}

In \cite{Bouss}, P. Bousseau proposed a deformation quantisation for some Poisson algebra {$H^0({\mathcal X},{\mathcal O}_{\mathcal X})$} structures connected with mirror duals of
 Looijenga pairs $(Y,D)$ where $Y$ is a smooth projective surface and $D$ some singular anticanonical   divisor. {Here ${\mathcal X} \to S$ is a (possibly singular) family of affine Poisson varieties.}
As examples he computed the deformation  quantisation of function algebras on affine {$r-$vertex} varieties $\mathbb V_r$ where $r$ is the number of irreducible components of the cubic 
divisor $D$ {of the A-side of mirror correspondence}.
{As we have seen,} when $r=1$ the {1-vertex} variety $\mathbb V_1$ is exactly the affine cone of the nodal curve embedded in the weighted projective space ${\mathbb WP}_{2,1,3}:$ 
$$
{\mathbb V_1}=\mathcal M_{\tilde\phi_{{2,1,3}_0} }={\rm Spec}A_{\tilde \phi_{{2,1,3}_0}}
$$ 
where $A_{\tilde \phi_{2,1,3}}$ is given in \eqref{eqref: 123tilde}.
The Poisson algebra on ${\mathbb V_1}$ is given by the brackets \eqref{Grosstheta_1}.

The Proposition 41 in \cite{Bouss} states that the {quantum algebra $\mathbb V_1^{(q)}$ is generated by three generators $\hat Y_1, \hat Y_2$ and  $\hat Y_3$ with the} relations
$$
\left\{
\begin{array}{c}
\sqrt{\hat q}\hat Y_3\hat Y_1 -  \frac{1}{\sqrt{\hat q}}\hat Y_1\hat Y_3 =  0 \\
\sqrt{\hat q}\hat Y_2\hat Y_3 -  \frac{1}{\sqrt{\hat q}}\hat Y_3\hat Y_2 = (\hat q-\hat q^{-1})\hat Y_1\\
\sqrt{\hat q}\hat Y_1\hat Y_2 -  \frac{1}{\sqrt{\hat q}}\hat Y_2\hat Y_1 = (\hat q^{3/2}-\hat q^{-3/2}) \hat Y_3^2
\end{array}
\right.
$$
and that the {$0$--level of the} central element:
$$
{\hat {\Omega}}_{2,1,3}(\hat Y)  = \hat Y_2 \hat Y_3 \hat Y_1 - \hat q^{1/2} \hat Y_1^2  - \hat q \hat Y_3^3
$$
gives the quantisation of  the algebra ({\ref{ref:nodcub2}}) (equivalent to \eqref{Grosstheta_1}).

{As another example,}  P. Bousseau 
considered also a deformation quantisation of  the function algebra on the { 2-vertex} ${\mathbb V_2}$ related to the mirror dual of the Looijenga pair $(Y,D)$ where the divisor has two connected components, namely for $\mathbb V_2 =\mathcal M_{ \phi_{1,1,2}^0}$, where $ \phi_{{1,1,2}_0}$ is given in \eqref{eq:phi-112-0} and the Poisson algebra 
is the Jacobian algebra on $\C[ y_1, y_2,y_3]$ with the brackets 
\eqref{nodcub_3}. 
{The Proposition 40 in \cite{Bouss} describes the quantum algebra $\mathbb V_2^{(q)}$ in a similar fashion with three generators $\hat Y_1, \hat Y_2$ and  $\hat Y_3$ with the relations
$$
\left\{
\begin{array}{c}
\sqrt{\hat q}\hat Y_3\hat Y_1 -  \frac{1}{\sqrt{\hat q}}\hat Y_1\hat Y_3 =  0 \\
\sqrt{\hat q}\hat Y_2\hat Y_3 -  \frac{1}{\sqrt{\hat q}}\hat Y_3\hat Y_2 =  0 \\
\sqrt{\hat q}\hat Y_1\hat Y_2 -  \frac{1}{\sqrt{\hat q}}\hat Y_2\hat Y_1 = (\hat q^{1/2}-\hat q^{-1/2}) \hat Y_3^2
\end{array}
\right.
$$
and the {$0$--level of  the} central element:
$$
{\hat {\Omega}}_{1,1,2}(\hat Y)  = \hat Y_1 \hat Y_2 \hat Y_3 - \hat q^{1/2} \hat Y_3^2
$$
gives the quantisation of \eqref{nodcub_3}. 
} 

{The deformational quantisation of the affine Poisson family ${\mathcal X} \to S$ elaborated by P. Bousseau was based on Gross--Haking--Keel approach to toric mirror conjecture, tropicalization and 
 quantum scattering diagrams and broken lines. On the other hand, if $Y$ in the Looijenga pair $(Y,D)$ is a del Pezzo surface of degree 1, 2 or 3 and the divisor $D$ a nodal cubic, Etingof, Oblomkov and Rains constructed the corresponding quantum algebras as spherical sub--algebras of their generalised DAHA (\cite{EOR}). In Subsection \ref{suse:EG-UP},  we obtained the generalised Sklyanin--Painlev\'e algebra as a degeneration of the $\widetilde E_6$ spherical sub--algebra and fit it within the Etingof--Ginzburg quantisation scheme (see Theorem \ref{th:EG-UZ}).} 
  
A natural question  posed by Bousseau is to compare  his  deformation quantisations with the scheme of quantisation proposed by Etingof and Ginzburg.  {In the next theorem we answer this question in the case of degenerated affine del Pezzo surfaces (the case  ${\mathbb V_1}$ and  ${\mathbb V_2}$).}

\begin{theorem}\label{compar}
The deformation  quantisations of the affine Poisson structures on $\mathbb V_{1,2}$ obtained in \cite{Bouss} coincide (after a proper rescaling) with the appropriate degenerations of the quantum Sklyanin-Painlev\'e algebras defined by relations \eqref{eq:q-comm-GSP}.
\end{theorem}
\proof
We start by observing that  the quantum  algebra corresponding to \eqref{ref:quadrdef Sklyanin3} is a degenerate case of the Calabi-Yau algebra 
$\C\langle X_1,X_2,X_3\rangle /J_{\Phi_{phys}}$ {with the potential } \eqref{vacalg}. Indeed, by analogy with the classical case, we introduce the coordinates $Y_i, i=1,2,3$ connected  to $X_1,X_2,X_3$ by the
following relations:
$$
X_1 = \frac{Y_1}{\sqrt{m_1}},\quad  X_2 =  \frac{Y_2}{\sqrt{m_1}},\quad X_3 = m_1 Y_3,
$$
to obtain
\begin{equation}\label{eq-phi-y-cone}
\begin{split}
\Phi_{phys} = Y_1Y_2Y_3 - qY_2Y_1Y_3 + \frac{\Lambda}{3}\left(m_1^3 Y_3^3 +\frac{Y_1^3+Y_2^3}{\sqrt{m_1^3}}\right)+\\
\qquad+ \frac{1}{2}Y_1^2 + \frac{m_2}{m_1} Y_2^2+ m_1 m_2 Y_3^2 +e_1Y_1 + e_2Y_3+ e_3Y_2,\end{split}
\end{equation}
which is by our discussion a  PBW non-homogeneous deformation of the Koszul generalised Sklyanin algebra. By putting $\Lambda=m_1^{-3}$, $m_2=0$ and $e_1=e_2=e_3=0$, we obtain
$$
{\Phi}_{m_1} = Y_1Y_2Y_3 - qY_2Y_1Y_3 + \frac{1}{3}\left( Y_3^3 +\frac{Y_1^3+Y_2^3}{m_1^3\sqrt{m_1^3}}\right)+ \frac{1}{2}Y_1^2,
$$
and in the limit $m_1\to\infty$ we obtain
$$
{\Phi}_{\infty}(Y)  = Y_1Y_2Y_3 - qY_2Y_1Y_3  + \frac{1}{3} Y_3^3 + \frac{1}{2}Y_1^2
$$
and the corresponding quantum algebra $\C\langle Y_1,Y_2,Y_3\rangle /J_{{\Phi}_{\infty}}$ has relations
\beq\label{eq-deformvacdeg}
\left\{
\begin{array}{cc}
Y_3Y_1 - qY_1Y_3 = & 0 \\
Y_2Y_3 - qY_3Y_2 = & Y_1\\
Y_1Y_2 - qY_2Y_1 = & Y_3^2
\end{array}
\right.
\eeq
This algebra has central element
\beq\label{eq-centrdeg}
{\Omega}^{m_1}_{0}(Y)  = Y_3Y_2 Y_1+ \frac{q}{q^2-1} Y_1^2 +  \frac{q^2}{q^3-1} Y_3^3
\eeq
and quantises  the coordinate ring of the cone over the nodal rational genus 1 curve or the coordinate ring of the affine surface \eqref{ref:quadrdef Sklyanin3}. 
But these are the same as \eqref{eq-deformvacdeg} and \eqref{eq-centrdeg} by setting 
$$
q=\frac{1}{\hat q}, \quad
\quad Y_1= \frac{(1=\hat q)(\hat q -1)^3(1+\hat q+\hat q^2)^2}{\hat q^{\frac{1}{4}}}
\hat Y_1,
\quad \hat Y_2 = \hat q^{\frac{7}{4}} Y_2,
\quad
\hat Y_3= (1-\hat q^2-\hat q^3 +\hat q^5) Y_3.
$$
We can degenerate the algebra \eqref{eq-deformvacdeg} further by rescaling the variables $Y_1,Y_2,Y_3$ and taking different limits. Namely, setting
$$
Y_1\to \epsilon_1  Y_1, \quad Y_2\to \epsilon_2  Y_2, \quad Y_3\to \epsilon_3 Y_3, 
$$
we obtain
\beq\label{eq-deformvacdegeps}
\left\{
\begin{array}{cc}
Y_3Y_1 - qY_1Y_3 = & 0 \\
Y_2Y_3 - qY_3Y_2 = &\frac{\epsilon_1}{\epsilon_2 \epsilon_3} Y_1\\
Y_1Y_2 - qY_2Y_1 = &\frac{\epsilon_3^2}{\epsilon_1 \epsilon_2} Y_3^2
\end{array}
\right.
\eeq
with central element
\beq\label{eq-centrdegeps}
{\Omega}^{m_1}_{0}(Y)  = Y_3Y_2 Y_1+ \frac{q}{q^2-1}\frac{\epsilon_1}{\epsilon_2 \epsilon_3} Y_1^2 +  \frac{q^2}{q^3-1} {\epsilon_3^2}{\epsilon_1 \epsilon_2}  Y_3^3
\eeq
Imposing $\epsilon_2=1$, $\epsilon_1=\epsilon_3^2$, $\epsilon_3$, in the limit $\epsilon_3\to 0$, we obtain
\beq\label{eq-deformvacdeg2}
\left\{
\begin{array}{cc}
Y_1Y_3 - qY_3Y_1 = & 0 \\
Y_2Y_3 - qY_3Y_2 = & 0\\
Y_2Y_1 - qY_1Y_2 = & Y_3^2
\end{array}
\right.
\eeq
and the central  element is given by 
\beq\label{eq-centrdeg2}
{\Omega}^{\infty}(Y)  = Y_2Y_3 Y_1+ \frac{1}{q^2-1} Y_3^2. 
\eeq
\endproof

Observe that by choosing different values and limits of $\epsilon_1,\epsilon_2,\epsilon_3$ in \eqref{eq-deformvacdegeps}, we can recognise the algebras given by the super-pontentials of non-commutative Painlev\'e cubics (PIV and PII) to which the next two subsections are dedicated.

\subsubsection{One non-zero mass and Painlev\'e IV}
We consider  the deformation provided by addition a single mass term to $\Phi_{smooth}$.
The corresponding potential (4.2 of  \cite{DBerLeigh}) reads (up to {symmetric group $\Sigma_3$-}action):
\beq\label{superpotsinglmass}
\Phi_{1m} = X_1X_2X_3 -qX_2X_1X_3 - \frac{m}{2}X_1^2.
\eeq

The corresponding ideal is defined by 
\beq\label{ncPIV}X_1X_2 - qX_2X_1 = 0 ;\quad X_2X_3 - qX_3X_2 =  mX_1;\quad X_3X_1 - qX_1X_3= 0\eeq
Taking the Poisson limit $q\to 1$ one gets the cubic Casimir :
$$\phi_{\rm cl,PIV} (x_1,x_2,x_3) = x_1x_2x_3 - \frac{m}{2}x_1^2 . $$

Once again, to link with some of our Painlev\'e cubics ( in the single mass case it will be the PIV cubic) we need to add the linear terms:
\beq\label{superpotsinglmasslin}
\Phi_{1,m} = X_1X_2X_3 - qX_2X_1X_3 - \frac{m}{2}X_1^2 + d_1X_1 +d_2X_2 +d_2X_3.
\eeq
Taking $d_2=d_3$ one gets
\beq\label{ncPIV-l}X_1X_2 - qX_2X_1 = d_2 ;\quad X_2X_3 - qX_3X_2 =  mX_1 + d_1;\quad X_3X_1 - qX_1X_3 = d_2\eeq
and the cubic Casimir ($q\neq \pm 1$)
$$\Phi_{\rm PIV}  = X_1X_2X_3 -qX_2X_1X_3 - \frac{m}{2}X_1^2 + \frac{1}{1-q}(d_1 X_1 +d_2(X_2+X_3)). $$
corresponds to the $PIV$ case in the table of cubics 1.
\subsubsection{NC Painlev\'e II}
Consider the potential 
\beq\label{superpotPII}
\Phi_{\rm PII} = X_1X_2X_3 - q X_2X_1X_3 + (q-1)(X_1 + \Omega_2 X_2 +X_3)
\eeq
The algebraic relations corresponding to the Jacobian ideal are
\beq\label{ncPII}X_1X_2 - qX_2X_1 =  (q-1) ;\quad X_2X_3 - qX_3X_2 =  (q-1) ;\quad X_3X_1 - qX_1X_3 = (q-1)\Omega_2\eeq
and the Poisson limit gives the Casimir  cubics for the Miwa-Jimbo Painlev\'e II cases:
$$\phi_{\rm cl,PII} = -x_1x_2x_3 +x_1+ \omega_2x_2 + x_3$$
The authors of \cite{DBerLeigh} argue that, in the  framework of study of ``orbifold singularities", one  should take the relation on the
moduli space
\beq\label{ChebDef}
x_1x_2x_3 - (e_1^n x_1 + e_2^n x_2+ e_3^n x_3) + 2(e_1e_2e_3)^{n/2}T_n(-\frac{w}{2(e_1e_2e_3)^{1/2}}) =0,
\eeq
where $T_n(w) = \cos(n \arccos w)$ is the $n$- th Chebyshev polynomial. 

Taking $e_1=e_2=e_3=\exp(\frac{i\pi}{n})$ and $n=1$ (which means $T_1(w) =w$) we have the expression
$$
x_1x_2x_3+ x_1 + x_2+ x_3  - w =0,
$$ so  the Miwa-Jimbo Painlev\'e Casimir cubic can be considered as  the $n=1$ member of the family \eqref{ChebDef}.

\subsection{Le Bryun -Witten algebras}\label{suse:LBW}
Our final  remark is that the  generalised Sklyanin-Painlev\'e algebra with potential $\Phi^\gamma$
gives an example of the {\it conformal $\mathfrak{sl}_2$- enveloping algebra} $U_{abc}(\mathfrak{sl}_2)$ (\cite{LB}). It corresponds to the choice of
the parameters $a=c=q;\quad m_1=m_2=1;\quad e_1=e_3=0, e_2=-1;-\gamma = b:$
\beq\label{CUabc}
\left\{
\begin{array}{cc}
X_1X_2 - qX_2X_1 = & X_3 \\
X_2X_3 - qX_3X_2 = & X_1 \\
X_3X_1 - qX_1X_3 = & -\gamma X_2^2 + X_2 +1
\end{array}
\right.
\eeq
This algebra corresponds to the generalised Sklyanin with $\beta=\gamma=0$, $a=b=c=q$, case (2) of Theorem \ref{gensklhom}

The central element is
$$
\Omega_{LBW}=(q^2-1) X_3 X_2 X_1 - \gamma\frac{1+q}{q(1+q+q^2)} X_2^3 + q X_1^2 +\frac{1}{q} X_2^2+ q X_3^2.
$$

It was proved by Le Bruyn in (\cite{LB})  that the conformal $\mathfrak{sl}_2$ enevloping algebras are {\it Auslander regular} and have the {\it Cohen-Macaulay property} as finitely generated (left) filtered rings.
We observe now that, following the results of Artin, Tate and Van den Bergh (\cite{ATV}), one can construct a cubic divisor $C \hookrightarrow \P^2$ for any three-dimensional Auslander-regular algebra and the
algebra is defined by the divisor and an automorphism $\sigma : C\to C.$

This divisor is defined by the  equation 
\beq\label{Cdiv}
[C]={\rm det}\left(\begin{array}{ccc}\gamma X_2& -qX_1& X_3\\X_1&0&-qX_2\\-qX_3& X_2&0\\
\end{array}\right)=0,
\eeq
where the determinant is calculated quantically as follows:
$$-\gamma q X_2^3  + (q^3  -1)X_1X_2X_3 =0$$
and defines a conic  ($-\gamma q X_2^2  + (q^3  -1)X_1X_3 =0$) and a line $X_2=0$. 
The automorphism $\sigma$ is given on the line by $\sigma (X_1:0:X_3) = (X_1:0:qX_3)$ and on the conic by $\sigma(X_1:X_2:X_3) = (qX:X_2:q^{-1}X_3).$
Thus, we see that if $\gamma \equiv 0, q^3\neq 1$ then the divisor gives the triangular configuration $X_1X_2X_3=0$ and if $q^3=1$ then the divisor degeneratres in a triple line.

\section{del Pezzo of degree $3$ and open problems}\label{se:conclusion}

In this section, we summarise our results concerning del Pezzo of degree $3$ and their quantisation in three tables and highlight some open problems for the future.

Let us start by describing Table 4. 

The first column contains a list of double affine Hecke algebras. The elliptic DAHA of type $\widetilde E_6$ is due to Rains, \cite{R1}, while the GDAHA of type $E_6^{(1)}$ is due to \cite{EOR}. The abbreviation ``Deg. GDAHA"  corresponds to some Whittaker degenerations of the $E_6^{(1)}$ GDAHA \cite{Cher1,M4}, the $\check CC_1$ DAHA is due to Cherednik \cite{Cher,Sa,NS}, while the abbreviation ``Deg. DAHA"  correspond to the algebras obtained in \cite{M2} by Whittaker degeneration.

The second column is the  polynomial $\phi$ such that $\mathcal M_\phi$ is the  center of the corresponding (elliptic or generalised or degenerate) DAHA for $q=1$: for the cases of Elliptic DAHA, this was conjectured in \cite{EOR}, for GDAHA it was proved in \cite{EOR}, for the $\check CC_1$ in \cite{Obl} and all other cases in \cite{M2,M4}.

As discussed is Section \ref{se:GHK-theta}, the  projective completion $\overline{\mathcal M}_\phi$ is a del Pezzo of degree $3$ with divisor $D_\infty$ - this is specified in the third column of table 4.

The table is split vertically by a double line - the whole right side of the table is due to H. Sakai \cite{sakai}, and we have used his notation here.  Before explaining what this double line represents, let us recall the definition of an {\it Okamoto pair} $(X,\Delta)$: this is a pair $(X,\Delta)$  where $X$ is a generalised Halphen surface, namely the blow up of $9$ points in $\mathbb P^2$ in non generic position, and $\Delta$ is a divisor that tells us the position of such $9$ points. Note that  $\Delta$ has the same configuration as a degenerate elliptic curve in the classification by Kodaira-Neron. In other words, the $9$ non generic points lie at the intersection between $\Delta$ and a generic elliptic curve in $\mathbb P^2$. The generalised Halphen surfaces are uniquely determined by their divisor $\Delta$ listed in the fifth column.  

Some of these divisors have multiple points on them. Starting from the fifth row, at the intersection of lines we always have a multiple point, this is denoted by an empty circle. The order of this point can be calculated by removing from the number $9$ the order of all other points. The single bullet points mean simple points, the bullets with a circle and a number next to them mean multiple points with the order specified by the number. 
In the last column we show the blow up of such divisor $\Delta$ at the multiple points.  

Sakai labels the generalised Halphen surface according to the affine Weyl group corresponding to the intersection matrix of the divisor. Note that in the case $A_0^{(1)}$ Sakai uses two notations according to the divisor, no star means $\Delta$ is a smooth elliptic curve, one star means  $\Delta$  is  rational curve with a node. These labels are given in the fourth column.

The first line of the table corresponds to the elliptic Painlev\'e equation, the next three lines to the multiplicative or $q$-difference Painlev\'e equations and the last eight rows correspond to the Painlev\'e differential equations. There are also additive difference Painlev\'e equations, which we give in Table 6 because the corresponding quantum algebra is not Calabi-Yau \cite{MoSm}. Finally, there are also high dimension multiplicative difference Painlev\'e equations the quantum description of which is postponed \cite{M4}.

In the case of the Painlev\'e differential equations, the left and right sides  of the table are related by the so-called
Riemann-Hilbert correspondence  - this was proved by several authors, a nice unified approach can be found in \cite{SvdP}. The basic idea is that the Okamoto pair corresponds to the space of initial conditions of the given equation, while the Looijenga pair corresponds to the monodromy manifold. 
 
Okamoto's theory of initial value spaces \cite{O79}, developed by Sakai \cite{sakai}, provides a beautiful unification of differential and discrete equations. Whether differential or discrete, initial values for any nonlinear equation, can be regular (meaning the solution will be analytic around the initial point) or can be unbounded (reflecting the existence of a singularity at the initial point). Okamoto compactified this space to the complex projective plane and showed that any subsequent indeterminacy can be removed by resolving the base points through blowup techniques from algebraic geometry. It is a miraculous fact that nine blowups leads to a regularisation of the whole space for all differential and discrete Painlev\'e equations.
For discrete Painlev\'e equations, there is no satisfactory concept of monodromy manifold - it is true that each additive discrete Painlev\'e equation comes from the Backl\"und transformations of one of the differential ones, so that one could use the monodromy manifold associated to the latter, however without a direct isomonodromic approach, interesting dynamical behaviour may be lost. Moreover for the multiplicative discrete Painlev\'e equations, a notion of monodromy manifold is completely missing. This leads us to

\begin{conj}
For the elliptic and multiplicative/additive discrete Painlev\'e equations, the Riemann Hilbert correspondence assigns to the generalised Halphen surface in Table 4 the corresponding Looijenga pair.
\end{conj}

Intuitively speaking, evidence for this conjecture is provided by the fact that the polynomials defining the divisors $D_\infty$ in the first four lines of Table 4 are the same as those defining the corresponding Halphen divisors  $\Delta$ - i.e. $D_\infty=\Delta$ for the first four lines in the table.

We list the quantum results in Table 5 and 6. All the quantum algebras in Table 5 are specialisations of the generalised Sklyanin-Painlev\'e algebra introduced in subsection \ref{se:pain-pot}.

We conclude by mentioning the relation between the quantum algebras in Table 5 and the matrix generalisations of the Painlev\'e equations. Building upon work by Retakh and the third author   \cite{RR}, in \cite{BCR} a set of non-commutative relations which are non-commutative analogues of monodromy data relations for the Painlev\'e II equation was constructed. 
The interesting feature of these non-commutative relations is that by taking the scalar degeneration of the non-commutative operator $q$, one obtains our quantum  Painlev\'e II monodromy variety. This observation opens the possibility of relating higher rank Elliptic/Generalised DAHA to the theory of matrix Painlev\'e equations.
\newpage

%\psset{unit=0.8}
\begin{pspicture}(-7,-12)(6,8)
\psframe(-7,-11)(6,8)

%vertical lines

\psline(-5.8,-11)(-5.8,8)
\psline(-2.6,-11)(-2.6,8)
\psline(-0.8,-11)(-0.8,8)
\psline(-0.75,-11)(-0.75,8)
%\psline(.4,-11)(.4,8)
\psline(.4,-11)(.4,8)
\psline(3,-11)(3,8)

%horizontal lines

\psline(-7,7)(6,7)
\psline(-7,5.5)(6,5.5)
\psline(-7,4)(6,4)

\psline(-5.8,2.5)(6,2.5)

\psline(-7,1)(6,1)

\psline(-7,-.5)(6,-.5)

\psline(-5.8,-2)(6,-2)
\psline(-5.8,-3.5)(6,-3.5)
\psline(-5.8,-5)(6,-5)
\psline(-5.8,-6.5)(6,-6.5)

\psline(-5.8,-8)(6,-8)
\psline(-5.8,-9.5)(6,-9.5)

% first row
\rput(-6.4,7.53){\tiny{DAHA}}
\rput(-4.3,7.7){\tiny{Center $\phi$}}
\rput(-4.3,7.43){\tiny{for $q=1$}}
\rput(-1.9,7.7){\tiny{del Pezzo}}
\rput(-1.8,7.3){\tiny{divisor $D_\infty$}}
\rput(-0.2,7.8){\tiny{Gen.}}
\rput(-0.2,7.53){\tiny{Halphen}}
\rput(-0.2,7.3){\tiny{surface}}
\rput(1.2,7.6){\tiny{Halphen}}
\rput(1.2,7.3){\tiny{divisor $\Delta$}}
\rput(4.4,7.53){\tiny{Blow up}}

% second row
\rput(-0.2,6.5){$A_0^{(1)}$}
\rput(-4.2,6.5){\tiny{$\tau x_1 x_2 x_3 +x_1^3 +x_2^3+x_3^3$}}
\rput(-4.2,6.1){\tiny{$ 
+  a_1x_1^2 +b_1x_2^2 +c_1x_3^2 +  $}}
\rput(-4.2,5.7){\tiny{$ 
 + a_2x_1 + b_2x_2 +c_2x_3 + d $}}
 
 \rput(-1.8,6.3){\tiny{$\tau x_1 x_2 x_3 +$}}
 
 \rput(-1.7,5.9){\tiny{$ +x_1^3 +x_2^3+x_3^3$}}
 
\rput(-6.35,6.5){\tiny{Elliptic}}
\rput(-6.3,6.1) {{$\widetilde E_6$}}
%\psframe(3.1,5.6)(5.8,6.7)

\rput(-2,2){\psset{unit=.7}\pscurve(4.4,6)(4.5,6.3)(5.,6.5)(5.2,6.8)
\pscurve(4.4,6)(4.5,5.7)(5.,5.5)(5.2,5.2)}

\rput(4.5,6.5){N.A.}

% third row
\rput(-.2,5.){$A_0^{(1)^\star}$}
\rput(-4.2,5.3){\tiny{$\tau  x_1x_2x_3 + x_1^3 +x_2^3  +$}}
\rput(-4.2,4.8){\tiny{$ 
+  a_1x_1^2 +b_1x_2^2 +c_1x_3^2 +  $}}
\rput(-4.2,4.3){\tiny{$ 
 + a_2x_1 + b_2x_2 +c_2x_3 + d $}}
  \rput(-1.8,5.){\tiny{$\tau  x_1x_2x_3 +$}}
    \rput(-1.8,4.6){\tiny{$+ x_1^3 +x_2^3  $}}

\rput(2.3,4.75) {\psset{unit=0.4}
\begin{pspicture}(-3,-3)(7,3)
%right loop
\psbezier(0,0)(1,2)(3,2)(3,0)
\psline(0,0)(-1.,-1.5)
\psbezier(0,0)(1,-2)(3,-2)(3,0)
\psline(0,0)(-1.,1.5)
\end{pspicture}}

\rput(4.5,4.8){N.A.}

\rput(-6.45,5){\tiny{GDAHA}}
\rput(-6.5,4.5){$E_6^{(1)}$}

%fourth row
\rput(-.2,3.5){$A_1^{(1)}$}
\rput(-4.2,3.8){\tiny{$\tau  x_1 x_2 x_3 +x_1^3 $}}
\rput(-4.2,3.3){\tiny{$ 
+  a_1x_1^2 +b_1x_2^2 +c_1x_3^2 +  $}}
\rput(-4.2,2.8){\tiny{$ 
 + a_2x_1 + b_2x_2 +c_2x_3 + d $}}

  \rput(-1.7,3.3){\tiny{$\tau  x_1 x_2 x_3 +x_1^3 $}}

\rput(1.8,3.2){$\begin{array}{c}
\xy /r0.15pc/:
{\ar@/^12pt/@{-}(42,19);(42,2)},
{\ar@/^12pt/@{-}(48,2);(48,19)},
\endxy
\end{array}$}

\rput(4.5,3.4){N.A.}

%fifth row
\rput(-.2,2){$A_2^{(1)}$}
\rput(-4.2,2.){\tiny{$x_1 x_2 x_3$}}
\rput(-4.2,1.6){\tiny{$ 
+  a_1x_1^2 +b_1x_2^2 +c_1x_3^2 +  $}}
\rput(-4.2,1.2){\tiny{$ 
 + a_2x_1 + b_2x_2 +c_2x_3 + d $}}
 
   \rput(-1.8,1.6){\tiny{$x_1 x_2 x_3$}}

\rput(2.2,1.8){\psset{unit=0.2}
\begin{pspicture}(-5,-3.5)(7,3.5)
\pcline[linewidth=1pt](-3.5,1.8)(2.5,-2.5)
\pcline[linewidth=1pt](-4,-3)(2,-1)
\pcline[linewidth=1pt](-3,2.8)(-2,-3)
\end{pspicture}}

\rput(4.5,2){N.A.}

\rput(-6.5,2.7){\tiny{Deg.}}
\rput(-6.4,2.3){\tiny{GDAHA}}

%sixth row
\rput(-.2,0.5){$D_4^{(1)}$}
\rput(-4.2,.7){\tiny{$x_1 x_2 x_3 -x_1^2-x_2^2-x_3 ^2+$}}
\rput(-4.2,.2){\tiny{$
+\omega_1 x_1+\omega_2 x_2+\omega_3 x_3+\omega_4$}}
\rput(-6.5,0.6){\tiny{DAHA}}
\rput(-6.5,0){$\check C C_1$}

\psline(1.8,-.3)(1.8,.9)
\psline(1,-.3)(2.,.9)
\psline(1.6,.9)(2.6,-.3)
\rput(1.5,.3){\tiny{$\bullet$}}
\rput(1.25,0){\tiny{$\bullet$}}
\rput(1.8,.3){\tiny{$\bullet$}}
\rput(1.8,0){\tiny{$\bullet$}}
\rput(2.1,.3){\tiny{$\bullet$}}
\rput(2.35,0){\tiny{$\bullet$}}
\pscircle[linewidth=.5pt](1.8,.65){.08}

\psline(3.5,0.5)(5,0.5)
\psline(3.7,0.8)(3.7,-.3)
\psline(4,0.8)(4,-.3)
\psline(4.3,0.8)(4.3,-.3)
\psline(4.8,0.8)(4.8,-.3)

   \rput(-1.8,0.4){\tiny{$x_1 x_2 x_3$}}

%seventh row
\rput(-.2,-1){$D_5^{(1)}$}
\rput(-4.2,-1.){\tiny{$x_1 x_2 x_3 -x_1^2-x_2^2+$}}
\rput(-4.2,-1.4){\tiny{$\omega_1 x_1+\omega_2 x_2+\omega_3 x_3+\omega_4$}}
\psline(4.5,-.8)(5.4,-1.3)
\psline(3.8,-1.4)(4.8,-.8)
\psline(5.3,-1.8)(5.3,-.9)
\psline(5.,-1.8)(5.,-.9)
\psline(4.3,-1.8)(4.3,-.9)
\psline(4.,-1.8)(4.,-.9)

\rput(0,-1.5){\psline(1.8,-.3)(1.8,.9)
\psline(1,-.3)(2.,.9)
\psline(1.6,.9)(2.6,-.3)
\rput(1.5,.3){\tiny{$\bullet$}}
\rput(1.25,0){\tiny{$\bullet$}}
\rput(1.8,0){\tiny{$\bullet$}}
\rput(2.1,.3){\tiny{$\bullet$}}
\rput(2.35,0){\tiny{$\bullet$}}
\pscircle[linewidth=.5pt](1.8,.65){.08}}

   \rput(-1.8,-1.3){\tiny{$x_1 x_2 x_3$}}

%%% eight row

\rput(-.2,-2.5){$D_6^{(1)}$}
\rput(-4.2,-2.2){\tiny{$x_1 x_2 x_3 -x_1^2-x_2^2+$}}
\rput(-4.2,-2.8){\tiny{$\omega_1 x_1+\omega_2 x_2+\omega_4$}}
\psline(3.7,-2.2)(3.7,-3.1)
\psline(4.0,-2.2)(4,-3.1)
\psline(4.3,-2.2)(4.3,-3.1)
\psline(4.6,-2.2)(4.6,-3.1)
\psline(4.9,-2.2)(4.9,-3.1)

\psline(3.5,-2.8)(4.45,-2.8)
\psline(4.2,-2.35)(5.05,-2.35)

\rput(0,-3){\psline(1.8,-.3)(1.8,.9)
\psline(1,-.3)(2.,.9)
\psline(1.6,.9)(2.6,-.3)
\rput(1.5,.3){\tiny{$\bullet$}}
\rput(1.25,0){\tiny{$\bullet$}}
\rput(2.1,.3){\tiny{$\bullet$}}
\rput(2.35,0){\tiny{$\bullet$}}
\pscircle[linewidth=.5pt](1.8,.65){.08}}

   \rput(-1.8,-2.6){\tiny{$x_1 x_2 x_3$}}
   
% eight row

\rput(-6.5,-4.4){\tiny{Deg.}}
\rput(-6.4,-4.8){\tiny{DAHA}}

\rput(-.2,-4){$D_7^{(1)}$}
\rput(-4.2,-4){\tiny{ $x_1 x_2 x_3 -x_1^2-x_2^2+$}}
\rput(-4.2,-4.4){\tiny{$+\omega_1 x_1-x_2$}}

\psline(3.4,-4.5)(4.4,-4.5)
\psline(3.5,-4.7)(3.5,-3.7)
\psline(3.9,-4.7)(3.9,-3.7)

\psline(4.7,-4.5)(5.7,-4.5)
\psline(5.1,-4.7)(5.1,-3.7)
\psline(5.5,-4.7)(5.5,-3.7)

\psline(4.1,-4.7)(4.6,-3.7)
\psline(5,-4.6)(4.4,-3.7)

\rput(0,-4.5){
\psline(1,-.3)(2.,.9)
\psline(1.6,.9)(2.6,-.3)
\rput(1.5,.3){\tiny{$\bullet$}}
\rput(1.7,.3){\tiny{$2$}}
\pscircle[linewidth=.5pt](1.5,.3){.08}
\rput(1.25,0){\tiny{$\bullet$}}
\rput(1.4,0){\tiny{$2$}}
\pscircle[linewidth=.5pt](1.25,0){.08}
%\rput(2.1,.3){\tiny{$\bullet$}}
%\rput(2.35,0){\tiny{$\bullet$}}
\pscircle[linewidth=.5pt](1.8,.65){.08}}

   \rput(-1.8,-4){\tiny{$x_1 x_2 x_3$}}

\rput(-.2,-5.5){$D_8^{(1)}$}
\rput(-4.2,-6.){\tiny{
$x_1 x_2 x_3 -x_1^2-x_2^2-x_2$}}

\psline(3.5,-6.3)(3.5,-5.3)
\psline(3.9,-6.3)(3.9,-5.3)
\psline(5.3,-6.3)(5.3,-5.3)
\psline(5.6,-6.3)(5.6,-5.3)
\psline(3.4,-6.1)(4.4,-6.1)
\psline(4.1,-6.3)(4.6,-5.3)
\psline(5,-6.2)(4.4,-5.3)
\psline(4.7,-5.5)(5.7,-5.5)
\psline(4.5,-6.3)(5.,-5.3)

\rput(0,-6.){
\psline(1,-.3)(2.,.9)
\psline(1.6,.9)(2.6,-.3)
\rput(1.5,.3){\tiny{$\bullet$}}
\rput(1.66,.3){\tiny{$4$}}
\pscircle[linewidth=.5pt](1.5,.3){.08}
%\rput(2.1,.3){\tiny{$\bullet$}}
%\rput(2.35,0){\tiny{$\bullet$}}
\pscircle[linewidth=.5pt](1.8,.65){.08}}

   \rput(-1.8,-6){\tiny{$x_1 x_2 x_3$}}
%%%

\rput(-.2,-7){$E_6^{(1)}$}
\rput(-4.2,-7){\tiny{$x_1 x_2 x_3 -x_1^2+$}}
\rput(-4.2,-7.4){\tiny{$\omega_1 x_1+\omega_2 x_2+\omega_3 x_3+\omega_4$}}

\psline(3.3,-6.9)(3.85,-6.9)
\psline(4.1,-6.9)(4.75,-6.9)
\psline(5,-6.9)(5.65,-6.9)
\psline(3.5,-7.5)(5.5,-7.5)
\psline(3.6,-6.8)(3.6,-7.7)

\psline(4.3,-6.8)(4.3,-7.7)

\psline(5.2,-6.8)(5.2,-7.7)

\rput(0,-7.6){
\psline(1,-.3)(2.,.9)
\psline(1.6,.9)(2.6,-.3)
\rput(1.5,.3){\tiny{$\bullet$}}
%\pscircle[linewidth=.5pt](1.5,.3){.08}
\rput(1.25,0){\tiny{$\bullet$}}
%\pscircle[linewidth=.5pt](1.25,0){.08}
%\rput(2.1,.3){\tiny{$\bullet$}}
\rput(2.35,0){\tiny{$\bullet$}}
\pscircle[linewidth=.5pt](2.35,0){.08}
\rput(2.5,0){\tiny{$2$}}
\pscircle[linewidth=.5pt](1.8,.65){.08}}

%%%
\rput(-.2,-8.5){$E_7^{(1)}$}
\rput(-4.2,-8.5){\tiny{$ x_1 x_2 x_3 - x_1^2+$}}
\rput(-4.2,-8.9){\tiny{$ +\omega_1  x_1-x_2-1$}}
% linee orizzontali da sx a dx
\psline(3.2,-9.)(3.6,-9.)
\psline(3.5,-8.5)(4.3,-8.5)
\psline(4.1,-9)(4.9,-9)

% linee verticali da sx a dx
\psline(3.3,-8.4)(3.3,-9.1)
\psline(4.7,-8.4)(4.7,-9.1)
\psline(3.9,-8.4)(3.9,-9.1)
\psline(3.55,-8.4)(3.55,-9.1)
\psline(4.2,-8.4)(4.2,-9.1)

   \rput(-1.8,-7.2){\tiny{$x_1 x_2 x_3$}}
   \rput(-1.8,-8.5){\tiny{$x_1 x_2 x_3$}}
      \rput(-1.8,-10){\tiny{$x_1 x_2 x_3$}}

\rput(0,-9){
\psline(1,-.3)(2.,.9)
\psline(1.6,.9)(2.6,-.3)
\rput(1.5,.3){\tiny{$\bullet$}}
%\pscircle[linewidth=.5pt](1.5,.3){.08}
\rput(1.25,0){\tiny{$\bullet$}}
%\pscircle[linewidth=.5pt](1.25,0){.08}
%\rput(2.1,.3){\tiny{$\bullet$}}
%\rput(2.35,0){\tiny{$\bullet$}}
%\pscircle[linewidth=.5pt](2.35,0){.08}
\pscircle[linewidth=.5pt](1.8,.65){.08}}

%% last row 
\rput(-.2,-10){$E_8^{(1)}$}
\rput(-4.2,-10){\tiny{$x_1 x_2 x_3-x_1-x_2+1$}}
% linee orizzontali da sx a dx
\psline(3.2,-10.5)(3.6,-10.5)
\psline(3.5,-10)(4.3,-10)
\psline(4.1,-10.5)(4.9,-10.5)
\psline(4.6,-10)(5.2,-10)
% linee verticali da sx a dx
\psline(4.7,-9.9)(4.7,-10.6)
\psline(3.9,-9.9)(3.9,-10.6)
\psline(3.55,-9.9)(3.55,-10.6)
\psline(4.2,-9.9)(4.2,-10.6)
\psline(5.1,-9.9)(5.1,-10.6)

\rput(0,-10.5){
\psline(1,-.3)(2.,.9)
\rput(1.5,.3){\tiny{$\bullet$}}
\rput(1.7,.3){\tiny{$9$}}
\pscircle[linewidth=.5pt](1.5,.3){.08}}
%\rput(1.25,0){\tiny{$\bullet$}}
%\pscircle[linewidth=.5pt](1.25,0){.08}
\rput(0,-11.5){{\textsc{Table 4.} Results for del Pezzo of degree 3.}}
\end{pspicture}

%% QUANTUM TABLE

\begin{figure}[h]
%\psset{unit=0.8}
\begin{pspicture}(-6,-1)(6,8)
\psframe(-6,-.5)(6,8)

%vertical lines
\psline(-2.6, -.5)(-2.6,8)
\psline(1.2,-.5)(1.2,8)
\psline(4.5,-.5)(4.5,8)

%horizontal lines
\psline(-6,7)(6,7)
\psline(-6,5.5)(6,5.5)
\psline(-6,4)(6,4)
\psline(-6,2.5)(6,2.5)
\psline(-6,1)(6,1)

% first row
\rput(-4.3,7.5){\tiny{Polynomial  $\phi$}}
\rput(-1.1,7.5){\tiny{Quantum relations}}
\rput(2.6,7.5){\tiny{potential}}
\rput(5.2,7.7){\tiny{Central}}
\rput(5.2,7.5){\tiny{ element}}

% second row
\rput(-4.2,6.5){\tiny{$x_1 x_2 x_3 +x_1^3 +x_2^3+x_3^3$}}
\rput(-4.2,6.1){\tiny{$ 
+  a_1x_1^2 +b_1x_2^2 +c_1x_3^2 +  $}}
\rput(-4.2,5.7){\tiny{$ 
 + a_2x_1 + b_2x_2 +c_2x_3 + d $}}
 \rput(-1.4,6.5){\eqref{eq:t-comm}}
 \rput(2.6,6.5){{$\Phi_{EG}+\Psi_{EG}$}}
\rput(2.7,6){{see \eqref{homSklcub}, \eqref{homSklcubp}}}
\rput(5,6.5){\eqref{eq:omegaEG}}

% third row
\rput(-4.2,5.3){\tiny{$ x_1x_2x_3 + x_1^3 +x_2^3  +$}}
\rput(-4.2,4.8){\tiny{$ 
+  a_1x_1^2 +b_1x_2^2 +c_1x_3^2 +  $}}
\rput(-4.2,4.3){\tiny{$ 
 + a_2x_1 + b_2x_2 +c_2x_3 + d $}}

\rput(-1.4,5.){\eqref{eq:p,q,r-comm}}
\rput(-1.4,4.4){with $\gamma=0$}
\rput(2.6,5.){{$\Phi_{\alpha,\beta,0}+\Psi_{EG}$}}
\rput(2.7,4.4){{see \eqref{nonhomSklcub}, \eqref{homSklcubp}}}
\rput(5,4.8){\eqref{eq:cen-EG-gen}}

%fourth row

\rput(-4.2,3.8){\tiny{$x_1 x_2 x_3 +x_1^3 $}}
\rput(-4.2,3.3){\tiny{$ 
+  a_1x_1^2 +b_1x_2^2 +c_1x_3^2 +  $}}
\rput(-4.2,2.8){\tiny{$ 
 + a_2x_1 + b_2x_2 +c_2x_3 + d $}}
\rput(-1.4,3.6){\eqref{eq:p,q,r-comm}}
\rput(-1.4,3.){with $\beta=\gamma=0$}
\rput(2.6,3.6){{$\Phi_{\alpha,0,0}+\Psi_{EG}$}}
\rput(2.7,3.){{see \eqref{nonhomSklcub}, \eqref{homSklcubp}}}
\rput(5,3.6){\eqref{eq:cen-EG-gen}}
\rput(5,3.2){with}
\rput(5,2.8){$\beta=0$}
%fifth row

\rput(-4.2,2.){\tiny{$x_1 x_2 x_3$}}
\rput(-4.2,1.6){\tiny{$ 
+  a_1x_1^2 +b_1x_2^2 +c_1x_3^2 +  $}}
\rput(-4.2,1.2){\tiny{$ 
 + a_2x_1 + b_2x_2 +c_2x_3 + d $}}
 
   \rput(-1.4,2){\eqref{eq:t-comm}}
       \rput(-1.4,1.6){with $t=0$}
 \rput(2.6,2.2){$\Phi_{EG}+\Psi_{EG}$}
       \rput(2.6,1.8){with $t=0$}
\rput(2.7,1.4){{see \eqref{homSklcub}, \eqref{homSklcubp}}}

\rput(5,2){\eqref{eq:omega-eg0}}

%sixth row
\rput(-4.2,.7){$\phi_P^{(d)}$, \tiny{$d=PVI, \dots,PI$}}
\rput(-4.2,.2){see \eqref{eq:mon-mf}}
\rput(-.8,0.5){\eqref{eq:q-comm} with $\epsilon^{(d)}_i,\Omega_i^{(d)}$}
\rput(-1.2,0.0){ in  \eqref{eq:epsilon} \eqref{eq:omega}}

 \rput(2.6,0.7){$\Phi_{\mathcal{UZ}}$}
       \rput(2.6,0.2){see \eqref{qEGcub}}

\rput(5,0.5){\eqref{q-cubics}}

\rput(0,-1){{\textsc{Table 5.} Quantum counterpart of Table 4 \tiny{(we have squashed the last eight lines of Table 4 into one)}.}}
\end{pspicture}
\end{figure}

\newpage
%%%% TABLE of Mori and Smith

\begin{figure}[h]
\label{mori-smith}\rput(0,-6){
\begin{pspicture}(-6,-1)(6,8)
\psframe(-6,2.5)(6,8)

%vertical lines
\psline(-2.6, 2.5)(-2.6,8)
\psline(1.2,2.5)(1.2,8)
\psline(2.5,2.5)(2.5,8)

%horizontal lines
\psline(-6,7)(6,7)
\psline(-6,5.5)(6,5.5)
\psline(-6,4)(6,4)

% first row
\rput(-4.3,7.5){\tiny{Polynomial  $\phi$}}
\rput(-1.1,7.5){\tiny{Quantum relations}}
\rput(1.8,7.7){\tiny{Halphen}}
\rput(1.8,7.3){\tiny{surface}}
\rput(4.2,7.5){\tiny{Divisor $\Delta$}}

% first column

\rput(-4.3,6.5){$x_1^3-x_2^2 x_3$}
\rput(-4.3,4.8){$x_2^2 x_3-x_1^2 x_2$}
\rput(-4.3,3.3){$x_1^3+ x_3^3$}

% second column

\rput(-1.1,6.6){$x_1^2=x_2^2=0$ }
\rput(-1.,6.1){$x_3 x_2+x_2 x_3=0$}

\rput(-.8,5){$  x_2 x_3+x_3 x_2 -x_1^2=0$}
\rput(-.8,4.5){$ x_2^2=x_2 x_1+x_1 x_2=0$}

\rput(-.8,3.3){$x_1^2 =x_2^2=0$}

% third column 
\rput(2,6.5){$A_0^{(1)^{\ast\ast}}$}
\rput(2,4.8){$A_1^{(1)^{\ast}}$}
\rput(2,3.3){$A_2^{(1)^{\ast}}$}

% fourth column 
\psarc(4.2,6.7){0.5}{180}{270}
\psarc(4.2,5.7){0.5}{90}{180}

\psarc(4.2,5.2){0.5}{180}{360}
\psarc(4.2,4.2){0.5}{0}{180}

\psline(4.1,2.7)(4.1,3.7)
\psline(4.4,2.7)(3.8,3.7)
\psline(3.8,2.7)(4.4,3.7)

\rput(0,2){\textsc{Table 6.} Non Calabi-Yau cases}
\end{pspicture}}
\end{figure}

\newpage

\bibliography{qc} % e.g., if it's test.bib, just put NAMEOFBIBFILE=test -- do not include the .BIB extension.

\begin{thebibliography}{10}

\bibitem{ATV}
M.~Artin, J.~Tate, and M.~Van~den Bergh.
\newblock Modules over regular algebras of dimension {$3$}.
\newblock {\em Invent. Math.}, 106(2):335--388, 1991.

\bibitem{AS}
Michael Artin and William~F. Schelter.
\newblock Graded algebras of global dimension {$3$}.
\newblock {\em Adv. in Math.}, 66(2):171--216, 1987.

\bibitem{DBerLeigh}
David Berenstein, Vishnu Jejjala, and Robert~G. Leigh.
\newblock Marginal and relevant deformations of {$N=4$} field theories and
  non-commutative moduli spaces of vacua.
\newblock {\em Nuclear Phys. B}, 589(1-2):196--248, 2000.

\bibitem{BT}
Roland Berger and Rachel Taillefer.
\newblock Poincar\'{e}-{B}irkhoff-{W}itt deformations of {C}alabi-{Y}au
  algebras.
\newblock {\em J. Noncommut. Geom.}, 1(2):241--270, 2007.

\bibitem{BCR}
M.~Bertola, M.~Cafasso, and V.~Rubtsov.
\newblock Noncommutative {P}ainlev\'{e} equations and systems of {C}alogero
  type.
\newblock {\em Comm. Math. Phys.}, 363(2):503--530, 2018.

\bibitem{Bouss}
Pierrick Bousseau.
\newblock Quantum mirrors of log-calabi-yau surfaces and higher genus curve
  counting.
\newblock {\em arXiv:1808:07336 v1}, 2018.

\bibitem{BG}
Alexander Braverman and Dennis Gaitsgory.
\newblock Poincar\'{e}-{B}irkhoff-{W}itt theorem for quadratic algebras of
  {K}oszul type.
\newblock {\em J. Algebra}, 181(2):315--328, 1996.

\bibitem{BridgeSmith}
Tom Bridgeland and Ivan Smith.
\newblock Quadratic differentials as stability conditions.
\newblock {\em Publ. Math. Inst. Hautes \'{E}tudes Sci.}, 121:155--278, 2015.

\bibitem{ChM}
Leonid Chekhov and Marta Mazzocco.
\newblock Shear coordinate description of the quantized versal unfolding of a
  {$D_4$} singularity.
\newblock {\em J. Phys. A}, 43(44):442002, 13, 2010.

\bibitem{ChMR}
Leonid~O. Chekhov, Marta Mazzocco, and Vladimir~N. Rubtsov.
\newblock Painlev\'{e} monodromy manifolds, decorated character varieties, and
  cluster algebras.
\newblock {\em Int. Math. Res. Not. IMRN}, (24):7639--7691, 2017.

\bibitem{Cher}
Ivan Cherednik.
\newblock Double affine {H}ecke algebras, {K}nizhnik-{Z}amolodchikov equations,
  and {M}acdonald's operators.
\newblock {\em Internat. Math. Res. Notices}, (9):171--180, 1992.

\bibitem{Cher1}
Ivan Cherednik.
\newblock Whittaker limits of difference spherical functions.
\newblock {\em Int. Math. Res. Not. IMRN}, (20):3793--3842, 2009.

\bibitem{Den}
F~Denef.
\newblock Les houches lecture in constructing string vacua.
\newblock In {\em String Theory and the Real World. From Particle Physics to
  Astrophysics}, Proc. Summer School in Theor. Physics, pages 483--610. North
  Holland, 2008.

\bibitem{Dolg}
Vasiliy Dolgushev.
\newblock The {V}an den {B}ergh duality and the modular symmetry of a {P}oisson
  variety.
\newblock {\em Selecta Math. (N.S.)}, 14(2):199--228, 2009.

\bibitem{EtGinzb}
Pavel Etingof and Victor Ginzburg.
\newblock Noncommutative del {P}ezzo surfaces and {C}alabi-{Y}au algebras.
\newblock {\em J. Eur. Math. Soc. (JEMS)}, 12(6):1371--1416, 2010.

\bibitem{EOR}
Pavel Etingof, Alexei Oblomkov, and Eric Rains.
\newblock Generalized double affine {H}ecke algebras of rank 1 and quantized
  del {P}ezzo surfaces.
\newblock {\em Adv. Math.}, 212(2):749--796, 2007.

\bibitem{FN}
Hermann Flaschka and Alan~C. Newell.
\newblock Monodromy- and spectrum-preserving deformations. {I}.
\newblock {\em Comm. Math. Phys.}, 76(1):65--116, 1980.

\bibitem{GShV}
Michael Gekhtman, Michael Shapiro, and Alek Vainshtein.
\newblock {\em Cluster algebras and {P}oisson geometry}, volume 167 of {\em
  Mathematical Surveys and Monographs}.
\newblock American Mathematical Society, Providence, RI, 2010.

\bibitem{GrPlessRoan}
B.~R. Greene, M.~R. Plesser, and S.-S. Roan.
\newblock New constructions of mirror manifolds: probing moduli space far from
  {F}ermat points.
\newblock In {\em Essays on mirror manifolds}, pages 408--448. Int. Press, Hong
  Kong, 1992.

\bibitem{Gri}
Pierre-Paul Grivel.
\newblock Une histoire du th\'{e}or\`eme de {P}oincar\'{e}-{B}irkhoff-{W}itt.
\newblock {\em Expo. Math.}, 22(2):145--184, 2004.

\bibitem{GrossHK}
Mark Gross, Paul Hacking, and Sean Keel.
\newblock Mirror symmetry for log {C}alabi-{Y}au surfaces {I}.
\newblock {\em Publ. Math. Inst. Hautes \'{E}tudes Sci.}, 122:65--168, 2015.

\bibitem{IS}
Natalia Iyudu and Stanislav Shkarin.
\newblock Three dimensional {S}klyanin algebras and {G}r\"{o}bner bases.
\newblock {\em J. Algebra}, 470:379--419, 2017.

\bibitem{IS1}
Natalia Iyudu and Stanislav Shkarin.
\newblock Classification of quadratic and cubic pbw algebras on three
  generators.
\newblock 2018.

\bibitem{MJ1}
Michio Jimbo and Tetsuji Miwa.
\newblock Monodromy preserving deformation of linear ordinary differential
  equations with rational coefficients. {II}.
\newblock {\em Phys. D}, 2(3):407--448, 1981.

\bibitem{K1}
Tom~H. Koornwinder.
\newblock The relationship between {Z}hedanov's algebra {${\rm AW}(3)$} and the
  double affine {H}ecke algebra in the rank one case.
\newblock {\em SIGMA Symmetry Integrability Geom. Methods Appl.}, 3:Paper 063,
  15, 2007.

\bibitem{LB}
Lieven Le~Bruyn.
\newblock Conformal {${\rm sl}_2$} enveloping algebras.
\newblock {\em Comm. Algebra}, 23(4):1325--1362, 1995.

\bibitem{M2}
Marta Mazzocco.
\newblock Confluences of the {P}ainlev\'{e} equations, {C}herednik algebras and
  {$q$}-{A}skey scheme.
\newblock {\em Nonlinearity}, 29(9):2565--2608, 2016.

\bibitem{M4}
Marta Mazzocco.
\newblock Whittaker degenerations of gdaha.
\newblock {\em in progress}, 2019.

\bibitem{MolRag}
A.~I. Molev and E.~Ragoucy.
\newblock Symmetries and invariants of twisted quantum algebras and associated
  {P}oisson algebras.
\newblock {\em Rev. Math. Phys.}, 20(2):173--198, 2008.

\bibitem{MoSm}
Izuru Mori and S.~Paul Smith.
\newblock The classification of 3-{C}alabi-{Y}au algebras with 3 generators and
  3 quadratic relations.
\newblock {\em Math. Z.}, 287(1-2):215--241, 2017.

\bibitem{NS}
Masatoshi Noumi and Jasper~V. Stokman.
\newblock Askey-{W}ilson polynomials: an affine {H}ecke algebra approach.
\newblock In {\em Laredo {L}ectures on {O}rthogonal {P}olynomials and {S}pecial
  {F}unctions}, Adv. Theory Spec. Funct. Orthogonal Polynomials, pages
  111--144. Nova Sci. Publ., Hauppauge, NY, 2004.

\bibitem{Obl}
Alexei Oblomkov.
\newblock Double affine {H}ecke algebras of rank 1 and affine cubic surfaces.
\newblock {\em Int. Math. Res. Not.}, (18):877--912, 2004.

\bibitem{OdRat}
Aleksandr~V. Odesski.
\newblock Rational degeneration of elliptic quadratic algebras.
\newblock In {\em Infinite analysis, {P}art {A}, {B} ({K}yoto, 1991)},
  volume~16 of {\em Adv. Ser. Math. Phys.}, pages 773--779. World Sci. Publ.,
  River Edge, NJ, 1992.

\bibitem{Od84}
A.~V. Odesski\u{\i}.
\newblock An analogue of the {S}klyanin algebra.
\newblock {\em Funktsional. Anal. i Prilozhen.}, 20(2):78--79, 1986.

\bibitem{FO}
A.~V. Odesski\u{\i} and B.~L. Fe\u{\i}gin.
\newblock Sklyanin's elliptic algebras. {T}he case of a point of finite order.
\newblock {\em Funktsional. Anal. i Prilozhen.}, 29(2):9--21, 95, 1995.

\bibitem{OdRubpol}
A.~V. Odesski\u{\i} and V.~N. Rubtsov.
\newblock Polynomial {P}oisson algebras with a regular structure of symplectic
  leaves.
\newblock {\em Teoret. Mat. Fiz.}, 133(1):3--23, 2002.

\bibitem{O79}
Kazuo Okamoto.
\newblock Sur les feuilletages associ\'{e}s aux \'{e}quations du second ordre
  \`a points critiques fixes de {P}. {P}ainlev\'{e}.
\newblock {\em Japan. J. Math. (N.S.)}, 5(1):1--79, 1979.

\bibitem{POR}
Giovanni Ortenzi, Vladimir Rubtsov, and Serge~Rom\'{e}o Tagne~Pelap.
\newblock On the {H}eisenberg invariance and the elliptic {P}oisson tensors.
\newblock {\em Lett. Math. Phys.}, 96(1-3):263--284, 2011.

\bibitem{Pain}
P.~Painlev\'{e}.
\newblock M\'{e}moire sur les \'{e}quations diff\'{e}rentielles dont
  l'int\'{e}grale g\'{e}n\'{e}rale est uniforme.
\newblock {\em Bull. Soc. Math. France}, 28:201--261, 1900.

\bibitem{Pol}
A.~Polishchuk.
\newblock Poisson structures and birational morphisms associated with bundles
  on elliptic curves.
\newblock {\em Internat. Math. Res. Notices}, (13):683--703, 1998.

\bibitem{PP}
Alexander Polishchuk and Leonid Positselski.
\newblock {\em Quadratic algebras}, volume~37 of {\em University Lecture
  Series}.
\newblock American Mathematical Society, Providence, RI, 2005.

\bibitem{Prid}
Stewart~B. Priddy.
\newblock Koszul resolutions.
\newblock {\em Trans. Amer. Math. Soc.}, 152:39--60, 1970.

\bibitem{R}
E.~Rains.
\newblock Computer calculation.

\bibitem{R1}
E.~Rains.
\newblock Elliptic double affine hecke algebras.
\newblock {\em arXiv:1709.02989}, 2017.

\bibitem{RR}
Vladimir Retakh and Vladimir Rubtsov.
\newblock Noncommutative {T}oda chains, {H}ankel quasideterminants and the
  {P}ainlev\'{e} {II} equation.
\newblock {\em J. Phys. A}, 43(50):505204, 13, 2010.

\bibitem{RoanIsing}
Shi-shyr Roan.
\newblock Mirror symmetry of elliptic curves and {I}sing model.
\newblock {\em J. Geom. Phys.}, 20(2-3):273--296, 1996.

\bibitem{Sa}
Siddhartha Sahi.
\newblock Nonsymmetric {K}oornwinder polynomials and duality.
\newblock {\em Ann. of Math. (2)}, 150(1):267--282, 1999.

\bibitem{sakai}
Hidetaka Sakai.
\newblock Rational surfaces associated with affine root systems and geometry of
  the {P}ainlev\'{e} equations.
\newblock {\em Comm. Math. Phys.}, 220(1):165--229, 2001.

\bibitem{BSh}
Boris Shoikhet.
\newblock The {PBW} property for associative algebras as an integrability
  condition.
\newblock {\em Math. Res. Lett.}, 21(6):1407--1434, 2014.

\bibitem{Smith}
S.~Paul Smith.
\newblock ``{D}egenerate'' 3-dimensional {S}klyanin algebras are monomial
  algebras.
\newblock {\em J. Algebra}, 358:74--86, 2012.

\bibitem{SYZ}
Andrew Strominger, Shing-Tung Yau, and Eric Zaslow.
\newblock Mirror symmetry is {$T$}-duality.
\newblock {\em Nuclear Phys. B}, 479(1-2):243--259, 1996.

\bibitem{TerU}
Paul Terwilliger.
\newblock The universal {A}skey-{W}ilson algebra and {DAHA} of type
  {$(C^\vee_1,C_1)$}.
\newblock {\em SIGMA Symmetry Integrability Geom. Methods Appl.}, 9:Paper 047,
  40, 2013.

\bibitem{GinzbCY}
Ginzburg V.A.
\newblock Calabi-yau algebras.
\newblock 2006.

\bibitem{SvdP}
Marius van~der Put and Masa-Hiko Saito.
\newblock Moduli spaces for linear differential equations and the
  {P}ainlev\'{e} equations.
\newblock {\em Ann. Inst. Fourier (Grenoble)}, 59(7):2611--2667, 2009.

\bibitem{VZ}
Luc Vinet and Alexei Zhedanov.
\newblock Quasi-linear algebras and integrability (the {H}eisenberg picture).
\newblock {\em SIGMA Symmetry Integrability Geom. Methods Appl.}, 4:Paper 015,
  22, 2008.

\bibitem{ChW1}
Chelsea Walton.
\newblock Representation theory of three-dimensional {S}klyanin algebras.
\newblock {\em Nuclear Phys. B}, 860(1):167--185, 2012.

\end{thebibliography}

\bibliographystyle{plain}

\end{document}